\documentstyle[12pt]{article}
\begin{document}
\title{Quasi-invariant and pseudo-differentiable measures on
a non-Archimedean Banach space.I. Real-valued measures.}
\author{Sergey V. Ludkovsky.}
\date{09.05.2001}
\maketitle
\begin{abstract}
Quasi-invariant and pseudo-differentiable measures on a
Banach space $X$ over a non-Archimedean locally compact 
infinite field with a non-trivial valuation
are defined and constructed. Measures are considered with values
in $\bf R$. Theorems and criteria are formulated and
proved about
quasi-invariance and pseudo-differentiability of measures relative
to linear and non-linear operators on $X$. 
Characteristic functionals of measures are studied.
Moreover, the non-Archimedean
analogs of the Bochner-Kolmogorov and Minlos-Sazonov theorems
are investigated. Infinite products of measures also are considered.
Convergence of quasi-invariant and pseudo-differentiable
measures in the corresponding spaces of measures
is investigated.
\end{abstract}
\section{Introduction.}
There are few works about integration in a classical 
Banach space, that is over the field $\bf R$ of real numbers
or the field $\bf C$ of complex numbers
\cite{boui,chri,cons,dal,sko,vah}.
On the other hand, for a non-Archimedean Banach space $X$
(that is over a
non-Archimedean field) this theory is less developed.
An integration in $X$ is a very important part of the 
non-Archimedean analysis.
 The progress of quantum mechanics and different branches of
modern physics related, for example, with theories of elementary
particles lead to the necessity of developing integration theory
in a non-Archimedean Banach space \cite{ardrvo,djdr,ish,vla,yojan}.
It may also be useful for the development of non-Archimedean analysis.
Non-Archimedean functional analysis develops rapidly
in recent years  and has many principal differences
from the classical functional analysis
\cite{khrum,roo,sch1,sch2,schrn,vla}. 
Topological vector spaces over non-Archimedean fields
are totally disconnected, classes of smoothness
for functions and compact operators are defined
for them quite differently from that of the classical case,
also the notion of the orthogonality of vectors
has obtained quite another meaning.
In the non-Archimedean case analogs of the Radon-Nikodym
theorem and the Lebesgue theorem about convergence
are true under more rigorous and another conditions.
Especially strong differences are for measures 
with values in non-Archimedean fields, because classical
notions of $\sigma $-additivity and quasi-invariance have lost 
their meaning. 
\par On the other hand the development of the 
non-Archimedean functional analysis and its applications in
non-Archimedean quantum 
mechanics \cite{vla,khrum,yojan} leads to the necessity
of solving such problems. For example, problems related with
quantum mechanics on manifolds are related with diffeomorphism groups
and loop groups,
their representations and measures on them \cite{ish,lu3,lubp2}. 
In articles \cite{lu1,lu3,luseabm,lubp2} 
quasi-invariant measures on diffeomorphism and loop groups
and also on non-Archimedean manifolds were constructed.
These measures were used for the investigation of irreducible representations 
of topological groups \cite{lu3,luseabm,lubp}.
The theorems proved in this work enlarge classes of  
measures on such groups and manifolds,
this also enlarges classes of representations.
For example, theorems of the Minlos-Sazonov type
characterize measures with the help of characteristic functionals
and compact operators. In the non-Archimedean case
compact operators are more useful, than nuclear operators
in the classical case. Theorems of the Bochner-Kolmogorov and
Kakutani type characterize products of measures
and their absolute continuity relative to others measures.
\par In this article measures are considered on Banach spaces,
though the results given below can be developed
for more general topological vector spaces, for example,
it is possible to follow the ideas of works
\cite{macoma,maarma,mamaz}, in which were considered
non-Archimedean analogs of the Minlos-Sazonov theorems for real-valued 
measures on topological vector spaces over non-Archimedean fields
of zero characteristic. But it is impossible to make in one article.
In this article, apart from articles of M\c{a}drecki,
measures are considered also with values in non-Archimedean fields,
for the cases of real-valued measures also Banach spaces
over non-Archimedean fields $\bf K$ of characteristic
$char ({\bf K})>0$ are considered.
It is well-known, that a real-valued measure $m$ on a locally compact 
Hausdorff totally disconnected
Abelian topological group $G$ is called the Haar measure, if 
\par $(H)$ $m(x+A)=m(A)$ for each $x\in G$ and  
each Borel subset $A$ in $G$. \\
For the $s$-free group $G$ a measure $m$ 
with values in a non-Archimedean field $\bf K_s$ satisfy condition $(H)$
only for an algebra of clopen (closed and open) subsets $A$,
where a field $\bf K_s$ is a finite algebraic extension of $\bf Q_s$.
Indeed, in the last case if a measure is 
locally finite and $\sigma $-additive on
the Borel algebra of $G$, then it is purely atomic with atoms
being singletons, so it can not be invariant
relative to the entire Borel algebra (see Chapters 7-9 \cite{roo}). 
A lot of definitions and theorems given below are the non-Archimedean
analogs of classical results. Frequently their formulations and
proofs differ strongly. If proofs differ slightly, only general
circumstances are given in the non-Archimedean case.
\par In \S 2 weak distributions, 
characteristic functions of measures and their properties
are defined and investigated.
The non-Archimedean analogs of the Minlos-Sazonov and 
Bochner-Kolmogorov theorems are given. Quasi-measures 
also are considered. 
In \S 3 products of measures are considered together with their 
density functions. 
In the present paper broad classes of quasi-invariant measures 
are defined and constructed.
Theorems about quasi-invariance of measures under definite 
linear and non-linear transformations
$U: X\to X$ are proved. \S 4 contains a notion of pseudo-differentiability
of measures. This is necessary, because for functions
$f: {\bf K}\to \bf R$ there is not
any notion of differentiability (there is not such non-linear
non-trivial $f$). 
There are given criteria for the pseudo-differentiability. 
In \S 5 there are
given theorems about converegence of measures with
taking into account their quasi-invariance and pseudo-differentiability,
that is, in the corresponding spaces of measures.
The main results are Theorems 2.27, 2.35, 3.4,
3.20, 3.24, 3.25, 4.2, 4.3, 4.5, 4.7, 5.7-5.10.
\par In the first part of this article real-valued measures are considered.
In the second part measures with values in non-Archimedean fields are
considered. This is caused by differences in definitions, formulations of
statements and their proofs in two such cases. When differences are 
small only briefly matters are dicussed in the second case.
\par {\bf Notations.} Henceforth, $\bf K$
denotes a locally compact infinite field with a non-trivial norm,
then the Banach space $X$ is over $\bf K$. 
In the present article measures on $X$ have values in $\bf R$.
We assume that either $\bf K$ is a finite algebraic extension of 
the field of $p$-adic numbers $\bf Q_p$ or 
$char ({\bf K})=p$ and $\bf K$ is isomorphic with a field
$\bf F_p{(\theta )}$ of formal power series
consisting of elements $x=\sum_ja_j\theta ^j$,
where $a_j\in \bf F_p$, $|\theta |=p^{-1}$, 
$\bf F_p$ is a finite field of $p$ elements,
$p$ is a prime number. These imply that 
$\bf K$ has the Haar measures with values in $\bf R$ 
\cite{roo}. If $X$ is a Hausdorff topological space
with a small inductive dimension $ind (X)=0$, then
\par $E$ denotes an algebra of subsets of $X$, 
as a rule $E\supset Bf(X)$ for real-valued measures, where
\par $Bf(X)$ is a Borel $\sigma $-field of $X$ in \S 2.1;
\par $Af(X,\mu )$ is the completion of $E$ by a 
measure $\mu $ in \S 2.1;
\par ${\sf M}(X)$ is a space of norm-bounded measures on $X$ in \S 2.1;
\par ${\sf M}_t(X)$ is a space of Radon norm-bounded measures in \S 2.1;
\par $c_0(\alpha ,{\bf K})$ is a Banach space
and $P_L$ is a projector (fixed relative to a chosen basis) in \S 2.2;
\par $\mu _L$ is a projection of a measure $\mu $ in \S 2.2;
\par $\{ \mu _{L_n}: n \} $ is a sequence of weak distributions in \S 2.2;
\par $B(X,x,r)$ is a ball in \S 2.2;
\par $\chi _{\xi }$ is a character with values in ${\bf T}\subset \bf C$ 
in \S 2.6;
\par $\theta (z)=\hat \mu $ is a characteristic functional in \S 2.6;
\par $\delta _0$ is the Dirac measure in \S 2.8;
\par $\mu _1*\mu _2$ is a convolution of measures in \S 2.11;
\par $\psi _{q,\mu }$ and $\tau _q$ in \S 2.14;
\par $C(X,{\bf K})$ is a space of continuous functions from $X$ 
into $\bf K$ in \S 2.6;
\par $X^*$ is the topological dual space  of $X$ \cite{nari};
\par $\hat {\sf C}(Y,\Gamma )$, $\tau (Y)$ in \S 2.26;
\par $\sf B_{+}$, $\sf C_{+}$ in \S 2.29;
\par $\nu \ll \mu $, $\nu \sim \mu $, $\nu \perp \mu $ in \S 2.36.
\section{Weak distributions and families of measures.}
\par  {\bf 2.1.} For a Hausdorff topological space $X$ with
a small inductive dimension $ind (X)=0$ \cite{eng}
the Borel $\sigma$-field is denoted $Bf(X)$. 
Henceforth, measures $\mu $ are given on a measurable space $(X,E)$.
The completion of $E$
relative to $\mu $ is denoted by $Af(X,\mu )$. The total variation of $\mu $
with values in $\bf R$ on a subset $A$ is denoted by 
$\| \mu |_A \| $ or $| \mu |(A)$ for $A \in Af(X,\mu )$.
If $\mu $ is non-negative and $\mu (X)=1$, then it is called a probability
measure.
\par A measure $\mu $ on $E$ is called Radon,
if for each $\epsilon >0$ there exists a compact subset $C\subset X$ 
such that $ \| \mu |_{(X\setminus C)} \| < \epsilon $.
Henceforth, ${\sf M}(X)$ denotes a space of norm-bounded
measures, ${\sf M}_t(X)$ is its subspace of Radon
norm-bounded measures.
\par  {\bf 2.2.} Each Banach space $X$ over $\bf K$ 
in view of Theorems 5.13 and 5.16 \cite{roo}
is isomorphic with $c_0(\alpha , {\bf K}):=
\{ x: \mbox{ } x=(x_j:$ $ j \in \alpha ),$
$card(j:$ $|x_j|_{\bf K}>b )< \aleph _0 \mbox{ for each } b>0 \} $, where
$\alpha $ is a set considered as an ordinal due to the Kuratowski-Zorn lemma, 
$card(A)$ denotes the cardinality of $A$,
$\| x \| := \sup (|x_j|: \mbox{ } j \in \alpha )$. A dimension of $X$ over
$\bf K$ is by the definition $dim_{\bf K}X := card (\alpha )$.
For each closed $\bf K$-linear subspace $L$ in $X$ there exists an operator
of a projection $P_L: X\to L$. Moreover, an orthonormal in the
non-Archimedean sense basis in $L$ has a completion to an orthonormal basis
in $X$ such that $P_L$ can be defined in accordance with a chosen basis.
\par If $A \in Bf(L)$, then $P_L^{-1}(A)$ is called a cylindrical subset
in $X$ with a base $A$, $B^L:=P_L^{-1}(Bf(L))$,
$B_0:=\cup (B^L:$ $ L\subset X, L  \mbox{ is a Banach subspace ,}
dim_{\bf K}X< \aleph _0 )$.
The minimal $\sigma $-algebra $\sigma B_0$ generated by $B_0$
coincides with $Bf(X)$, if $dim_{\bf K}X\le \aleph _0$. Henceforward, it is
assumed that $\alpha \le \omega _0$, where $\omega _0$ is the initial ordinal
with the cardinality $\aleph _0:=card({\bf N})$. Then there exists
an increasing sequence of Banach subspaces 
$L_n \subset L_{n+1}\subset ...$ such that
$cl(\cup [L_n: n])=X$, $dim_{\bf K}L_n=\kappa _n $ for each $n$, where
$cl(A)=\bar A$ denotes a closure of $A$ in $X$ for $A\subset X$. 
We fix a family of projections $P^{L_m}_{L_n}: L_m\to L_n$
such that $P^{L_m}_{L_n}P^{L_n}_{L_k}=P^{L_m}_{L_k}$
for each $m\ge n\ge k$.
A projection
of the measure $\mu $ onto $L$ denoted by
$\mu _L(A):=\mu (P_L^{-1}(A))$ for each $A\in Bf(L)$ compose the consistent
family:
$$(1) \mbox{  } \mu_{L_n}(A)=\mu _{L_m}(P_{L_n}^{-1}(A)\cap L_m)$$
for each 
$m \ge n$, since there are projectors $P_{L_n}^{L_m}$,
where $\kappa _n \le \aleph _0$ and there may be chosen
$\kappa _n< \aleph _0$ for each $n$.
\par An arbitrary family of measures $ \{ \mu _{L_n}: n \in {\bf N} \} $ 
having property $(1)$
is called a sequence of a weak distribution
(see also \cite{dal,sko}).
\par By $B(X,x,r)$ we denote a ball $\{y: \mbox{ }y \in X,
\mbox{ } \|x-y\| \le r \}$, which is clopen (closed and open) in $X$.
\par  {\bf 2.3. Lemma.} {\it A sequence of a weak distribution
$\{ \mu _{L_n}: n\}$
is generated by some measure $\mu $ on $Bf(X)$ if and only if
for each $c>0$ there exists $b>0$ such that $| |\mu _{L_n} |
(B(X,0,r) \cap L_n) - | \mu _{L_n} |( L_n)| \le c $ and $\sup_n
|\mu _{L_n} |(L_n) < \infty $.}
\par  {\bf Proof.} In the case of $\mu $ with values in $\bf R$
we can use a Hahn decomposition $\mu =\mu ^+ - \mu ^- $
and substitute everywhere in the proof of Lemma
1 \S 2\cite{sko} a Hilbert space over $\bf R$ onto the Banach space
$X$ over $\bf K$,
since $X$ is the Radon space in view of Theorem 1.2  \S
I.1.3 \cite{dal}, then $|\mu |(A)=\mu ^+(A) +\mu ^- (A)$ for $A \in
Bf(X)$.
\par {\bf 2.4. Definition and notations.} A 
function $\phi :X\to \bf R$
of the form $\phi (x)=\phi _S(P_S x)$ is called a cylindrical function if
$\phi _S$ is a $Bf(S)$-measurable
function on a finite-dimensional over
$\bf K$ space $S$ in $X$. For $\phi _S \in L^1 (S, \mu , {\bf R})$
for $\mu $ with values in $\bf R$ we can define an integral 
by a sequence of a weak distribution $\{ \mu _{S_n}: n \} $: 
$$\int_X \phi (x) \mu _* (dx)
:=\int \phi _{S_n}(x)\mu _{S_n} (dx) .$$ 
\par {\bf 2.5. Lemma.} {\it A subset $A\subset X=c_0(\omega _0,{\bf K})$
is relatively compact if and only if $A$ is bounded
and for each $c>0$ there exists a finite-dimensional over
$\bf K$ subspace $L\subset X$ such that $\bar A \subset L^c
:= \{ y \in X: \mbox{ } d(y,L):=\inf \{ \| x-y\| : x \in L\} \le c \} $.}
\par {\bf Proof.} If $A$ is bounded and for each $c>0$
there exists $L^c$ with $\bar A \subset L^c$, then there is a
sequence $ \{ k(j): \mbox{ } j \in {\bf N} \} \subset \bf Z$
such that $\lim _{j \to \infty } k(j)= \infty ,$  $\bar A \subset \{
x \in X: \mbox{ } |x_j| \le p^{-k(j)}, \mbox{ } j=1,2,...\}=:S$,
but $X$ is Lindel\"of, $S$ is sequentially compact,
hence $\bar A$ is compact (see \S 3.10.31 \cite{eng}). If
$\bar A$ is compact, then for each $c>0$ there exists a finite number
$m$ such that $\bar A \subset \bigcup_{j=1}^m B(X,x_j,c)$, where
$x_j \in X$. Therefore, $\bar A \subset L^c$ for $L=sp_{\bf K}(x_j:$
$j=1,...,m):=(x=\sum_{j=1}^m b_jx_j: b_j \in K).$
\par {\bf 2.6. Remarks and definitions.} 
When $char ({\bf K})=0$, then $\bf K$ as the $\bf Q_p$ linear space
is isomorphic with $\bf Q_p^n$, where $n \in {\bf N}:= \{ 1,2,... \} $. 
The topologically adjoint space
over $\bf Q_p$
(that is, of continuous linear functionals $f: {\bf K}\to \bf Q_p$)
is isomorphic with $\bf Q_p^n$ \cite{hew}.
For $x$ and $z \in \bf Q_p^n$ we denote by
$(z,x)$ the following sum $\sum_{j=1}^n x_jz_j$, where $x=(x_j:$
$j=1,...,n)$, $x_j \in \bf Q_p$. Each number $y\in \bf Q_p$ has  a
decomposition $y=\sum_l
a_lp^l$, where $\min (l:$ $a_l\ne 0)=:ord_p(y)> - \infty $
($ord(0):=\infty $) \cite{nari}, $a_l \in (0,1,...,p-1)$,
we define a symbol $\{ y \}_p:=\sum_{l<0} a_lp^l$ for
$|y|_p>1$ and $\{ y\} _p=0$ for $|y|_p \le 1$. 
\par For a locally compact field $\bf K$ with
a characteristic $char ({\bf K})=p>0$ let $\pi _j(x):=a_j$
for each $x=\sum_ja_j\theta ^j\in \bf K$ (see Notation).
For $\xi \in \bf K^*$ we denote $\xi (x)$ also by $(\xi ,x)$.
All continuous characters $\chi : {\bf K}\to \bf C$
have the form 
$$(1)\mbox{ }\chi _{\xi }(x)=\epsilon ^{z^{-1}\eta ((\xi ,x))}$$ 
for each $ \eta ((\xi ,x)) \ne 0$,
$\chi _{\xi }(x):=1$ for $ \eta ((\xi ,x))=0,$
where $\epsilon =1^z$ is a root of unity, $z=p^{ord(\eta ((\xi ,x)))},$ 
$\pi _j: {\bf K}\to \bf R$,
$\eta (x):= \{ x \} _p$ and $\xi \in {\bf Q_p^n}^*=
\bf Q_p^n$ for $char ({\bf K})=0$,
$\eta (x):=\pi _{-1}(x)/p$ and $\xi \in {\bf K}^*=\bf K$
for $char({\bf K})=p>0$, $x \in \bf K$, 
(see \S 25 \cite{hew}). 
Each $\chi $ is locally constant, hence $\chi : {\bf K}\to
\bf T$ is also continuous, 
where $\bf T$ denotes
the discrete group of all roots of $1$ (by multiplication).
\par For a measure $\mu $
there exists a characteristic functional (that is, called the
Fourier-Stieltjes transformation) $\theta =\theta _{\mu }: C(X,{\bf K}) 
\to \bf C$:
$$(2) \mbox{  } \theta (f):=\int _X \chi _e(f(x)) \mu (dx),$$
where either $e=(1,...,1)\in \bf Q_p^n$ for $char ({\bf K})=0$, 
or $e=1\in \bf
K^*$ for $char ({\bf K})=p>0$, $x \in X$, $f$ is in the space $C(X,{\bf K})$
of continuous functions from $X$ into $\bf K$, in particular for
$z=f$ in the topologically conjugated space $X^*$
over $\bf K$, $z:X \to \bf K$, $z \in X^*$, $\theta (z)=:
\hat \mu (z)$. 
It has the folowing properties:
$$(3a) \mbox{ }\theta (0)=1 \mbox{ for } \mu (X)=1$$ and $\theta (f)$
is bounded on $C(X,{\bf K})$;
$$(3b)\mbox{ }\sup_f | \theta (f)|=1\mbox{ for probability measures };$$ 
$$(4) \mbox{ } \theta (z) \mbox{ is weakly continuous, that is, } (X^*,
\sigma (X^*,X))\mbox{-continuous},$$
$\sigma (X^*,X)$ denotes a weak topology on $X^*$,
induced by the Banach space $X$ over $\bf K$. 
To each $x \in X$ there corresponds a continuous
linear functional $x^*:X^* \to \bf K$, $x^*(z):=z(x)$, moreover, $\theta
(f) $ is uniformly continuous relative to the norm on 
$$C_b(X,{\bf K}):=\{ f\in C(X,{\bf K}):
\| f\| :=\sup_{x\in X}|f(x)|_{\bf K}<\infty \} ;$$
$$(5) \mbox{ } \theta (z) \mbox{ is positive definite on }X^* \mbox{ and on }
C(X,{\bf K})$$
for $ \mu
\mbox{ with values in } [0, \infty ).$
\par  Property (4) follows from Lemma 2.3, boundedness and continuity of
$\chi _e$ and the fact that due to the Hahn-Banach theorem there is
$x_z \in X$
with $z(x_z)=1$ for $z \ne 0$ such that $z|_{(X \ominus L)}=0$ and
$$\theta (z)=\int_X \chi _e(P_L(x)) \mu(dx)=\int_L \chi _e(y) \mu _L(dy),$$
where $L=Kx_z$, also due to the Lebesgue theorem 2.4.9 \cite{fed}
for real measures. Indeed, for each $c>0$ there exists a
compact subset $S \subset X$ such that $|\mu |(X\setminus S)<c$,
each bounded subset $A \subset X^*$ is uniformly equicontinuous
on $S$ (see (9.5.4) and Exer. 9.202 \cite{nari}), that is,
$ \{ \chi _e(z(x)): $ $z \in A \} $ is the uniformly equicontinuous
family (by $x \in S$). On the other hand, $\chi _e(f(x))$
is uniformly equicontinuous on a bounded $A\subset C_b(X,{\bf K})$ by
$x\in S$.
\par Property (5) is accomplished, since 
$$\sum_{l,j=1}^N \theta (f_l
-f_j) \alpha _l \bar \alpha _j =\int_X |\sum_{j=1}^N \alpha _j
\chi _e(f_j(x))|^2 \mu (dx) \ge 0,$$ 
particularly, for $f_j=z_j\in X$, where $\bar \alpha _j$ is a complex 
conjugated number to $\alpha _j.$
\par   We call a functional $\theta $ finite-dimensionally concentrated,
if there exists $L \subset X$, $dim_{\bf K}L< \aleph _0 ,$ such that
$\theta |_{(X\setminus L)}=\mu (X)$. For each $c>0$ and $\delta >0$
in view of Theorem I.1.2 \cite{dal} and Lemma
2.5 there exists a finite-dimensional over $\bf K$ subspace $L$ and compact
$S \subset L^{\delta }$ such that 
$\| X\setminus S\| _{\mu } <c$. Let
$\theta ^L (z):=\theta (P_Lz)$. 
\par This definition is correct, since
$L\subset X$, $X$ has the isometrical embedding into $X^*$ as
the normed space associated with the fixed basis of $X$, 
such that functionals $z \in X$ separate points in $X$.
If $z \in L$, then $|\theta (z)-\theta ^L(z)| \le c\times b \times q$,
where $b=\| X\| _{\mu }$, $q$ is independent of $c$ and $b$.
Each characteristic functional $\theta ^L(z)$ is uniformly continuous by
$z \in L$ relative to the norm $\| *\|$ on $L$, since $|\theta ^L(z)
-\theta ^L(y)| $ $\le |\int_{S' \cap L} [\chi _e(z(x))-\chi _e(y(x))]$
$\mu _L(dx)|$ $+|\int_{L\setminus S'} [\chi _e(z(x))-\chi _e(y(x))]$
$\mu _L(dx)|$, where the second term does not exceed $2C'$ for 
$ \| L\setminus S'\| _{\mu _L}<c'$ for
a suitable compact subset $S' \subset X$ and $\chi _e(z(x))$ is 
an uniformly equicontinuous
by $x\in S'$ family relative to $z\in B(L,0,1)$.
\par Therefore,
$$(6) \mbox{  } \theta (z)=\lim _{n \to \infty } \theta _n(z)$$
for each finite-dimensional over $\bf K$ subspace $L$, where $\theta
_n(z)$ is uniformly equicontinuous and finite-dimensionally concentrated
on $L_n \subset X$, $z \in X$, $cl(\bigcup_nL_n)=X$, $L_n
\subset L_{n+1}$
for every $n$, for each $c>0$ there are $n$ and $q>0$
such that $|\theta (z)-\theta _j(z)| \le cbq$
for $z \in L_j$ and $j>n$, $q=const >0$ is independent of
$j$, $c$ and
$b$. Let $ \{ e_j:$ $j \in {\bf N} \} $ be the standard orthonormal 
basis in $X$,
$e_j=(0,...,0,1,0,...)$ with $1$ in $j$-th place. Using
countable additivity of $\mu $, local constantness of $\chi _e$,
considering all $z=be_j$ and $b \in \bf K$, we get that $\theta (z)$ on $X$
is non-trivial, whilst $\mu $ is a non-zero measure, since due to Lemma
2.3 $\mu $ is characterized uniquely by $ \{ \mu _{L_n} \} $. Indeed,
for $\mu $ with values in $\bf R$ a measure $\mu _V$ on $V$,
$dim_{\bf K}V< \aleph _0 ,$ this follows from
the properties of the Fourier transformation $F$ on spaces of 
generalized functions and also on $L^2(V,\mu _V,{\bf C})$ (see \S 7
\cite{vla}), where 
$$F(g)(z):=
\lim _{r \to \infty } \int_{B(V,0,r)} \chi _e(z(x))g(x)m(dx),$$
$z\in V,$ $m$ is the Haar measure on $V$ with values in $\bf R$.
Therefore, the mapping $\mu \mapsto \theta _{\mu }$ is injective.
\par  {\bf 2.7. Proposition.} {\it Let $X=\bf K^j$, $j \in \bf N$,
\par (a) $\mu $ and $\nu $ be real probability measures on $X$,
suppose $\nu $ is symmetric. Then $\int_X \hat \mu (x) \nu (dx)$ $=
\int_X \hat \nu(x) \mu (dx)$ $\in \bf R$ and for each
$0<l<1$ is accomplished the following inequality: \\ 
$\mu ([x \in X:$ 
$\hat \nu (x) \le l]) $ $\le \int_X(1-\hat \mu (x)) \nu (dx)/
(1-l)$. 
\par (b). For each real probability measure $\mu $
on $X$ there exists $r> p^3$ such that for each $R>r$ and $t>0$
the following inequality is accomplished: \\
$\mu ([x \in X:$ $\| x \| \ge tR]) \le $ $
c\int_X[1-\hat \mu (y \xi)] \nu (dy),$ \\
where $\nu (dx)=C\times
exp(-|x|^2)m(dx)$, $m$ is the Haar measure 
on $X$ with values in $[0, \infty )$,
$m(B(X,0,1))=1,$ $\nu (X)=1,$ $2>c=const \ge 1$ is independent from
$t,$ $c=c(r)$ is non-increasing whilst $r$ is increasing, $C>0$.}
\par {\bf Proof.} (a). Recall that $\nu $ is symmetric, if
$\nu (B)=\nu (-B)$ for each $B \in Bf(X)$.
Therefore, $\int_X \chi _e(z(x))\nu (dx)$ $=\int_X \chi _e(-z(x))
\nu (dx)$, that is equivalent to $\int _X [\chi _e(z(x))-\chi _e(-z(x))]
\nu (dx)=0$
or $\hat \nu (z) \in \bf R$. If $0<l<1$, then $\mu ([x \in X:$ $
\hat \nu (x) \le l])=\mu ([x:$ $1-\hat \nu (x) \ge 1-l])$
$\le \int_X(1-\hat \nu (x)) \mu (dx)/(1-l)=$ $\int_X(1-\hat \mu (x))
\nu (dx)/(1-l)$ due to the Fubini theorem.
\par  (b). Let $\nu (dx)=\gamma (x)m(dx),$ where $\gamma (x)=C\times exp(
-|x|^2),$ $C>0,$ $\nu (X)=1.$
Then $F(\gamma )(z)=:\hat \gamma (z) \ge 0,$ and $\hat \gamma (0)=1$
and $\gamma $ is the continuous positive definite function with
$\gamma
(z)\to 0$ whilst $|z| \to \infty .$ In view of (a):
$\mu ([x: \| x\| \ge tR])$ 
$\le \int_X [1-\hat \mu (y \xi )] \nu (dy)/(1-l),$ 
where $|\xi |=1/t,$ $t>0$, $l=l(R)$.
Estimating integrals, we get (b).
\par {\bf 2.8. Lemma.}{\it Let in the notation of Proposition 2.7
$\nu _{\xi }(dx)=\gamma _{\xi }(x)m(dx)$, $\gamma _{\xi }(x)=
C(\xi )exp(-|x \xi |^2),$ $\nu _{\xi }(X)=1,$ $\xi \ne 0,$ then
a measure $\nu _{\xi }$ is weakly converging to the Dirac measure
$\delta _0$ with the support in $0 \in X$ for $|\xi | \to \infty .$}
\par  {\bf Proof.} We have: $C(\xi )^{-1}=
C_q(\xi )^{-1}=\sum _{l \in \bf Z}[p^{lq}-p^{(l-1)q}]
exp(-p^{2l}|\xi |^2) < \infty ,$ where the sum by $l<0$ does not
exceed $1$, $q=jn$, $j=dim_{\bf K}X,$ $n=dim_{\bf Q_p}{\bf K}.$
Here $\bf K$ is considered as the Banach space 
$\bf Q_p^n$ with the following norm
$|*|_p$ equivalent to $|*|_{\bf K}$,
for $x=(x_1,...,x_j) \in X$ with $x_l \in \bf K$ as usually
$|x|_p=\max_{1\le l \le j}
|x_l|_p$, for $y=(y_1,...,y_n) \in \bf K$  with $y_l \in \bf Q_p$:
$|y|_p:=\max_{1\le l \le n} |y_l|_{\bf Q_p}$.
Further, $p^{l+s}
\sum_{x_l \ne 0} exp(2\pi i\sum_{i=l} ^{-s-1}x_i p^{i+s})$ $=
\int_{p^{l+s}}^1 exp(2 \pi i \phi )d\phi$ $+\beta (s),$ where
$s+l<0$,
$ \lim_{s \to -\infty }(\beta (s)p^{-s-l})=0,$ therefore, $sup[
|\hat \gamma _1(z)|_{\bf R} |z|_X:$ $ z \in X, $ $|z| \ge p^3] \le 2.$
Then taking $0\ne \xi \in \bf K$ and carrying out the substitution of variable
for continuous and bounded functions $f:X \to \bf R$ we get
$\lim _{|\xi | \to \infty }$
$\int_Xf(x) \nu _{\xi }(dx)=f(0).$ This means that $\nu _{\xi }$ is weakly
converging to $ \delta _0$ for $|\xi | \to \infty $.
\par {\bf 2.9. Theorem.}{\it Let $\mu _1$ and $\mu _2$ be 
measures in ${\sf M}(X)$ such that $\hat \mu _1(f)=\hat \mu _2(f)$
for each $f \in \Gamma $. Then $\mu _1 =\mu _2$, where $X=c_0(\alpha ,K),$
$\alpha \le \omega _0$, $\Gamma $ is a vector subspace in a space
of continuous functions $f: X\to \bf K$ separating points in $X$.}
\par {\bf Proof.} Let at first $\alpha < \omega _0$, then due to
continuity of the convolution $\gamma _{\xi }*\mu _j$ by $\xi $,
and Proposition 4.5 \S I.4\cite{vah} and Lemma 2.8 we get
$\mu _1=\mu _2,$ since the family $\Gamma $ generates $Bf(X)$.
Now let $\alpha =\omega _0$, $A=\{x \in
X:\mbox{ } (f_1(x),...,f_n(x)) \in S \},$ $\nu _j$ be an image of a measure
$\mu _j$ for a mapping $x \mapsto (f_1(x),...,f_n(x))$, where
$S \in Bf({\bf K^n})$, $ f_j \in
X \hookrightarrow X^*$. Then $\hat \nu _1(y)=\hat \mu _1(y_1f_1+
...+y_nf_n)=\hat \mu _2(y_1f_1+...+y_nf_n)=\hat \nu _2(y)$ for each
$y=(y_1,...,y_n) \in \bf K^n$, consequently, $\nu  _1=\nu _2$ on $E$.
Further we can use the Prohorov theorem 3.4 \S 1.3 \cite{vah},
since compositions of $f\in \Gamma $ with continuous functions
$g: {\bf K}\to \bf R$ generate a family of real-valued 
functions separating points of $X$.
\par  {\bf 2.10. Proposition.} {\it Let $\mu _l$ and $\mu $
be measures in ${\sf M}(X_l)$ and ${\sf M}(X)$ respectively,
where $X_l=
c_0(\alpha _l,{\bf K}),$ $\alpha _l \le \omega _0,$ $X=\prod_1^n X_l,$
$n \in \bf N.$ Then the condition $\hat \mu (z_1,...,z_n)=\prod_{l=1}^n
\hat \mu _l(z_l)$ for each $(z_1,...,z_n)\in X\hookrightarrow X^*$
is equivalent to $\mu =\prod_{l=1}^n \mu _l$.}
\par {\bf Proof.} Let $\mu =\prod_{l=1}^n \mu _l$, then
$\hat \mu (z_1,...,z_n)=\int _X \chi _e(\sum z_l(x_l))\prod_{l=1}^n
\mu _l (dx_l)$ $=\prod_{l=1}^n \int_{X_l} \chi _e(z_l(x_l))\mu _l(dx_l).$
The reverse statement follows from Theorem 2.9.
\par {\bf 2.11. Proposition.} {\it Let $X$ be a Banach space over $\bf K$;
suppose $\mu $,
$\mu _1$ and $\mu _2$ are probability measures on $X$. Then
the following conditions are equivalent: $\mu $ is the convolution of two 
measures $\mu _j$, $\mu=\mu _1*\mu _2$, and
$\hat \mu (z)=\hat \mu _1
(z)\hat \mu _2(z)$ for each $z \in X$.}
\par {\bf Proof.} Let $\mu =\mu _1 *\mu _2$. This means by the definition
that $\mu $ is the image of the measure $\mu _1 \otimes \mu _2$ for
the mapping $(x_1,x_2) \mapsto x_1 +x_2,$ $x_j \in X,$ consequently,
$\hat \mu (z)=\int_{X\times X} \chi _e(z(x_1+x_2))$ $(\mu _1 \otimes \mu _2)
(d(x_1,x_2))$ $=\prod_{l=1}^2 
\int_X \chi _e(z(x_l))\mu _l(dx_l)$ $=\hat \mu _1
(z)\hat \mu _2(z).$ On the other hand, if $\hat \mu _1\hat \mu _2=\mu ,$
then $\hat \mu =(\mu _1 *\mu _2)^{\wedge }$ and due to Theorem 2.9 above
for real measures we have $\mu =\mu _1*\mu _2.$
\par {\bf 2.12. Corollary.} {\it Let $\nu $ be a probability measure on
$Bf(X)$ and $\mu *\nu =\mu $ for each $\mu $,
then $\nu =\delta _0$.}
\par {\bf Proof.} If $z_0 \in X \hookrightarrow X^*$
and $\hat \mu (z_0) \ne 0,$  then from $\hat \mu (z_0)\hat \nu (z_0)=
\hat \mu (z_0)$ it follows that $\hat \nu _0(z_0)=1$. From the property
2.6(6) we get that there exists $m \in \bf N$ with $\hat \mu (z)
\ne 0$ for each $z$ with $\| z \| =p^{-m}$, since $\hat \mu (0)=1$.
Then $\hat \nu (z + z_0)=1$, that is, $\hat \nu |_{(B(X,
z_0,p^{-m}))}=1.$ Since $\mu $ are arbitrary we get
$\hat \nu |_X=1$, that is, $\nu =\delta _0$ due to \S 
2.6 and \S 2.9.
\par  {\bf 2.13. Corollary.} {\it Let $X$ and $Y$ be Banach spaces
over $\bf K$,
(a) $\mu $ and $\nu $ be probability measures on $X$ and $Y$
respectively, suppose $T: X\to Y$ is a continuous linear operator.
A measure $\nu $ is an image of $\mu $ for $T$ if and only if
$\hat \nu =\hat \mu \circ T^*,$ where $T^*: Y^* \to X^*$ 
is an adjoint operator. (b). A characteristic functional
of a real measure $\mu $ on $Bf(X)$ is real if and only if $\mu $ is
symmetric.}
\par {\bf Proof} follows from \S 2.6 and \S 2.9.
\par {\bf 2.14. Definition.} We say that a real probability measure
$\mu $ on $Bf(X)$ for a Banach space 
$X$ over $\bf K$ and $0<q< \infty $ has a weak
$q$-th order if $\psi _{q,\mu }(z)=\int_X|z(x)|^q\mu (dx) < \infty $
for each $z \in X^*$. The weakest vector topology in $X^*$ relative to which
all $(\psi _{q, \mu }: \mu)$ are continuous is denoted by $\tau _q$.
\par {\bf 2.15. Theorem.} {\it A characteristic functional $\hat \mu $
of a real probability Radon measure $\mu $ on $Bf(X)$ is continuous in
the topology $\tau _q$ for each $q>0.$}
\par  {\bf Proof.} For each $c>0$ there exists a compact $S
\subset X$ such that $\mu (S)>1-c/4$ and 
$$|1-\hat \mu (z)| \le |\int_S
(1-\chi _e(z(x)))\mu (dx)| + |\int_{X\setminus S} (1-\chi _e(z(x)))
\mu (dx)| \le |1-\hat \mu _c(z)|+c/2,$$ 
where $\mu _c(A)= ( \mu (A \cap S)/
\mu (S)$ and $A \in Bf(X)$; further analogously to the proof of
IV.2.3\cite{vah}.
\par  {\bf 2.16. Proposition.} {\it For a completely regular space
$X$ with $ind(X)=0$ the following statements are accomplished:
\par (a) if $(\mu _{\beta })$ is a bounded net of measures 
in ${\sf M}(X)$ that weakly converges to a measure
$\mu $ in ${\sf M}(X)$, then $(
\hat \mu _{\beta } (f))$ converges to $\hat \mu (f)$ for each continuous
$f: X \to \bf K$; if $X$ is separable and metrizable then $(\hat \mu _{\beta })$
converges to $\hat \mu $ uniformly on subsets that are
uniformly equicontinuous in $C(X,{\bf K})$; 
\par (b) if $M$ is a bounded dense family in a ball of the 
space ${\sf M}(X)$ for 
measures in ${\sf M}(X)$, then a family $(\hat \mu: $ $\mu \in M)$
is equicontinuous on a locally $\bf K$-convex space $C(X,{\bf K})$ in a topology
of uniform convergence on compact subsets $S \subset X$.}
\par {\bf Proof.} (a). Functions $\chi _e(f(x))$ 
are continuous and
bounded on $X$, where $\hat \mu (f)=\int_X \chi _e(f(x))\mu (dx)$. Then
(a) follows from the definition of the weak convergence and Proposition
1.3.9\cite{vah}, since $sp_{\bf C} \{ \chi _e(f(x)):$
$f\in C(X,{\bf K}) \} $ is dense in $C(X,{\bf C})$.
\par (b). For each $c>0$ there exists a compact subset 
$S \subset X$ such that
$| \mu |(S)>|\mu (X)|-c/4$.
Therefore, for $\mu \in \sf M$ and $f \in C(X,{\bf K})$ with
$|f(x)|_{\bf K}<c<1$ for $x \in S$ we get 
$$| \mu (X)-Re(\hat \mu (f)|=2
|\int_X [(\chi _e^{1/2}(f(x))-\chi _e^{1/2}(-f(x))/(2i)]^2\mu (dx)|<c/2.$$
Further analogously to
the proof of Proposition IV.3.1\cite{vah}, since $X$  is the
$T_1$-space and for each point $x$ and each closed subset $S$ in $X$
with $x \notin S$ there is a continuous function $h: X \to B({\bf K},0,1)$
such that $h(x)=0$ and $h(S)=\{ 1\} $.
\par {\bf 2.17. Theorem.} {\it Let $X$ be a Banach space over $\bf K$,
$\eta :\Gamma \to \bf C$ be a continuous positive definite function,
$(\mu
_{\beta })$ be a bounded weakly relatively compact net in the
space ${\sf M}_t(X)$ of Radon norm-bounded measures 
and there exists $\lim _{\beta } \hat \mu _{\beta } (f)=
\gamma (f)$ for each $f \in \Gamma $ 
and uniformly on compact subsets of the completion $\tilde \Gamma $,
where $\Gamma \subset C(X,{\bf K})$
is a vector subspace separating points in $X$. Then
$(\mu _{\beta })$ weakly converges to $\mu \in {\sf M}_t(X)$ with
$\hat \mu |_{\Gamma }=\gamma .$}
\par  {\bf Proof} is analogous to the proof of Theorem IV.3.1\cite{vah}
and follows from Theorem 2.9 above.
\par  {\bf 2.18. Theorem.} {\it (a). A bounded
family of measures in ${\sf M}({\bf K^n})$ 
is weakly relatively compact if and only if a family
$(\hat \mu: $ $\mu \in M)$ is equicontinuous on $\bf K^n$. 
\par (b). If $(\mu _j:$ $j \in {\bf N})$ 
is a bounded sequence of measures
in ${\sf M}_t({\bf K^n})$, $\gamma : {\bf K^n} \to \bf C$
is a continuous and positive definite function, 
$\hat \mu _j(y) \to \gamma (y)$ for each $y \in \bf K^n$,
then $(\mu _j)$ weakly converges to a measure
$\mu $ with $\hat \mu =\gamma .$ 
\par (c). A bounded sequence of measures
$(\mu _j)$ in ${\sf M}_t({\bf K^n})$ 
weakly convereges to a measure
$\mu $ in ${\sf M}_t({\bf K^n})$ if and only if for each $y \in \bf K^n$
there exists $\lim _{j \to \infty }\hat \mu _j(y)=\hat \mu (y).$
\par (d). If a bounded net $(\mu _{\beta })$ in ${\sf M}_t({\bf K^n})$
converges uniformly on each bounded subset
in $\bf K^n$, then $(\mu _{\beta })$ converges weakly to a
measure $\mu $ in ${\sf M}_t({\bf K^n})$, where $n \in \bf N$.}
\par {\bf Proof.} (a). This follows from the Prohorov theorem 1.3.6\cite{vah}
and Propositions 2.7, 2.16. 
\par (b). We have the following inequality:
$\lim_m \sup_{j>m}
\mu _j([x \in \bf K^n:$ $|x| \ge
tR]) $ $\le 2 \int_{\bf K^n} (1-Re(\eta (\xi y)))\nu (dy)$ 
with $|\xi |=1/t$ due to
\S 2.7 and \S 2.8.
In view of Theorem 2.17 $(\mu _j)$ converges weakly to
$\mu $ with
$\hat \mu =\gamma $. (c,d). These can be proved analogously to IV.3.2\cite{vah}.
\par  {\bf 2.19. Corollary.} {\it If $({\hat \mu }_{\beta })\to 1$ 
uniformly on some
neighbourhood of $0$ in $\bf K^n$ for a bounded net
of measures $\mu _{\beta } $ in ${\sf M}_t({\bf K^n})$, then
$(\mu _{\beta })$ converges weakly to $\delta _0$.}
\par  {\bf 2.20. Definition.}  A family of probability measures
$M \subset
{\sf M}_t(X)$ for a Banach space 
$X$ over $\bf K$ is called planely concentrated if
 for each $c>0$ there exists a $\bf K$-linear subspace $S \subset X$ with
$dim_{\bf K}S=n< \aleph _0 $ such that $\inf (\mu
(S^c): $ $\mu \in M)>1-c$. The Banach space ${\sf M}_t(X)$ is 
supplied with the following norm
$\| \mu \|:=|\mu |(X).$
\par {\bf 2.21. Lemma.} {\it Let $S$ and $X$ be the same as in \S 2.20;
$z_1,...,z_m \in X^*$ be a separating family of points in $S$. Then a set
$E:=S^c \cap (x \in X:$ $|z_j(x)| \le r_j; $ $j=1,...,m)$ is bounded
for each $c>0$ and $r_1,...,r_m\in (0,\infty ).$}
\par  {\bf Proof.} A space $S$ is isomorphic with $\bf K^n$, consequently,
$p(x)=\max(|z_j|:$ $j=1,...,m)$ is the norm in $S$ equivalent to
the initial norm.
\par {\bf 2.22. Theorem.} {\it Let $X$ be a Banach space 
over $\bf K$ with a family
$\Gamma \subset X^*$ separating points in $M\subset {\sf M}_t(X)$. Then $M$
is weakly relatively compact if and only if a family $ \{ \mu _z:$
$\mu \in M \} $ is weakly relatively compact for each $z \in \Gamma $
and $M$ is planely concentrated, where $\mu _z$ is an image measure 
on $\bf K$ of a measure $\mu $ induced by $z$.}
\par  {\bf Proof} follows from Lemmas 2.5, 2.21 and the Prohorov theorem
(see also Theorem 1.3.7\cite{vah} with a substitution $[-r_j,r_j]$
onto $B({\bf K},0,r_j)$).
\par {\bf 2.23. Theorem.} {\it For $X$ and $\Gamma $ the same as
in Theorem 2.22 a sequence $ \{ \mu _j:$ $j \in {\bf N} \} \subset
{\sf M}_t(X)$  is weakly convergent to $\mu \in {\sf M}_t(X)$
if and only if for each $z \in \Gamma $ there exists
$\lim _{j\to \infty } {\hat \mu }_j(z)={\hat \mu }(z)$ and a family
$\{ \mu _j \} $ is planely concentrated.}
\par  {\bf Proof}  follows from Theorems 2.17,18,22 (see also
Theorem IV.3.3\cite{vah}).
\par  {\bf 2.24. Proposition.} {\it Let $X$ be a weakly regular
space with $ind(X)=0$, $\Gamma \subset C(X,{\bf K})$ be a vector subspace
separating points in $X$, $(\mu _n:$ $ n \in {\bf N})$ $\subset
{\sf M}_t(X),$ $\mu \in {\sf M}_t(X)$, $\lim_{n \to \infty }
\hat \mu _n(f)=\hat \mu (f)$ for each $f \in \Gamma $. Then
$(\mu _n)$ is weakly convergent to $\mu $ relative to the weakest
topology $\sigma (X,\Gamma )$ in $X$ relative to which all
$f \in \Gamma $ are continuous.}
\par  {\bf Proof} follows from Theorem 2.18 and is analogous to
the proof of Proposition IV.3.3\cite{vah}.
\par  {\bf 2.25.} Let $(X,{\sf U})$ $=\prod_{\lambda }
(X_{\lambda },{\sf U}_{\lambda })$ be a product of measurable
completely regular Radon spaces $(X_{\lambda },{\sf U}_{\lambda })$
$=(X_{\lambda },{\sf U}_{\lambda },{\sf K}_{\lambda })$, where
${\sf K}_{\lambda }$ are compact classes approximating from below
each measure $\mu _{\lambda }$ on $(X_{\lambda },{\sf U}_{\lambda })$,
that is, for each
$c>0$  and elements $A$ of an algebra ${\sf U}_{\lambda }$
there is $S \in {\sf K}_{\lambda }$, $S \subset A$
with $\| A\setminus S\| _{\mu _{\lambda }} <c$. 
\par {\bf 2.26. Definition.} Let $X$ be a Banach space 
over $\bf K$, then a mapping
$f:X\to \bf C$ is called pseudocontinuous, if its restriction $f|_L$
is uniformly continuous for each subspace  $L\subset X$ with
$dim_{\bf K}L<\aleph _0$. Let $\Gamma $ be a family of mappings $f:Y\to \bf K$
of a set $Y$ into a field $\bf K$. We denote by $\hat
{\sf C}(Y,\Gamma )$ the minimal $\sigma $-algebra (that is called
cylindrical) generated by an algebra ${\sf C}(Y,\Gamma )$
of subsets of the form $C_{f_1,...,f_n;E}:= \{ x\in X:$ $(f_1(x),...,f_n(x))
\in S \} $, where $S \in Bf({\bf K^n})$, $f_j \in \Gamma $. 
We supply $Y$ with
a topology $\tau (Y)$ which is generated by a base
$(C_{f_1,...,f_n;E}:$ $f_j \in \Gamma ,$
$ E $ $\mbox{is open in}$ $\bf K^n)$.
\par  {\bf 2.27. Theorem. Non-Archimedean analog of the Bochner-Kolmogorov
theorem.} {\it Let $X$ be a Banach space 
over $\bf K$, $X^a$ be its algebraically
adjoint $\bf K$-linear space (that is, of all linear mappings
$f:X\to \bf K$ not necessarily continuous). A mapping $\theta : X^a\to \bf C$
is a characteristic functional of a probability measure
$\mu $ with values in $\bf R$ and is defined on
$\hat {\sf C}(X^a,X)$ if and only if $\theta $ satisfies conditions
2.6(3,5) for $(X^a,\tau (X^a))$ and is pseudocontinuous on $X^a$.}
\par  {\bf Proof.} (I). For $dim_{\bf K}X=card(\alpha )<\aleph _0$
a space $X^a$ is isomorphic with $\bf K^{\alpha }$, hence the statement of
theorem for a measure $\mu $ follows from
Theorems 2.9 and 2.18 above, since
$\theta (0)=1$ and $|\theta (z)|\le 1$ for each $z$.
\par (II). We consider now the case of $\mu $ with values in $\bf R$
and $\alpha < \omega _0$. In \S 2.6 
(see also \S 2.16-18,24) it was proved that $\theta =\hat \mu $
has the desired properties for real probability measures $\mu $.
On the other hand, there is $\theta $ which satisfies the conditions of
the theorem. Let $\theta _{\xi } (y)=\theta (y)h_{\xi }(y)$, where
$h_{\xi }(y)=F[C(\xi )exp(-\| x \xi \| ^2)](y)$ (that is, the Fourier transform
by $x$), $\nu _{\xi }({\bf K^{\alpha }})=1$, $\nu _{\xi }(dx)=C(\xi )
exp(-\| x \xi \| ^2)m(dx)$ (see Lemma 2.8), $\xi \ne 0$.
Then $\theta _{\xi }(y)$ is positive definite and is uniformly continuous
as a product of two such functions. Moreover, $\theta
_{\xi }(y) \in L^1({\bf K^{\alpha }},m,{\bf C})$. For $\xi \ne 0$
a function $f_{\xi }(x)=$ $\int_{\bf K^{\alpha }} \theta _{\xi }(y)
\chi _e(x(y))m(dy)$ is bounded and continuous, the function
$exp(-\| x \xi \| ^2)=:s(x)$ is positive definite. Since $\nu _{\xi }$
is symmetric and weakly converges to $\delta _0$,
hence there exists $r>0$ such that for each $|\xi |>r$ we have
${\hat \gamma }_{\xi }(y)=$ $\int _{\bf K^{\alpha }} C(\xi )exp(-\| $ $x \xi
\| _p^2)\chi _e(y(x))m(dx)$ $=\int [\chi _e((y(x))+\chi _e(-y(x))]2^{-1}
exp(-\| x\xi \| _p^2)C($ $\xi )m(dx)/2$ $>1-1/R$ for $|y| \le R$,
consequently, ${\hat \gamma }_{\xi }(y)={\hat \zeta }_{\xi }^2(y)$
for $|y| \le R$, where ${\hat \zeta }_{\xi }$ is positive definite
uniformly continuous and has a uniformly continuous extension on
$\bf K^{\alpha }$. Therefore, for each $c>0$ there exists $r>0$ such that
$\| \nu _{\xi }-\kappa _{\xi }*\kappa _{\xi } \| <c$ for each
$|\xi |>r$, where $\kappa _{\xi }(dx)=\zeta _{\xi }(x)m(dx)$
is a $\sigma $-additive non-negative measure.
Hence due to corollary from Proposition IV.1.3\cite{vah}
there exists $r>0$
such that $\int_{\bf K^{\alpha }} \theta _{\xi }(y) \chi _e(-x(y))
\nu _j(dy) \ge 0$ for each $|j|>r$, consequently,
$f_{\xi }(x)=\lim _{|j| \to \infty } \int _{\bf K^{\alpha }}
\theta _{\xi }(y)$ $ \chi _e(-x(y))\nu _j(dy) \ge 0$.
From the equality $F[F(\gamma _{\xi })(-y)](x)=\gamma _{\xi }(x)$
and the Fubini theorem it follows that 
$$\int f_{\xi }\chi _e(y(x))h_j
(x)m(dx)=\int \theta _{\xi }(u+y)\nu _j(du).$$
For $y=0$ we get
$$\lim _{|\xi | \to \infty } \int f_{\xi }(x)m(dx)=\int f(x)m(dx)
=\lim_{|\xi | \to \infty } \lim_{|j| \to \infty } \int f_{\xi }
(x)h_j(x)m(dx)$$ $$\mbox{and  }\lim_{|\xi |\to \infty } \lim_{|j|
\to \infty } |\int_{\bf K^{\alpha }}
\theta _{\xi }(u)\nu _j(du)| \le 1.$$ From Lemma 2.8 it follows that
$\hat f(y)=\theta (y)$, since by Theorem 2.18 $\theta =\lim _{|\xi |
\to \infty } \theta _{\xi }$ is a characteristic function of a probability
measure on $Bf({\bf K^{\alpha }})$, where $f(x)=
\int_{\bf K^{\alpha }} \theta (y)\chi _e(-x(y))m(dy)$.
\par  (III). Now let $\alpha =\omega _0$. It remains to show that
the conditions imposed on $\theta $ are sufficient, because their necessity
follows from the modification of 2.6 (since $X$ has an algebraic embedding
into $X^a$).
The space $X^a$ is isomorphic with ${\bf K}^{\Lambda }$ which is the space
of all $\bf K$-valued functions defined on the Hamel basis $\Lambda $
in $X$. The Hamel basis exists due to the Kuratowski-Zorn lemma
(that is, each finite system of vectors in $\Lambda $ is linearly independent
over $\bf K$, each vector in $X$ is a finite linear combination over
$\bf K$ of elements from $\Lambda $). Let $J$ be a family of all non-void
subsets in $\Lambda $. For each $A \in J$ there exists a functional
$\theta _A: {\bf K}^A \to \bf C$ such that $\theta _A(t)=
\theta (\sum_{y \in A}t(y)y)$ for $t \in {\bf K}^A$. From the conditions
imposed on $\theta $ it follows that $\theta _A(0)=1$, $\theta _A$ is
uniformly continuous and bounded on ${\bf K}^A$, 
moreover, it is positive definite
(or due to 2.6(6) for each $c>0$ there are $n$ and $q>0$ such that
for each $j>n$ and $z \in {\bf K}^A$ the following inequality is satisfied:
$$(i)\mbox{  }|\theta _A(z)-\theta _j(z)| \le cbq,$$
moreover, $L_j \supset {\bf K}^A$, $q$ is independent on $j$, $c$ and $b$.
From (I,II) it follows that on $Bf({\bf K}^A)$ 
there exists a probability measure
$\mu _A$ such that $\hat \mu _A =\theta _A$. The family of measures
$ \{ \mu _A:$ $A \in J \} $ is consistent and bounded, since
$\mu _A=\mu _E \circ
(P_E^A)^{-1}$, if $A \subset E$, where $P_E^A: {\bf K}^E \to
{\bf K}^A$ are the natural projectors. Indeed, in the case of measures with
values in $\bf R$ each $\mu _A$ is the probability measure. 
\par  In view of Theorem 1.1.4\cite{dal} on
a cylindrical $\sigma $-algebra of the space ${\bf K}^{\Lambda }$
there exists the unique measure $\mu $ such that $\mu _A=\mu
\circ (P^A)^{-1}$ for each $A \in J$, where $P^A: {\bf K}^{\Lambda }\to 
{\bf K}^A$
are the natural projectors. From $X^a={\bf K}^{\Lambda }$ it follows that
$\mu $ is defined on $\hat {\sf C}(X^a,X)$.
For $\mu $ on ${\hat C}(X^a,X)$ there exists its extension
on $Af(X,\mu )$ (see \S 2.1).
\par  {\bf 2.28. Definition. \cite{sch2} } A continuous linear operator
$T:X\to Y$ for Banach spaces
$X$ and $Y$ over $\bf K$ is called compact, if
$T(B(X,0,1))=:S$ is a compactoid, that is, for each neighbourhood
$U \ni 0$ in $Y$ there exists a finite subset
$A \subset Y$ such that $S \subset U
+co(A)$, where $co(A)$ is the least $\bf K$-absolutely convex subset
in $V$ containing $A$ (that is,
for each $a$ and $b\in \bf K$ with $|a|\le 1$,
$|b| \le 1$ and for each $x,y \in V$ the following inclusion
$ax+by \in V$ is accomplished).
\par  {\bf 2.29.} Let $\sf B_+$ be a subset of non-negative
functions which are $Bf(X)$-measurable and let $\sf C_+$ be its subset
of non-negative cylindrical functions. By ${\hat {\sf B}}_+$ we denote
a family of functions $f \in \sf B_+$ such that $f(x)=\lim_{n}g_n(x)$,
$g_n\in \sf C_+$, $g_n \ge f$. For $f \in {\hat {\sf B}}_+$ let
$$\int_X f(x)\mu _*(dx):=\inf_{g \ge f, \mbox{  }g \in {\sf C_+}}
\int_Xg(x)\mu _*(dx).$$
\par  {\bf 2.30. Lemma.} {\it A sequence of a weak distribution
$(\mu _{L_n})$
of probability Radon measures is generated by a real probability measure
$\mu $ on $Bf(X)$ of a Banach space 
$X$ over $\bf K$ if and only if there exists
$$(i)\mbox{  }\lim_{|\xi |\to \infty }\int_X G_{\xi }(x)\mu_*(dx)=1,$$
where $\int_X G_{\xi }(x)\mu _*(dx):=S_{\xi }(\{ \mu _{L_n}:n \})$  and \\
$S_{\xi }(\{ \mu _{L_n} \}):=\lim_{n \to \infty }$ $\int_{L_n}
F_n(\gamma _{\xi ,n})(x)$ $\mu _{L_n}(dx),$
$\gamma _{\xi ,n}(y):=\prod_{l=1}^{m(n)} \gamma _{\xi }(y_l)$, \\
$F_n$ is a Fourier transformation by $(y_1,...,y_n)$, $y=(y_j:$ $j \in \bf N)$,
$y_j \in \bf K$,
$\gamma _{\xi }(y_l)$ are the same as in Lemma 2.8 for ${\bf K}^1$; here
$m(n)=dim_{\bf K}L_n< \aleph _0$, $cl(\bigcup_n L_n)=X=c_0(\omega _0,K)$.}
\par  {\bf Proof.} If a sequence of weak ditributions 
is generated by a measure
$\mu $, then in view of Conditions 2.6(3-6), Lemmas 2.3, 2.5, 2.8, Propositions
2.10 and 2.16, Corollary 2.13, the Lebesgue convergence theorem
and the Fubini theorem, also from the proof of Theorem 2.27 and the Radon
property of $\mu $
it follows that there exists $r>0$ such that 
$$\int_X G_{\xi }
(x)\mu _*(dx)=\int_X G_{\xi }(x)\mu (dx)=\lim_{n \to \infty }
\int_{L_n} \gamma _{\xi ,n}(y) \hat \mu _{L_n}(y) m_{L_n}(dy),$$
since $\lim _{j \to \infty }
x_j=0$ for each $x=(x_j:$ $ j) \in X$. In addition, $\lim _{|\xi |
\to \infty }
S_{\xi }(\{ \mu _{L_n} \})$ $=\int_X \mu (dx)=1$. Indeed,
for each $c>0$ and $d>0$ there exists a compact $V_c \subset X$
with $\| \mu |_{(X\setminus
V_c)} \| <c$ and there exists $n_0$ with $V_c \subset L_n^d$ for each
$n>n_0$.
Therefore, choosing suitable sequences of $c(n)$, $d(n)$, $V_{c(n)}$ and
$L_{j_n}$ we get that $[\int_{L_n} \gamma _{\xi ,n}(y) \hat \mu _{L_n}
(y)m_{L_n}(dy):$ $n \in {\bf N}]$ is a Cauchy sequence, where
$m_{L_n}$ is the real Haar measure on $L_n$, the latter is considered as
$\bf Q_p^{m(n)b}$, $b=dim_{\bf Q_p}{\bf K}$, $m(B(L_n,0,1))=1$.
Here we use $G_{\xi }(x)$ for a formal expression of the limit
$S_{\xi }$ as the integral. Then $G_{\xi }(x)$
$(mod$ $\mu )$ is defined evidently as a function for $\mu $ or
$ \{ \mu_{L_n}:$ $n \} $ with a compact support, also for  $\mu $ with
a support in a finite-dimensional subspace $L$ over $\bf K$ in $X$.
By the definition $supp(\mu _{L_n}:$ $n)$ is compact, if
there is a compact $V
\subset X$ with $supp(\mu _{L_n}) \subset P_{L_n}V$ for each
$n$.
That is, Condition (i) is necessary.
\par On the other hand, if (i) is satisfied, then for each
$c>0$ there exists $r>0$ such that $|\int_X G_{\xi }(x) \mu _*(dx)-1|<
c/2$, when
$|\xi | >r$, consequently, there exists $n_0$ such that for each
$n>n_0$ the following inequality is satisfied:
$$|1-\int_XF_n(\gamma _{\xi ,n})(x)\mu _*(dx)|
\le | \| \mu |_{(L_n\cap B(X,0,R))} \| -1| +$$
$$\sup_{|x|>R} |F_n(\gamma _{\xi ,n})(x)| \| \mu _{L_n}|_{
(L_n\setminus B(X,0,R))} \| .$$
Therefore, from
$\lim_{R \to \infty } \sup _{|x| >R} |F_n(\gamma _{\xi ,n})(x)|=0$
and from Lemma 2.3 the statement of Lemma 2.30 follows.
\par {\bf 2.31. Notes and definitions.}
Suppose $X$ is a locally convex space over a locally compact field
$\bf K$ with non-trivial non-Archimedean valuation
and $X^*$ is a topologically adjoint space.
The minimum $\sigma $-algebra with respect to which
the following family $ \{ v^* : v^*\in X^* \} $
is measurable is called a $\sigma $-algebra of cylindrical sets.
Then $X$ is called a $RS$-space
if on $X^*$ there exists a topology $\tau $ such that the 
continuity of each positive definite function $f: X^*\to \bf C$
is necessary and sufficient for $f$ to be a characteristic 
functional of a non-negative measure. Such topology
is called the $R$-Sazonov type topology.
The class of $RS$-spaces contains
all separable locally convex spaces over $\bf K$.
For example, $l^{\infty }(\alpha ,{\bf K})=c_0(\alpha ,{\bf K})^*$,
where $\alpha $ is an ordinal \cite{roo}.
In particular we also write $c_0({\bf K}):=c_0(\omega _0,{\bf K})$ and
$l^{\infty }({\bf K}):=l^{\infty }(\omega _0,{\bf K})$,
where $\omega _0$ is the first countable ordinal.
\par Let $n_{\bf K}(l^{\infty },c_0)$ denotes the weakest topology on 
$l^{\infty }$ for which all functionals $p_x(y):=\sup_n|x_ny_n|$ are 
continuous, where $x=\sum_nx_ne_n\in c_0$ and $y=\sum_ny_ne^*_n
\in l^{\infty }$, $e_n$ is the standard base in $c_0$.
Such topology $n_{\bf K}(l^{\infty },c_0)$ is called the normal topology.
The induced topology on $c_0$ is denoted by $n_{\bf K}(c_0,c_0)$.
\par {\bf 2.32. Theorem.} {\it Let $f: l^{\infty }({\bf K})
\to \bf C$ be a functional such that
\par $(i)$ $f$ is positive definite,
\par $(ii)$ $f$ is continuous in the normal topology
$n_{\bf K}(l^{\infty },c_0)$, then $f$ is the characteristic functional
of a probability measure on $c_0({\bf K})$.}
\par {\bf Proof.} The case of the topological vector space $X$ over
$\bf K$ with $char({\bf K})>0$ and a real-valued measure
$\mu $ can be proved analogously to \cite{mamaz} for $char ({\bf K})=0$
due to \S 2.6 and \S \S 2.25-2.30 (see also the proofs in Part II
for $\bf K_s$-valued measures).
\par {\bf 2.33. Theorem.} {\it Let $\mu $ be a probability measure
on $c_0({\bf K})$, then $\hat \mu $ is continuous in 
the normal topology $n_{\bf K}(l^{\infty },c_0)$ on $l^{\infty }$.}
\par {\bf Proof.} In view of Lemma 2.3
for each $\epsilon >0$ there exists $S(\epsilon )\in c_0$ such that
$\| \mu |_{L(0,S(\epsilon ))} \| \ge 1 -\epsilon $,
where $L(y,z):= \{ x\in c_0:$ $|x_n-y_n|\le |z_n|,$ for each 
$n\in {\bf N} \} $. Therefore, 
$$|1-{\hat \mu }(x)|
\le \epsilon +\| 2\pi \eta (\xi x) \|_{C^0(L(0,S(\epsilon )))} 
\| \mu |_{L(0,S(\epsilon ))} \| ,$$ 
hence there exists a constant
$C>0$ such that $|1-{\hat \mu }|\le \epsilon +Cp_{S(\epsilon )}(x)$.
\par {\bf 2.34. Corollary.} {\it The normal topology $n_{\bf K}
(l^{\infty },c_0)$ is the $R$-Sazonov type topology
on $l^{\infty }({\bf K})$.}
\par  {\bf 2.35. Theorem. Non-Archimedean analog of the
Minlos-Sazonov theorem.}
{\it For a separable Banach space 
$X$ over $\bf K$ the following two conditions are equivalent:
$$(I)\mbox{ }\theta : X \to {\bf T}
\mbox{ satisfies conditions }2.6(3,4,5)
\mbox{ and }$$ for each $c>0$ there exists a compact operator
$S_c:X\to X$
such that $|Re(\theta (y)-\theta (x))| <c$ for $|\tilde z (S_cz)|<1$;
$$(II)\mbox{ }\theta \mbox{ is a characteristic functional of a 
probability Radon measure } \mu $$ on $E$, where $\tilde z$ is an element
$z \in X
\hookrightarrow X^*$ considered as an element of $X^*$ under the
natural embedding associated with the standard base of 
$c_0(\omega _0,{\bf K})$, $z=x-y$, $x$
and $y$ are arbitrary elements of $X$.}
\par {\bf Proof.} $(II\to I)$. For a positive definite function $\theta $
generated by a probability measure $\mu $ in view of the inequality
$|\theta (y)-\theta (x)|^2 \le 2 \theta (0)(\theta (0)-Re(\theta
(y-x))$ (see Propositions IV.1.1(c)\cite{vah}) and using the
normalization of a measure $\mu $ by $1$
we consider the case $y=0$. For each $r>0$ we have: $|Re(\theta (0)-
\theta (x))|=\int_X (1-[\chi _e(x(u))+\chi _e(-x(u))]/2) \mu (du)$ $\le
\int_{B(X,0,r)} 2 [\chi _e^{1/2}(x(u))-\chi _e^{1/2}(-x(u))/(2i)]^2 \mu (du)$ 
$+2\int_{
X\setminus B(X,0,r)} $ $\mu (du)$ $\le 2 \pi ^2 \int_{B(X,0,r)}
\eta (x(u))^2 \mu (dx)$ $+2 \mu([x:$ $\| x\| >r])$. In view of the
Radon property of the space $X$ and Lemma 2.5 for each $b>0$ and
$\delta >0$ there are a finite-dimensional over $\bf K$ subspace
$L$ in $X$ and a compact subset $W \subset X$ such that
$W \subset L^{\delta }$, $\| \mu |_{(X\setminus W)} \| <b$, 
hence $\| \mu |_{(X\setminus L^{\delta })} \| <b.$
\par We consider the following expression:
$$J(j,l):=2 \pi ^2 \int_{B(X,0,r)}
\eta (e_j(u)) \eta (e_l(u)) \mu (du),$$ 
where $(e_j)$ is the orthonormal basis
in $X$ which contains the orthonormal basis 
of $L=\bf K^n$, $n=dim_{\bf K}L$. Then we choose
sequences $b_j=p^{-j}$ and
$0< \delta _j <b_j$, subspaces $L_j$ and
$r=r_j$ such that $b_jr_j<1,$ $W_j \subset B(X,0,r_j)$, $0<r_j<r_{j+1}<
\infty $ for each $j \in \bf N$ and the orthonormal basis 
$(e_j)$ corresponding to the sequence
$L_j \subset L_{j+1} \subset ... \subset X$.
We get, due to finiteness of $n_j:=dim_{\bf K}L_j$, that
$\lim_{j+l \to \infty } J(j,l)=0$,
since $\| \mu |_{ \{ x: \| x\| >r_j \} } \| <b_j$,  
$\eta (x(u))=0$ for $x \in X \ominus
L_j$ with $\| x\| <b_j$, $u \in B(X,0,r_j)$. Then we define $g_{j,l}
:=\min \{ d:$ $d \in \Gamma _{\bf K}$  and  $d \ge |J(j,l)| \} $, evidently,
$g_{j,l} \le p|J(j,l)|$ and there are $\xi _{j,l} \in \bf K$
with $|\xi _{j,l}|_{\bf K}=g_{j,l}$. Consequently, the family $(\xi _{j,l})$
determines a compact operator $S: X \to X$ with $\tilde e_j(Se_l)=
\xi _{j,l}t$
due to Theorem 1.2\cite{sch2}, where $t=const \in \bf K$, $ t\ne 0$. Therefore,
$|Re(\theta (0) -\theta (z))|<c/2 +|\tilde z(Sz)|<c$ 
$|\theta (0)-\theta (z)|<c/2 + |{\tilde z}(Sz)|<c$.
We choose $r$ such that $\| \mu |_{(X\setminus B(X,0,r))} \| <c/2$ 
with $S$ corresponding to $(r_j:$ $j)$, where $r_1=r$,
$L_1=L$, then we take $t \in \bf K$ with $|t|c=2$.
\par  $(I \to II)$. Without restriction of generality we may take
$\theta (0)=1$
after renormalization of non-trivial $\theta $. 
In view of Theorem 2.32 as in \S 2.6 
we construct using $\theta (z)$ a consistent family
of finite-dimensional distributions $\{ \mu _{L_n} \} $. Let $m_{L_n}$
be a real Haar measure 
on $L_n$ which is considered as $\bf Q_p^a$ with
$a=dim_{\bf K}L_ndim_{\bf Q_p}{\bf K}$, $m(B(L_n,0,1))=1$. In view of
Proposition 2.7 and Lemmas 2.8, 2.30: $\int_{L_n} G_{\xi }(x) \mu
_{L_n} (dx)$ $=\int _{L_n} \gamma _{\xi ,n}(z) \theta (z)m_{L_n}(dz)$,
consequently, 
$$1-\int_{L_n} F_n(\gamma _{\xi ,n})(x) \mu _{L_n} (dx)=
\int \gamma _{\xi ,n}(z) (1-\theta (z)) m_{L_n} (dz)=:I_n(\xi ).$$
There exists an orthonormal basis 
in $X$ in which $S_c$ can be reduced to the following
form $S_c=SC\hat S_cE$ (see Appendix), where $\hat S_c=diag(s_j:$
$j \in {\bf N})$ in the orthonormal basis 
$(f_j :j)$ in $X$ and $S$ transposes a finite
number of vectors in the orthonormal basis.
That is, $|\tilde z(\hat S_cz)|=\max_j |s_j|\times |z_j|^2$.
In the orthonormal basis $(e_j: j)$ adopted to
$(L_n:$ $n)$ we have $|\tilde z(S_cz)|=\max_{j,l
\in \bf N}(|s_{j,l}|\times |z_j| \times |z_l|)$, $\| S_c\| =\max_{j,l}
|s_{j,l}|$, where $S_c=(s_{j,l}:$ $j,l\in {\bf  N})$ in the 
orthonormal basis $(e_j)$, $r=const>0$. 
In addition, $p^{-1}|x|_K \le |x|_p \le p|x|_{\bf K}$ for each
$x \in \bf K$.
If $S_c$ is a compact operator such that  $|Re(\theta (y)-
\theta (x))|<c$ for $|\tilde z(S_cz)|<1$, $z=x-y$, then 
$|Re(1-\theta (x))|< c+2|\tilde x(S_cx)|$ and 
$$I_n(\xi ) \le \int \gamma _{\xi ,n}
(z) [c+2|\tilde z(S_cz)|_K] m_{L_n} (dz)\le c+b\| S_c\| /|\xi |^2,$$
$b=const $ is independent from $n$, $\xi $ and $S_c$, 
$$b:=p\times 
\sup_{|\xi |>r}
|\xi |^2 \int_{L_n} \gamma _{\xi ,n}(z) |z|_p^2 m_{L_n} (dz)< \infty .$$
Due to the formula of changing variables in integrals (A.7\cite{sch1})
the following equality is valid: $J_n(\xi )=I_n(\xi )J_n(1)/
[I_n(1)|\xi |^2]$ for $|\xi |
\ne 0$, where $$J_n(\xi )=\int_{L_n} \gamma _{\xi ,n}(z)|z|_p^2 m_{L_n}
(dz).$$ Therefore, 
$$1-\int_X G_{\xi }(x) \mu _*(dx) \le c+b\| S_c \| /|\xi |^2.$$ 
Then taking the limit with $|\xi | \to \infty $ and then
with $c \to +0$ with the help of Lemma 2.30 we get the statement $(I\to II)$.
\par  {\bf 2.36. Definition.}
 Let on a completely regular space
$X$ with $ind (X)=0$ two non-zero real-valued
measures $\mu $ and $\nu $ are given.
Then $\nu $ is called absolutely continuous relative to $\mu $ if $\nu (A)=0$
for each $A \in Bf(X)$ with $\mu (A)=0$
and it is denoted $\nu \ll \mu $. Measures $\nu $ and $\mu $ are
singular to each other if there is $F\in Bf(X)$ with
$|\mu |(X\setminus F)=0$ and $|\nu |(F)=0$ and it is denoted
$\nu \perp \mu $. If $\nu \ll \mu $ and $\mu \ll \nu $ then they are called
equivalent, $\nu \sim \mu $.
\par {\bf 2.37.} For a Banach space 
$X$ over $\bf K$ and $\sigma $-algebra ${\sf B}
\supset Bf(X)$ with a real probability measure $\mu $ and
a $\sigma $-subalgebra $\sf B_0$ $\subset \sf B$ analogously to
\S 13\cite{sko} a function $\bar \mu (A|x)$ satisfying three conditions:
\par $(a)$ $\bar \mu (A|x)$ is $\sf B_0$-measurable by $x$ for each
$A \in \bf \sf B$; 
\par $(b)$ $\mu (A\cap A_0)=\int_{A_0}\bar \mu
(A|x)\mu (dx)=:\mu _A(A_0)$ for each $A_0 \in \sf B_0$;
\par $(c)$ $\bar \mu (A|x)$ is a measure by $A \in \sf B$ for almost all
$x$ relative to a measure $\mu $, then it is called a conditional
measure corresponding to $\mu $ relative to a $\sigma $-algebra
$\sf B_0$.
\par  Then we define $\mu (A|x)$ by formula (b) and we then redefine
it on the set of measure zero relative to $\mu $ for each $A$. Let
$\theta (z|x)=
\theta _1(z|x)+i\theta _2(z|x)$ with $\int_{A_0} [[\chi _e(\tilde z
(x))+\chi _e(-\tilde z(x))]/2] \mu (dx)$ $=\int_{A_0}\theta _1(z|x) \mu (dx)$;
$\int_{A_0} [\chi _e(\tilde z(x))-\chi _e(-\tilde z(x))]/(2i)] \mu (dx)$ $=\int_{A_0}
\theta _2($ $z|x) \mu (dx)$, where $\theta _j$ are measurable functions,
$z \in X$, $\tilde z \in X'$ corresponds to $z$ under the natural
embedding of $X$ into $X'$, $x \in X$, $A_0 \in \sf B_0$.
Then we choose analogously to \cite{sko} sets $A_{n,k}$ with
$diam(A_{n,k})<p^{-n}$, using \S \S 2.7, 2.8, 2.30 and
the proof of \S 2.27, substituting $exp(-(x,x)/n)$
on $G_{\xi }(x)$, defining sets $B_n$, $G$, $G_1,$ we get
for each $c>0$ that there exists $r>0$ such that 
$$Re(1-\theta (z|x))
\le \int_XG_{\xi }(y)(1-[[\chi _e(\tilde z(y))+\chi _e(\tilde z(x))]/2]
\mu (dy|x) +2\int_X(1-G_{\xi }(y)) \mu (dy|x) $$ 
$$\le \int (2\pi \eta ^2(\tilde z
(y)) G_{\xi }(y) \mu (dy|x)/2 +c/2$$  
for each $|\xi |>r$
and $x \in G_1 \cap G$, since from $\| z\| \times \| y\| \le 1$
it follows that $(1-
\chi _e(\tilde z(y))=0$ and $|1-exp(it)| \le |t|$ for $t \in \bf R$,
$i=(-1)^{1/2}$.
\par For $C(z_1,z_2)(x):=\int_X\eta (\tilde z_1(y))
\eta (\tilde z_2(y))G_{\xi }(y)\mu (dy|x)$ we have $C(0,0)=0$ and
$C(z,z)(x)=\int_{X\setminus B(X,0,r)} \eta ^2(\tilde z(y))G_{\xi }
(y) \mu (dy|x)$ for $\|z\|=1/r>0$. Since $\theta (z|x)=
\int_X \chi _e(\tilde z(y)) \mu (dy|x)$, considering projectors
$P_{L_n}:X\to L_n$ onto finite-dimensional subspaces over $\bf K$,
$L \subset X$, we get $\theta _{L}(z|x)$ corresponding to
$\mu _{L}(dx)$ and $\mu _L(dy|x)$. At the same time $\theta _L$
satisfies Conditions 2.6(3,4,5), consequently,  $\mu _L(A|x)$ are
$\sigma $-additive
by  $A$ for fixed $x$, since $\theta (z|x)$
are continuous by $z$ for $x \in G\cap G_1$. From the Prohorov theorem
IX.4.2.1\cite{boui}, Lemma 2.5, $x \in G\cap G_1$ and $C(e_j,e_l)(x)=
\int_{L_{j,l}\setminus B(L_{j,l},0,1)}\eta (y_l) \eta (y_j)
 F_2(\gamma _{\xi ,2})(y_j,y_l) \mu _{L_{j,l}}(dy_jdy_l|x)$ it follows that
$\lim_{j+l\to \infty }\int_{L_{j,l}} \psi (x)C(e_j,e_l)(x)
\mu _{L_{j,l}}(dx)$ $=\lim _{j+l\to \infty }\int _{L_{j,l}}
\psi (x) \int_{L_{j,l}} \eta (y_j)$ $\eta (y_l)\chi _e(\tilde z(y))
\mu _{L_{j,l}}(dy|x)\mu _{L_{j,l}}(dx)=0$ for each bounded
$\sf B_0$-measurable function $\psi (x)$, whence $\lim_{n \to
\infty } \mu ([x:$ $|C(e_j,e_l)(x)|>c$ for some $j+l>n])=0$
for each $c>0$, where $L_{j,l}=Ke_j\oplus Ke_l$, $y_j \in \bf K$.
\par   Let $G_2:=[x:$ for each $c>0$ there exists $n>0$ with
$|C(e_j,e_l)(x)|<c$ for each $j+l>n] \subset X\setminus (\bigcup _{c>0}
\bigcap _{n=1}^{\infty } \bigcup _{j+l>n}[x:$ $|C(e_j,e_l)(x)|>c])$,
then $\mu (G_2)=1$. Therefore, for each $x \in G\cap G_1 \cap G_2$
there exists a compact operator $D: X \to X$ such that
$|C(z_1,z_2)(x)| \le |\tilde z_2(Dz_1)|$ for each $z_1$, $z_2\in X$.
Using $\chi _e(\tilde z(x))$ (instead of $exp(2\pi i(z,x)$ from \cite{sko})
and Theorem 2.35 we find a conditional measure $\bar \mu (A|x)$.
\par  {\bf 2.38. Remark.} For $\mu : Bf(X)\to \bf R$
martingales, semi-martingales and sub-martingales are defined analogously
\S 14\cite{sko}. Then in the latter case proofs of Lemmas
1,2 and Theorems 1,2 differ slightly (with substitution of a Hilbert space on
a Banach space $X$ over $\bf K$ in all formulations). 
\section{ Quasi-invariant measures.}
\par In this section after few preliminary statements there are given
the definition of a quasi-invariant measure 
and the theorems about quasi-invariance of measures
relative to transformations of a Banach space $X$ over $\bf K$.
\par {\bf 3.1.} Let $X$ be a Banach space over $\bf K$, $(L_n:n)$
be a sequence of subspaces, $cl (\bigcup _n L_n)=X$, $L_n
\subset L_{n+1}$ for each $n$, $\mu ^j$ be probability measures,
$\mu ^2 \ll \mu ^1,$ $(\mu ^j_{L_n})$ be sequences of weak distributions,
also let there exist derivatives
$\rho _n(x)=\mu ^2_{L_n}(dx)/ \mu ^1_{L_n}(dx)$ and the following limit
$\rho (x):=\lim _{n \to \infty }\rho _n(x)$ exists.
\par {\bf Theorem.} {\it If $\mu ^j$
are real-valued and in addition $\int_X \rho (x)\mu ^1(dx)=1$ with
$\rho \in L^1(\mu ^1)$,
then this is equivalent to the following: there exists
$\rho (x)=\mu ^2(dx)/\mu ^1(dx)$
$(mod $ $\mu ^1)$.}
\par  {\bf Proof} may be obtained from the proof of theorem
1 \S 15\cite{sko} by simple modification substituting a Hilbert space
over $\bf R$ onto the Banach space $X$ over $\bf K$ (see also \S 2 above).
\par  {\bf 3.2. Theorem.} {\it Measures $\mu ^j: Bf(X)\to
\bf R$, $j=1,2$, for a Banach space $X$ over $\bf K$ are orthogonal
$\mu ^1 \perp \mu ^2$
if and only if $\rho (x)=0 $ $(mod $ $\mu ^1)$.}
\par  {\bf Proof.} In the case of real-valued $\mu ^j$
the proof differs only slightly from the proof of Theorem 2 \S 15\cite{sko}
(see also the proof in Part II for $\bf K_s$-valued measures).
\par  {\bf 3.3. Note.} For real-valued measures $\mu ^j$ on $Bf(X)$
for a Banach space $X$ over $\bf K$ (instead of a Hilbert space) 
using the above
given statements Theorems 3-6 and Corollary in \S 15\cite{sko} may be
reformulated and proved. For real-valued $\mu _j$ see 
also Kakutani Theorem 4.1 \S II.4.6\cite{dal}.
\par  {\bf 3.4. Theorem.} {\it Let $\nu $, $\mu $, $\nu _j$, $\mu _j$
be probability measures with values in $\bf R$, $X=\prod_{j=1}^{\infty }
X_j$ be a product of completely regular spaces $X_j$ with $ind (X_j)=0$
and $\mu \ll \nu $, then $\mu _j \ll \nu _j$
for each $j$ and $\prod _{j=1}^{\infty } \beta _j$ converges to $\beta $,
$\infty
> \beta > 0$, where $\beta _j:=\| (\rho _j)^{1/2} \| _{L^1(X_j, \nu _j)}$,
$\rho _j(x)=\mu _j(dx)/ \nu _j(dx)$.}
\par {\bf Proof.} The case of measures with values in $\bf R$ may be
obtained by a modification of Theorem 1 \S 16\cite{sko} (see
also the proof in Part II for $\bf K_s$-valued measures).
\par {\bf 3.5. Definition.}
Let $X$ be a Banach space over $\bf K$, $Y$ be a completely regular space with
$ind(X)=0$, $\nu : Bf(Y) \to \bf R$,
$\mu ^y: Bf(X) \to \bf R$ for each $y \in Y$, suppose $\mu ^y(A)
\in L^1(Y, \nu ,{\bf R})$ for each $A \in Bf(X)$. Then we define:
$$(i)\mbox{ } \mu (A)=\int_Y \mu _y(A) \nu (dy).$$
A measure $\mu $ is called mixed. 
We define measures $\pi ^j$ by the formula:
$$(ii)\mbox{ }\pi ^j(A\times C)=\int_C \mu ^{j,y}(A) \nu ^j(dy),$$
where $j=1,2$ and $\mu ^{y,j}$ together with $\nu ^j$ are defined
as above $\mu ^y$
and $\nu $. For $\mu ^y$ and $\nu $ with values in $\bf R$ 
definitions of $\mu $ and $\pi ^j$, Theorems 1,2 and Lemma 1  \S 18\cite{sko}
may be modified for the non-Archimedean case of the Banach space 
$X$ over $\bf K$.
\par  {\bf 3.6. Definition.} For a Banach space $X$ 
over $\bf K$ an element $a \in X$
is called an admissible shift of a measure $\mu $, 
if $\mu _a \ll \mu $, where $\mu _a(A)= \mu (S_{-a}A)$
for each $A$ in $Bf(X)$, $S_aA:=a+A$, $\rho (a,x):=\rho _{\mu }
(a,x):=\mu _a(dx)/ \mu (dx)$, $M_{\mu }:=[a \in X:$ $\mu _a \ll
\mu ]$ (see \S \S 2.1 and 2.36).
\par The proofs of the properties given below in \S 3.7
differ slightly from \S 19\cite{sko}, when a Hilbert space
is substituted on a Banach space $X$ over $\bf K$ (see the proofs
in Part II for $\bf K_s$).
\par  {\bf 3.7. Properties of $M _{\mu }$ and $\rho $ from \S 3.6.}
\par {\bf I.} {\it The set $M _{\mu }$ is a semigroup by addition, $\rho (a+b,
x)=\rho (a,x)\rho (b,x-a)$ for each $a,b \in M_{\mu }$.}
\par  {\bf II.} {\it If $a \in M_{\mu }$, $\rho (a,x) \ne 0$
$(mod$ $\mu )$, then $\mu _a \sim \mu $, $-a \in M_{\mu }$ and
$\rho (-a,x) =1/ \rho(a,x-a)$.}
\par  {\bf III. } {\it If $\nu \ll \mu $ and $\nu (dx)/ \mu (dx)=
g(x)$, then $M_{\mu } \cap M_{\nu }=M_{\mu } \cap [a:$ $\mu ([x:$
$g(x)=0,$ $g(x-a)\rho _{\mu }(a,x) \ne 0])=0]$ and $\rho _{\nu }
(a,x)=[g(x-a)/g(x)]\rho _{\mu }(a,x)$ $(mod $ $\nu )$ for $a \in
M_{\mu }\cap M_{\nu }$.}
\par  {\bf IV.} {\it If $\nu \sim \mu $, then $M_{\nu }=M_{\mu }$.}
\par  {\bf V.} Also (V-VIII) \S 19\cite{sko} may be slightly reformulated
and proved for real-valued measures.
\par {\bf 3.8. Definition and notes.} A linear operator $U$ on $X$
has a matrix representation $u(i,j)=\tilde e_i(
Ue_j)$, where $(e_j: j)$ is a non-Archimedean orthonormal basis in
in $X$, $\tilde e_i$ are vectors of the topologically adjoint space
$X^*$ under the natural embedding (see \S \S 2.2, 2.6).
For a compact operator
$U-I$ we get that there exists $\lim_{n \to \infty }
det\{ \hat r_n(U) \} \in \bf K$ which we adopt as the definition of
$det\{ U \}$, where $r_n:X\to \bf K^n$ is the projection operator with
$\hat r_n(U):=\{ u(i,j): i,j=1,...,n \}$ (see also Appendix,
use that $\lim_{i+j \to \infty } u(i,j)=0$).
\par  {\bf 3.9. Theorem.} {\it Let $X$ be a Banach space 
over $\bf K$, $\mu $ be
a probability real-valued measure and $T:X \to X$ be a compact operator
such that $Re(1-\hat \mu (z)) \to 0$ for $|\tilde
z(T^2z)| \to 0$ and $z \in X$, $\tilde z \in X'$ corresponding to $z$.
Then $M_{\mu } \subset (TX)^{\sim }$, where $Y^{\sim }$  is a
completion by $\| *\| _Y$ of a normed space $Y$.}
\par  {\bf Proof.} For each $a \in M_{\mu }$
and $b>0$ we have $J(z)=\int_X [1-[\chi _e(\tilde z(x))
+\chi _e(-\tilde z(x))]/2]
\rho (a,x) \mu (dx)$ $\le J_1(z)+J_2(z)$,
where $J_1(z):=b \int_X[1-[\chi _e(\tilde z(x))+\chi _e(-\tilde z(x))]/2]
\mu (dx)$,
$J_2(z)=2\int_{(x: \mbox{ } \rho (a,x)>b)} \mu (dx)$,
$J_1(z) \to 0$ for $|\tilde z(T^2z)| \to 0$, $J_2(z) \to 0$
for $b \to + \infty $, consequently, $\lim_{|\tilde z(T^2z)|
\to \infty } Re(\mu _a(X)- \chi _e(\tilde z(a)) \hat \mu (z))=0$.
From $|1-\hat \mu (z)|^2 \le 2Re(1-\hat \mu (z))$ and $\mu (X)=
\mu _a(X)=1$ it follows that $\chi _e(\tilde z(a))\to 1$ and
$\tilde z(a)
\to 0$ for $|\tilde z(T^2z)| \to 0$. If $|\tilde z(a)|< \epsilon $
for $|\tilde z(T^2z)| < \delta $, then $|\tilde z(a)|^2 < \epsilon ^2
|\tilde z(T^2z)|/ \delta$. Let $(e_j: j)$ be orthonormal basis in $X$ such that
$T=(T_{i,j}:$  $i,j \in {\bf N})$, $\lim_{i+j\to \infty }T_{i,j}=0,$
we denote $D_j:=\sup_k |T_{k,j}|,$ hence $|a_j|^2< \epsilon ^2
D_j^2/ \delta $,
$a \in (TX)^{\sim }$, $\| a\| _{TX}:=\sup_j D_j|a_j|< \infty $
(see Appendix).
\par  {\bf 3.10. Corollary.} {\it Let $X$ and $\mu $ be the same
as in \S 3.9, then $L\cap M_{\mu }$ is a set of the first category
in $X$ for each infinite-dimensional linear subspace $L \subset
X$ over $\bf K$.}
\par  {\bf Proof} follows from the fact that $(TB(X,0,p^j))^{\sim }
=:A_j$ is nowhere dense in $L \cap B(X,0,p^j)$ for each $j \in \bf N$,
since $A_j$ is compact in $X$ \cite{sch2}.
\par    The proof of Theorem 2 \S 19\cite{sko} can not be transferred
on non-Archimedean $X$,
there exists a compact $V$ in $X=c_0(\omega _0,{\bf K})$ and its linear span
$Y=sp_{\bf K}V$ such that an ultranorm
$p_E(x):=\inf[|a|:$ $x \in aE]$ (see Exer. 6.204 and 5.202\cite{nari})
produces from $Y$ non-separable and non-Radonian Banach space
$l^{\infty }(\omega _0,{\bf K})$
\cite{dal,roo}, where $E=co(V)$ in $Y$. 
\par {\bf Note.} Moreover, as it will be seen below
from the proofs of Theorems 3.20 and 4.2 
it follows the existense of a probability
quasi-invariant measure $\mu : Bf(X) \to {\bf R}$ 
with $M_{\mu } \supset H+G_T$, where $G_T$ 
is a compact subgroup in
$X$ with $\mu (G_T)>0$ such that
$\mu (L)=0$
for real-valued $\mu $ and each linear subspace $L$ in $X$
with $dim_{\bf K}L< \aleph _0$ (see Part II also respectively).
\par   {\bf 3.11. Definition.} For a Banach space 
$X$ over $\bf K$ and a measure $\mu : Bf(X)
\to {\bf R}$, 
$a \in X$, $\| a\| =1,$ a vector $a$
is called an admissible direction, if $a \in M_{\mu }^{\bf K}:=[z:$
$\| z\| _X=1,$ $\lambda z \in M_{\mu }$ and $\rho (\lambda z,x) \ne 0$
$(mod$ $\mu )$ (relative to $x$) and for each $\lambda \in {\bf K}]
\subset X$. Let $a \in
M_{\mu }^{\bf K}$ we denote by $L_1:={\bf K}a,$ $X_1=X\ominus L_1$, 
$\mu ^1$ and
$\tilde \mu ^1$ are the projections of $\mu $ onto $L_1$ and $X_1$
respectively, $\tilde
\mu =\mu ^1 \otimes \tilde \mu ^1$ be a measure on $Bf(X)$,
given by the the following equation
$\tilde \mu (A\times C)=\mu ^1(A)\tilde \mu ^1(C)$
on $Bf(L_1) \times Bf(X_1)$ and extended on $Bf(X)$, 
where $A\in Bf(L_1)$, $C \in Bf(X_1)$.
\par  {\bf 3.12. Theorem.} {\it For $a$ to be in $M_{\mu }^{\bf K}$
with a real-valued probability measure $\mu $ it is nesessary and sufficient
that the following conditions be satisfied (i-iv):
(i) $\mu \ll \tilde
\mu $; (ii) $\mu ^1 \ll m$, where $m$ is a real-valued Haar measure on $L_1$;
(iii) $\mu ^1(dx)/ \mu (dx):=h(x) \ne 0$ for $|a(x)|>0$;
(iv) $\mu ([x:$ $g(x- \lambda a) \ne 0,$ $g(x)=0])=0$ for
$g(x)=\mu (dx)/ \tilde \mu (dx)$ whilst $|\lambda |>0$.}
\par  {\bf Proof.} From $\tilde \mu _{\lambda a}= \mu ^1
_{\lambda a} \otimes \tilde \mu ^1$ and \S \S 3.10 II,III
sufficiency follows.
The proof of necessity is analogous to that of Theorem 1 \S 20\cite{sko}
(see \S 3.7). Indeed, from $\rho (\lambda z,x) \ne 0$
$(mod$ $\mu )$ it follows that $\mu ^1 \sim m$ due to Proposition 11
\S VII.1.9\cite{boui}.
\par  {\bf 3.13. Note.} For $\bf K_s$-valued measures Theorem
3.12 is untrue, since there are cases with $\mu $
that are not absolutely continuous relative to $\tilde \mu $.
For example, for a probability measure  $\mu =\mu _1+\mu _2\ne 0$,
$\mu _2(X)=0$ with $\mu ^1_1 $ equivalent to the Haar measure $m$ on $L_1$,
$\mu ^1_2=0$, when for $\mu _2$ all atoms $(a_{j,2})$
are points and for each $x \in L_1$ is accomplished
$\sum_{P_{L_1}(a_{j,2})=x} \mu _2(a_{j,2})=0$ and atoms  $(a_{j,1})$
of $\mu _1$ are such that $P_{L_1}(a_{j,1}) \ne P_{L_1}(a_{l,2})$
for each $j,l$. We can choose $(a_{j,s})$ such that there exists
a compact $S \subset X$ with $\mu _1(S)=0$, $\mu _2(S) \ne 0$,
but $\tilde \mu (S)=0$, since $\mu ^1_2=0$.
\par {\bf 3.14. Definition and notes.} A measure $\mu : Bf(X)\to \bf R$ 
for a Banach space $X$
over $\bf K$ is called a quasi-invariant measure 
if $M_{\mu }$ contains a $\bf K$-linear
manifold $J_{\mu }$ dense in $X$. 
\par From \S 3.7 and Definition 3.11
it follows that $J_{\mu } \subset M_{\mu }^{\bf K}$. 
\par Let $(e_j:$
$j \in {\bf N})$ be an orthonormal basis in $X$, $H=sp_{\bf K}(e_j:$ $j)$. 
We denote
$\Omega ({\bf R})=[\mu | \mu $ is a real measure with a finite total
variation on $Bf(X)$ and $H \subset J_{\mu }]$. 
\par A measure $\mu \in \Omega ({\bf R})$ that can not
be represented as a sum of two singular to each other measures from
$\Omega ({\bf R})$
is called an extremal measure. $Bf(X)$-measurable function $h:X \to \bf R$
is called invariant for $\mu \in \Omega ({\bf R})$, if
$h(a+x)=h(x)$ $\mu $-a.e. by $x$ for each $a \in H$.
\par {\bf 3.15. Lemma.} {\it If $\mu \in \Omega ({\bf R})$
and $h$ is a bounded invariant function for $\mu $, then
$h(x)=\int h(y) \mu (dy, {\sf B}^n|x)$ $\mu $-a.e. for each
$n$, where ${\sf B}^n=P_{X^n}^{-1}(Bf(X^n)),$ $H_n=sp_K
(e_1,...,e_n)$, $X^n=X\ominus H_n$, ${\sf B_n}=P_{H_n}^{-1}
(Bf(H_n))$, $P_{X^n}:X \to X^n$, $P_{H^n}:X \to H_n$ is projector
$\mu (*,{\sf B}^n|x)$ is a conditional measure.}
\par {\bf Proof.} From \S 2.37 the existence of a conditional measure follows.
Lemmas 1, 2 and corollary 1 \S 23\cite{sko} may be transferred
on a Banach space $X$ over
$\bf K$. Let us take $h_{\xi }(x)=\int_{H_n}h(x+y) \gamma _{\xi ,n}(x)
m_{H_n}(dx)$, where $\gamma _{\xi ,n}$ corresponds to $H_n$ and
$m_{H_n}$ are the same as in \S 2.30. Choosing $|\xi _l| \to \infty $
whilst $l \to \infty $ as in Lemma 3 \S 23\cite{sko} we get
the statement of Lemma 3.15.
\par  {\bf 3.16. Theorem.} {\it A measure $\mu \in \Omega ({\bf R})$
is an extremal measure if and only if $\mu $ is equivalent to
$\nu \in \sf R$, where ${\sf R}:=[\nu \in \Omega ({\sf R})|$
for each $n$ there exists $m>n$ such that $\nu (A\cap B)=\nu (A)
\nu (B)$ for each $A \in {\sf B}_n,$ $B \in {\sf B}^m]$.}
\par {\bf Proof} may be obtained by a slight modification of the proof of
Theorem 1 \S 23\cite{sko} with the help of Theorems 3.1, 3.2, 3.11,
Note 3.3 and \S \S
3.5, 2.38, since $\mu : Bf(X)\to \bf R$. In $\bf R$ the notion of
the infinite number of crossings of a segment $[\beta _1, \beta _2]
\subset \bf R$
by a sequence $(\alpha _j:$ $j \in {\bf N})$ is kept (that is,
there are $k_1<k_2<...$ such that $\alpha _{k_1} \ge \beta _2,$ $\alpha
_{k_2}\le \beta _1$,..., $\alpha
_{k_{2n-1}}\ge \beta _2,$ $ \alpha _{k_{2n}} \le \beta _1$). Then
a martingale $(\phi _n)$ with values in $\bf R$
may be considered, $B_{\beta _1, \beta _2 }:=
[x \in X:$ $(\phi _n(x):$ $n \in {\bf N})$ crosses an 
infinite number of times
a segment $(\beta _1, \beta _2)\subset \bf R$ $]$, where $\beta _2>
\beta _1 \ge 0$, $n$ is fixed, $\rho _n^m(x|y)= \mu _n^m(dx|y)/
\mu ^m(dx)$, $\mu _n^m(*|y):=\mu (*,{\sf B}_n| y)$ on ${\sf B}^m$
for $n \le m$,
$\mu ^m:=\mu |{\sf B}^m$. In addition corollaries from Lemma 3 \S 23
\cite{sko} may be slightly modified and proved for 
the Banach space $X$ over $\bf K$ due to Lemma 3.15
and Theorem 2.35. 
Then the Lebesgue measure on $\bf R^b$ must be substituted into
the Haar measure on $\bf Q_p^b$ with values in $[0, \infty )$. 
[Here the situation is
analogous to Theorem 1 \S 16\cite{sko} about absolute continuity of
real-valued product-measures
$\mu $ on a Hilbert space $X_0$ and its more general in some respects variant
Theorem II.4.1\cite{dal} on products of measurable spaces. Indeed,
$Bf(X_0)$ for $X_0$
in a topology of $\|*\|$ and a weak topology $\sigma (X_0, X_0^*)$
coincide.]
\par {\bf 3.17. Theorem. } {\it For a Banach space $X$ 
over $\bf K$ and each $\mu \in
\Omega ({\bf R})$ there exists $[\mu ^y:$ $y \in Y]\subset \sf R$ and
a measure $\nu $ on ${\sf B}^{\infty }$ such that $\mu (A)=\int_X \mu ^y(A)
\nu (dy)$ for each $A \in Bf(X)$, where ${\sf B}^{\infty }:=(\bigcap_n
\bar {\sf B}^n\cap Bf(X)$, $\bar {\sf B}^n$ is a completion of ${\sf B}^n$
by a measure $\mu $, $\nu =\mu |{\sf B}_{\infty },$ $\mu ^y(A)$ is
$Af(X, \mu )$-measurable (see the notation in \S \S 2.1, 3.14 and 3.16).}
\par  {\bf Proof.} From \S 2.37 for a conditional measure
$\mu ^y(A):=\mu (A,{\sf B}
^{\infty }|y)$ corresponding to $\mu $ relative to
${\sf B}^{\infty }$ it follows that $\mu (A)=\int \mu(A, {\sf B}^{\infty }
|y) \mu (dy)$. Let $m $ be the real-valued Haar measure on $Bf({\bf K})$, 
$m\otimes
\nu $ be a measure on $Bf({\bf K})\times {\sf B}^{\infty }$, then as in
the proof
of Theorem 2 \S 23\cite{sko} we establish that $\mu ^y_{\tau a}
\sim \mu ^y$ for $m\otimes \nu $-a.e. pairs $(\tau ,y) \in {\bf K}\times Y$,
$Y=X$. Therefore, for $\nu $-a.e. $y$ the set $S_y:=[\tau:$ $\mu ^y
_{\tau a} \sim \mu ^y,$ $ \tau \in {\bf K}]$ has $m({\bf K}\setminus S_y)=0$
and $S_y \in Af({\bf K},m)$, moreover, $S_y$ is the additive group,
consequently,
$\nu ([y:$ $S_y \ne {\bf K}])=0$. Then we may continue analogously to \S
23\cite{sko}.
\par  {\bf 3.18. Theorem.} {\it If $\mu : Bf(Y) \to \bf R$ 
is a $\sigma $-finite measure on $Bf(Y)$, $Y$
is a complete separable ultrametrizable $\bf K$-linear subspace
such that $co(S)$ is nowhere dense in $Y$ for each compact
$S \subset Y$, where $\bf K$ is an infinite non-discrete
non-Archimedean field with a multiplicative ultranorm
$|*|_{\bf K}$. Then from $J_{\mu }=Y$
it follows that $\mu =0$.}
\par  {\bf Proof.} Since $\mu $ is $\sigma $-finite, then
there are $(Y_j:$ $j \in H)\subset Bf(Y)$
such that $Y=\bigcup_{j\in H}
Y_j$ and $0< \| \mu |Bf(Y_j) \| \le 1$ for each $j$, where
$H \subset
\bf N$, $Y_j \cap Y_l= \emptyset $ for each $j \ne l$. If $card(
H)=\aleph _0$, then we define a function $f(x)=1/[2^j|\mu |(Y_j)] $
for real-valued $\mu $. 
Then we define a measure $\nu (A)=\int_A f(x) \mu (dx)$, $A
\in Bf(Y)$. Therefore, $|\nu |(Y) \le 1$ 
and $J_{\nu }=Y$, since $f \in L^1(Y, \mu ,{\bf R})$. 
Hence it is sufficient to consider
$\mu $ with $\| \mu \| \le 1$ and  $|\mu |(Y)=1$.
For each $n \in \bf N$ in view of the Radonian property of
$Y$ there exists a compact $X_n \subset Y$ such that $|\mu |(
Y\setminus X_n)<1/n$. In $Y$ there is a countable everywhere
dense subset $(x_j:$ $j \in {\bf N})$, hence $Y=\bigcup _{j
\in \bf N}B(Y,x_j,r_l)$
for each $r_l>0$, where $B(Y,x,r_l)=[y \in Y:$ $d(x,y) \le r_l]$,
$d$ is an ultrametric in $Y$, i.e. $d(x,z) \le max(d(x,y),$ $d(y,z))$,
$d(x,z)=d(z,x)$, $d(x,x)=0$, $d(x,y)>0$ for $x \ne y$ for each $x,y,z
\in Y$. Therefore, for each $r_l=1/l,$ $l\in \bf N$ there exists
$k(l)
\in \bf N$ such that $|\mu |(X_{n,l})>1-2^{-l-n}$, 
where $X_{n,l}:=\bigcup_{j=1}
^{k(l)}B(Y,x_j,r_l)$, consequently, $| \mu |(Y\setminus X_n) \le
2^{-n}$ for $X_n:=
\bigcap_{l=1}^{\infty }X_{n,l}$. The subsets $X_n$ are compact,
since $X_n$ are closed in $Y$ and the metric $d$ on $X_n$ is
completely bounded and $Y$ is complete (see Theorems 3.1.2 and 4.3.29
\cite{eng}). Then $0<|\mu |(X)
\le 1$ for $|\mu |(Y\setminus X)=0$
for $X:=sp_{\bf K}(\bigcup_{n=1}^{\infty } X_n).$
\par  The sets $\tilde Y_n=co(Y_n)$ are nowhere dense in $Y$ for
$Y_n=\bigcup_{l=1}^nX_l$, consequently, $sp_{\bf K}Y_n$ are nowhere dense
in $Y$. Moreover, $(Y\setminus \bigcup_{n=1}^{\infty }Y_n) \ne \emptyset $
is dense in $Y$ due to the Baire category theorem (see 3.9.3 and 4.3.26
\cite{eng}).
Therefore, $y+X \subset Y\setminus X$ for $y \in Y\setminus X$ and from
$J_{\mu }
=Y$ it follows that $|\mu |(X)=0$, since $|\mu |
(y+X)=0$ (see \S \S 2.38 and 3.14 above). 
Hence we get the contradiction, consequently, $\mu =0$.
\par  {\bf 3.19. Corollary.} {\it If $Y$ is a Banach space or a complete
countably-ultranormable infinite-dimensional over $\bf K$ space,
$\mu : Bf(Y) \to \bf R$ and $J_{\mu }=Y$, then $\mu =0$.}
\par  {\bf Proof.} The space $Y$ is evidently complete and
ultrametrizable, since its topology is given by a countable
family of ultranorms. Moreover, $co(S)$ is nowhere dense in $Y$
for each compact $S$ in $Y$, since $co(S)=cl(S_{bc})$ is compact
in $Y$ and does not contain in itself any open subset in $Y$ due to
(5.7.5)\cite{nari}
(for $Y$ over $\bf R$ and real-valued measures see Theorem
4 \S V.5.3\cite{gevi}).
\par {\bf 3.20. Theorem.} { \it Let $X$ be a separable Banach space over
a locally compact infinite field $\bf K$ with a nontrivial valuation
such that either ${\bf K}\supset \bf Q_p$ or $char ({\bf K})=p>0$. 
Then there are probability measures 
$\mu $ on $X$ with values in ${\bf R}$
such that $\mu $ are quasi-invariant 
relative to a dense $\bf K$-linear subspace
$J_{\mu }$}.
\par {\bf Proof.} Let $S(j,n):=p^{j}B({\bf K},0,1)
\setminus p^{j+1} B({\bf K},0,1)$
for $j\in \bf Z$ and $j\le n$,
$S(n,n):=p^nB({\bf K},0,1)$, $w$ be the Haar measure on $\bf K$
considered as the additive group (see \cite{hew,roo}) with
values in ${\bf R}$.
Then for each $c>0$ and $n \in {\bf N}$ there are measures
$m$ on $Bf({\bf K})$ such that
$m(dx)=f(x)v(dx)$, $\mid f(x)\mid >0$ for each $x \in \bf K$
and $\mid m(p^{n}B({\bf K},0,1))-1\mid <c$,
$m({\bf K})=1$, $\mid m\mid (E)\le 1$ for each $E\in Bf({\bf K})$ ,
where $v=w$, $v(B({\bf K},0,1))=1$. Moreover,
we can choose $f$ such that a density $m_{a}(dx)/m(dx)=:d(m;a,x)$
be continuous by
$(a,x)\in \bf K^{2}$  and for each $c'>0$,  $x$ and
$\mid a\mid \le p^{-n}$ : $\mid d(m;a,x)-1\mid <c'.$
For this we can define $f$ for $v=w$, for example, 
to be at the beginning locally constant
such that $f(x)=a(j,n)$ for $x\in S(j,n)$, where $a(j,n)=r^{n(j-n)}
(1-r^{-n})(1-1/p)p^{-n}$ for $j<n$, $a(n,n)=(1-r^{-2n})p^{-n}$
and $ m(E):=\sum \{ a(j,n) v(p^{n-j}(E\cap S(j,n))/v(p^{n-j}$
$S(j,n)): j \in {\bf Z}, j\le n \} $.
Then we can take $g(x)=f(x)+h(x)$ and $y(dx):=g(x)v(dx)$,
and a continuous $h(x):K\to \bf R$,
with $\sup \{ \mid h(x)/f(x)\mid : x\in K \} \le c"$ and
$0<c"\le 1/p^n$.
\par  Let $\{ m(j;dx) \}$ be a family of measures on
$\bf K$ with the corresponding sequence $\{ k(j) \}$ such that
$k(j)\le k(j+1)$ for each $j$ and
$\lim_{i\to \infty }k(i)= \infty $, where
$m(j;dx)$ corresponds to the partition $[S(i,k(j))]$. The Banach space 
$X$ is isomorphic with $c_0
(\omega _0,{\bf K})$ \cite{roo}. It has the orthonormal basis
$\{ e_j:j=1,2,... \} $ and the projectors
$P_jx=(x(1),...,x(j))$ onto $\bf K^j$, where $x=x(1)e_1+x(2)e_2+...$.
Then there exists a cylindrical measure $\mu $ generated by a consistent family
of measures $y(j,B)=b(j,E)$ for $B=P^{-1}_jE$ and $E\in Bf({\bf K}^j)$ 
\cite{boui,dal}
where $b(j,dz)=\otimes [m(j;dz(i)): i=1,...,j], z=(z(1),...,z(j))$.
Let $L:=L(t,t(1),...,t(l);l):=\{ x: x \in X \mbox{ and } \mid x(i)\mid
\le p^{a}, a=-t-t(i) \mbox{ for } i=1,...,l, \mbox{ and } a=-k(j)
\mbox{ for } j>l \} $, then $L$ is compact in $X$, since $X$ is
Lindel\"of and $L$ is sequentially compact \cite{eng}.
Therefore, for each $c>0$ there exists
$L$ such that $\mid \mu (X\setminus L)\mid <c$, since there is $l\in {\bf N}$
with $\mid 1-\prod [m(j;p^{k(j)}B(K,0,1)): j>l]\mid <c/2$ (or due to
the choice of $a(j,n)$).
\par In view of the Prohorov theorem \S IX.4.2\cite{boui} for real
measures and due to Lemma
2.3 $\mu $ has the countably-additive extension on $Bf(X)$, 
consequently, also on the complete $\sigma $-field
$Af(X,\mu )$ and $\mu $ is the Radon measure.
\par Let $z' \in sp_{\bf K} \{ e_j:j=1,2,... \} $
and $z"= \{ z(j): z(j)=0
\mbox{ for } j\le l \mbox{ and } z(j)\in S(n,n), j=1,2,...,
n=k(j) \} $, $l \in {\bf N}$, $z=z'+z"$.
In view of the Kakutani theorem \cite{dal}
and also I.1.4, II.4.1\cite{roo}
there are $m(j;dz(j)$ such that $\rho _{\mu }(z,x)=\prod \{ d(j;z(j),
x(j)): j=1,2,...\} =
\mu _z(dx)/\mu (dx)\in L^1(X,\mu ,{\bf R})$ 
for each such $z$ and $x\in X$, where $d(j;*,*)=
d(m(j;*),*,*)$ and $\mu _z(X)=\mu (X)=1$.
\par When $char({\bf K})=0$, then in this proof it is also possible
to consider $S(j,n):=\pi ^jB({\bf K},0,1)\setminus \pi ^{j+1}B({\bf K},0,1)$,
where $\pi \in \bf K$, $p^{-1}\le |\pi |<1$, $\pi $ is the generator of the
valuation group of $\bf K$.
\par {\bf 3.21. Note.} For a given $m=w$ new suitable
measures may be constructed, if to use images of measures
$m^{g}(E)=m(g^{-1}(E))$ such that for a diffeomorphism
$g\in Diff^1({\bf K})$ (see \S A.3) we have $m^{{g}^{-1}}(dx)/m(dx)=
\mid (g'(g^{-1}(x))\mid _{\bf K}$, where $|*| _{\bf K}= mod_{\bf K}(*)$ 
is the modular function of the field $\bf K$ associated with
the  Haar meassure on $\bf K$, at the same time  
$|*| _{\bf K}$ is the multiplicative norm in $\bf K$
consistent with its uniformity \cite{wei}.
Indeed, for $\bf K$ and $X=\bf K^j$ with $j \in \bf N$ and 
the Haar measure $v=w$ or $v=w'$ on $X$, $v_X:=v$
with values either in $\bf R$ or in $\bf K_s$ for $s\ne p$
and for a function $f \in L^1(X,v,{\bf R})$ or $f \in
L(X,v,{\bf K_s})$, respectively, due to Lemma 4 and Theorem 4
\S I.2\cite{wei}, Theorems 9.2 and A.7\cite{sch1}
we have: $\int_{g(A)} f(x)v(dx)$ $=\int_Af(g(y)) |g'(y)|_{\bf K} v(dy)$, where
$mod_{\bf K}(\lambda )v(dx):=v(\lambda dx)$, $\lambda \in \bf K$.
For a construction of new measures images of $v$ may be used, for example,
for analytic $h$ on $\bf K$ (see 43.1 \cite{sch1}) with $g=Ph$,
where $P$ is the antidifferentiation such that $g'=h$. Defining it
at first locally we can then extend it on all $\bf K$.
\par Henceforward, quasi-invariant measure $\mu $ on
$Bf(c_0(\omega _0,{\bf K}))$
constructed with the help of projective limits or sequences
of weak distributions 
of probability measures $(\mu _{H(n)}:$ $n)$ are considered,
for example, as in Theorem 3.20 such that
\par  (i) $\mu _{H(n)}(dx)=f_{H(n)}(x)v_{H(n)}(dx)$, $dim_{\bf K}H(n)
=m(n)< \aleph _0$ for each $n \in \bf N$, where $f_{H(n)} \in
L^1(H(n),v_{H(n)},{\bf R})$,
$H(n) \subset H(n+1) \subset ...$, $cl(\bigcup_nH(n)=c_0(
\omega _0,{\bf K})$, if it is not specified in another manner.
\par For real-valued probability quasi-invariant measures
for a sequence of a weak distribution
$(\mu _{H(n)})$ generated by $\mu $ in view of Proposition
11 \S VII.1.9
\cite{boui}, Lemma 2.3 and Definition 3.14 this condition is satisfied,
since $cl(J_{\mu })=
c_0(\omega _0,{\bf K})$.
\par  As will be seen below such measures $\mu $ are quasi-invariant
relative to families of the cardinality ${\sf c}=card ({\bf R})$
of linear and non-linear transformations $U: X\to X$. Moreover, for each
$V$ open in $X$  we have $\mu (V) >0$, when $f_{H(n)}(x)
>0$ for each $n \in \bf N$ and $x \in H(n)$.
\par  Let $\mu $ be a probability 
quasi-invariant measure satisfying (i) and  $(e_j:$ $j)$
be an orthonormal basis in $ M_{\mu }$, 
$H(n):=sp_{\bf K}(e_1,...,e_n)$, we denote by
\par  $\hat \rho
_{\mu }(a,x)=\hat \rho (a,x)=\lim_{n \to \infty } \rho ^n(P_na,P_nx)$,
$\rho ^n(P_na,P_nx):=f_{H(n)}(P_n(x-a))/f_{H(n)}(P_nx)$ for each
$a$ and $x$ for which this limit exists and $\hat \rho (a,x)=0$
in the contrary case, where $P_n: X\to H(n)$ are chosen 
consistent projectors. Let
$\rho (a,x)=\hat \rho (a,x)$, if $\int \hat \rho (a,x) \mu (dx)=1$, 
$\rho (a,x)$ is not
defined when $\hat \rho (a,x) \ne 1$. 
If for some another basis $(\tilde e_j:$ $j)$ and $\tilde \rho $
is accomplished
\par  (ii) $\mu  (S)=1$, $S:=
\bigcap_{a \in M_{\mu }}[x:$ $\rho
(a,x)=\tilde \rho (a,x)])$,
then $\rho (a,x) $ is called regularly dependent from a basis.
\par  {\bf 3.22. Lemma.} {\it Let $\mu $ be a probability measure,
$\mu : Bf(X)\to {\bf R}$, $X$ be a Banach space over $\bf K$,
suppose that for each basis $(\tilde e_j:$ $j)$ in
$M_{\mu }$ a quasi-invariance factor $\tilde \rho $ satisfies the following
conditions: 
\par $(1)$ if $\tilde \rho (a_j,x)$, $j=1,...,N$, are defined for
a given $x \in X$ and for each $\lambda _j \in \bf K$ then a function
$\tilde \rho (\sum_{j=1}^n \lambda _ja_j,x)$ is continuous by
$\lambda _j$, $j=1,...,N$; 
\par $(2)$ there exists an increasing sequence of subspaces
$H(n) \subset M_{\mu }$, $cl(\bigcup _nH(n))=X$, with projectors
$P_n: X \to H(n)$, $B \in Bf(X)$, $\mu (B)=0$
such that $\lim_{n\to \infty }\tilde \rho (P_na,x)=\tilde \rho (a,x)$
for each $a \in M_{\mu }$ and $x \notin B$ for which is defined
$\rho (a,x)$. Then $\rho (a,x)$ depends regularly from the basis.}
\par  {\bf Proof.} There exists a subset $S$ dense in each
$H(n)$, hence $\mu (B')=0$  for
$B'=\bigcup_{a \in S}[x:$ $\rho (a,x)\ne \tilde \rho (a,x)$].
From (1) it follows that $\tilde \rho (a,x)=\rho (a,x)$ on each $H(n)$
for $x \notin B'$. From $sp_{\bf K}S \supset H(n)$ and (2) it follows that
$\rho (a,x)=\tilde \rho (a,x)$ for each $a \in M_{\mu }$ and $x \in
X\setminus (B'\cup B)$, consequently, Condition 3.21(ii) is satisfied.
\par {\bf 3.23. Lemma.} {\it If a probability quasi-invariant measure 
$\mu :Bf(X)\to \bf R$
satisfies Condition 3.21(i), then there exists a compact operator
$T: X\to X$ such that $M_{\mu } \subset (TX)^{\sim }$, where $X$ is
the Banach space over $\bf K$.}
\par {\bf Proof.} This follows from Theorem 3.9 (see also Part II).
\par  {\bf 3.24.} Let $X$ be a Banach space over $\bf K$,
$|*|_{\bf K}=mod_{\bf K}(*)$, $U: X\to X$ be an invertible linear operator,
$\mu : Bf(X) \to \bf R$ be a probability quasi-invariant measure.
\par The uniform convergence of a (transfinite) sequence of functions on
$Af(V,\nu )$-compact subsets of a topological space
$V$ is called the Egorov condition, where $\nu $ is a measure on  $V$.
\par {\bf Theorem.} {\it Let pairs $(x-Ux,x)$
and $(x-U^{-1},x)$ be in $dom(\tilde \rho (a,x))$, where $dom(f)$ denotes
a domain of a function $f$, $\tilde \rho (x-Ux,x)>0$, $\tilde
\rho (x-U^{-1}x,x)>0$ $(mod $ $\mu )$.
Then $\nu \sim \mu $ and
$$(i) \mbox{ }\nu (dx)/ \mu (dx)=|det(U)|_{\bf K}\tilde \rho (x-U^{-1}x,x),$$
if $\rho $ depends regularly from the base, then $\tilde \rho $ may be
substituted by $\rho $ in formula $(i)$, where $\nu (A):=
\mu (U^{-1}A)$ for each $A \in Bf(X)$.}
\par  {\bf Proof.} In view of Theorem 3.9 (and Lemma 3.23) there
exists a compact operator $T: X\to X$ such that $M_{\mu }
\subset (TX)^{\sim }$, consequently, $(U-I)$ is a compact operator,
where $I$ is the identity operator. From the invertibility of
$U$ it follows that $(U^{-1}-I)$ is also compact, moreover, there
exists $det(U) \in \bf K$. Let
$g$ be a continuous bounded function, $g: \tilde H(n) \to \bf R$, 
whence $\int_X\phi (x)\nu (dx)$ $=\int_{\tilde H(n)}
g(x)[f_{H(n)}(U^{-1}x)/f_{\tilde H(n)}(x)]|det(U_n)|_{\bf K} \mu _{\tilde
H(n)}(dx)$, for $\phi (x)=g(\tilde P_nx)$, where subspaces exist such
that $\tilde H(n) \subset X$, $(U^{-1}-I)\tilde H(n) \subset \tilde H(n)$,
$cl(\bigcup_n \tilde H(n))=X$, $U_n:=\hat r_n(U)$, $r_n=\tilde P_n:
X \to \tilde H(n)$ (see \S \S 3.8 and 3.21), $\tilde H(n) \subset \tilde
H(n+1)\subset ...$ due to compactness of $(U-I)$. In view of the Fatou
theorem $J_m \ge J_{m, \rho }$, where $J_m:=\int_Xg(\tilde P_mx)
\nu (dx)$ and $J_{m,\rho }:=\int_Xg(\tilde P_mx)\tilde \rho (x-U^{-1}x,x)
|det(U)|_{\bf K} \mu (dx)$. Indeed, there exists $n_0$ such that
$|u(i,j)-\delta _{i,j}| \le 1/p$ for each $i$ and $j>n_0$,
consequently, $|det(U_n)|_{\bf K}=|det(U)|_{\bf K}$ for each
$n>n_0$, where $\delta _{i,j}$ is the Kronecker delta-symbol.
Then due to the Fubini Theorem 
there exists $\lim_{n\to \infty }[\mu _{\tilde
H(n)}(d\tilde P_nx)/\nu _{\tilde H(n)}(d\tilde P_nx)]$
$=\mu (dx)/ \nu (dx)$ $(mod $ $\nu )$. Further analogously to the proofs
of Theorems 1 and 2 \S 25\cite{sko}.
\par {\bf 3.25.} Let $X$ be a Banach space over $\bf K$,
$|*|_{\bf K}=mod_{\bf K}(*)$ with a probability quasi-invariant 
measure $\mu : Bf(X)\to \bf R$,
also let $U$ fulfils the following conditions:
$$(i)\mbox{ }U(x) \mbox{ and } U^{-1}(x)
\in C^1(X,X) \mbox{ (see }\mbox{ \S A.3)};$$
$$ (ii)\mbox{ } (U'(x)-I)\mbox{ is compact for each }x\in X;$$
$$(iii)\mbox{ }(x-U^{-1}(x)) \mbox{ and }(x-U(x)) \in J_{\mu }
\mbox{ for }\mu -\mbox{a.e. }x\in X;$$
$$ (iv)\mbox{ for }\mu \mbox{-a.e. }x \mbox{ pairs }
(x-U(x);x) \mbox{ and }(x-U^{-1}(x);x)$$
are contained in a domain of $\rho (z,x)$ such that
$\rho (x-U^{-1}(x),x)\ne 0$, $\rho (x-U(x),x)\ne 0$ $(mod$ $\mu )$;
$$ (v)\mbox{ }\mu (S')=1,$$
where $S':=([x: \mbox{ }\rho (z,x) \mbox{ is defined and continuous by }
z\in L ])$
for each finite-dimensional $L\subset J_{\mu }$;
$$ (vi)\mbox{ there exists }S \mbox{ with }\mu (S)=0
\mbox{ and for each }$$ $x\in X\setminus S$
and for each $z$ for which there exists $\rho (z,x)$ satisfying the
following condition:
$\lim_{n \to \infty }\rho (P_nz,x)=\rho (z,x)$ and the convergence
is uniform for each finite-dimensional $L\subset J_{\mu }$ by $z$ in
$L\cap [x\in J_{\mu }:$ $\mid x\mid \le c]$, where $c>0$,
$P_n: X\to H(n)$  are projectors onto
finite-dimensional subspaces $H(n)$ over $\bf K$ such that
$H(n)\subset H(n+1)$ for each $n \in \bf N$ and $cl \cup \{ H(n): n \}=X$;
$$ (vii)\mbox{ there exists }n \mbox{ for which for all }j>n
\mbox{ and }
x\in X  \mbox{ mappings}$$
$V(j,x):=x+P_j(U^{-1}(x)-x)$ and $U(j,x):=
x+P_j(U(x)-x)$ are invertible and
$\lim_j \mid \det U'(j,x)\mid =\mid \det U'(x)\mid $,
$\lim_j \mid \det V'(j,x)
\mid =1/ \mid \det U'(x)\mid $.
\par {\bf Theorem.} {\it The measure 
$\nu (A):=\mu (U^{-1}(A))$ is equivalent to $\mu $ and
$$ (i)\mbox{ }\nu (dx)/ \mu (dx)=\mid \det U'(U^{-1}(x))\mid _{\bf K}
\rho (x-U^{-1}(x),x).$$}
\par {\bf Proof.} I. Let at first $U$ be linear.
In general, for a linear operator $U$ with
compact $B=U-I$ there is the following decomposition $U=SCDE$, where
$C^t$ and $E$ are upper triangular infinite matrices,
$D=diag\{ d(j):
j\in {\bf N} \}$, operators $C-I$, $D-I$ and $E-I$
are compact in the corresponding orthonormal basis $\{ e_j: j \}$ in $X$, $S$ transposes
a finite number of vectors in orthonormal basis
(see \S A.2). Moreover, there are $det(C)=det(E)=1$,
$det(U)=det(D)\ne 0$.
\par   II. Let $V_n$ be a diagonal (or upper or lower triangular)
operator on $X$ such that $(V_n-I)(X)=P_nL$, where $dim_{\bf K}L=k<\aleph _0$,
$\lim_{n \to \infty }\| V_n-U_1\| =0$, $U_1$ is a diagonal (or
lower or upper triangular) operator, there exists $n_0$
such that $\| e_j-P_ne_j\|\le
1/p$ for $(e_j)$ in $L$, consequently, $\| \sum_j \lambda _jP_ne_j\|
=\max_j |\lambda _j|$ for $\lambda _j \in \bf
K$ and $dim_{\bf K}P_nL=k$
for $n>n_0$, in addition,  $\lim_{n \to \infty }\sup_{x \in L, \| x\| \le 1}
\| x-P_nx\| =0$. Then $\lim_n \tilde P_n^{-1}(x-V_n^{-1}x)=x-
U_1^{-1}x$, where $\tilde P_nx=P_nx\in P_nL$ for each $x \in L$.
Due to Conditions (vi, vii) we get $\lim_n \rho (x-V_n^{-1}x,x)=
\rho (x-U_1^{-1}x,x)$ $(mod$ $\mu )$.
\par   From the Fatou theorem and \S 3.21 it follows that,
$J_1 \ge J_{1, \rho }$,
where $J_1=\int_Xf(U_1x) \mu (dx)$, $J_{1,\rho }:=\int_Xf(x)\rho (x-
U_1^{-1}x,x)|det U_1|_{\bf K}\mu (dx)$ for continuous bounded function
$f:X\to [0, \infty )$. Analogously for $U_1^{-1}$
instead of $U_1$. Using instead of $f$ the function $\Phi (U_1^{-1}
x):=f(x)\rho _{\mu }(x-U_1^{-1}x,x)$ and Properties 3.7 we get
that $\rho _{\mu }(U_1x-x,U_1x)\rho _{\mu }(x-U_1x,x)=1$ $(mod$ $\mu )$.
Therefore, for $U=U_1U_2$ with diagonal $U_1$ and upper triangular $U_2$
and lower triangular $U_3$ operators with finite-dimensional over
$\bf K$ subspaces $(U_j-I)X$, $j=1,2,3$, the following
equation is accomplished
$$\int_Xf(Ux)\mu (dx)=\int_X f(x)\rho _{\mu }(x-U^{-1}x,x)
|det U|_{\bf K} \mu (dx).$$ If $(S^{-1}U-I)X=L$,
then from the
decomposition given in (I) $U=SU_2U_1U_3,$ we have $(U_j-I)X=L$, 
$j=1,2,3$ due to formulas from \S A.1, since
corresponding non-major minors are equal to zero.
\par  (III). If $U$ is an arbitrary linear operator satisfying
the conditions of this theorem, then from (iv-vi) and (I, II)
for each continuous bounded function $f: X \to [0, \infty )$
we have $J \ge J{\rho }$, where $J:=\int_Xf(U(x))\mu (dx)$ and $J_{\rho }:=
\int_Xf(x)\rho _{\mu }(x-U^{-1}(x),x)|det U|_{\bf K}\mu (dx)$. 
Analogously for
$U^{-1}$, moreover, $\rho (x-U^{-1}(x),x)|det U|_{\bf K}=:h(x)\in L^1(X,
\mu ,{\bf R})$, $h(x)>0$ $(mod$ $\mu )$, since there exists
$det U$.
\par  (IV). Now let $U$ be linear and $(U-I)(X)=L$, $dim_{\bf K}L=k<\aleph _0$,
$L\subset J_{\mu }$.
Suppose $U$ is polygonal, which means that there exists a partition
$X= \cup \{ Y(i): i=1,...,l \}$, $U(x)=a(i)+
V(i)x$ for $x \in Y(i)$, where $Y(i)$ are closed subsets,
$Int Y(i)\cap  Int Y(j)=\emptyset $ for each $i\ne j$, $a(i) \in X$
and $V(i)$ are linear operators.
Then $U^{-1}$ is also polygonal, $U'(x)=V(j)$ for $x \in Y(j)$ and
$\int_X f(a(i)+V(i)x)\mu (dx)=\int f(a(i)+x)\rho _{\mu }(x-V^{-1}(i)x,x)$
$\times |det(V(i))|_{\bf K} \mu (dx)$
for each real Borel measurable function $f$ 
and each $i$. From $a(j) \in M_{\mu }$
and \S 3.7 we get $\int_Xf(a(j)+V(j)x) \mu (dx)$ 
$=\int_X f(x)\rho (x-V(j)^{-1}
(x-a(j)),x) |det V(j)|_{\bf K} \mu (dx)$. Let $H_{k,j}:=[x\in X:$ $V(k)^{-1}x
=V(j)^{-1}x]$, assume without loss of generality that $V(k) \ne V(j)$ or
$a(k) \ne a(j)$ for each
$k \ne j$, since $Y(k) \ne Y(j)$ (otherwise they may be united).
Therefore, $H_{k,j} \ne X$. If $\mu (H_{k,j})>0$,
 then from $X \ominus H_{k,j} \supset \bf K$ it follows
that $M_{\mu }
\subset H_{k,j}$, but $cl(H_{k,j})=H_{k,j}$ and $cl(M_{\mu })=X$. This
contradiction means that $\mu (A)=0$,
where $A=[x:$ $V(k)^{-1}(x-a(k))=V(j)^{-1}(x-a(j))]$. Then
$\int_Xf(U(x)) \mu (dx)$ $=\int_Xf(x)\rho (x-U^{-1}(x),x)$
$|det U'(x)|_{\bf K}^{-1}\mu (dx)$.
\par  (V). In view of the Hahn-Banach theorem for $\bf K$
\cite{nari,roo} there are linear continuous functionals $\tilde e_j$
such that there exists orthonormal basis $( e_j:$ $j) $ in $X$ with
$\tilde e_i(e_j)=\delta _{i,j} $. Let
$s(i,j;x):=\tilde e_i(U(j,x)-x)$ and $s(i;x)=\tilde e_i(U(x)-x)
=\lim_j s(i,j;x)$
($\lim$ is taken in $C(X,{\bf K})$), consequently,
$\det U'(j,x)=\det(ds(i,j;x)e(k): i,k=1,...,j)$. Then for the construction
of the sequence $\{ U(j,*) : j \}$ it is sufficient to construct a sequence
of polygonal functions $\{ a(i,j;x)\} $, that is $a(i,j;x)=l_k(i,j)(x)+a_k$
for $x \in Y(k)$, where $l_k(i,j)$ are linear functionals,  $a_k\in \bf K$,
$Y(k)$ are closed in $X$, $Int(Y(j))\cap Int(Y(k))=\emptyset $ for each
$k \ne j$, $\bigcup_{k=1}^m Y(k)=X$, $m<\aleph _0$. For each $c>0$
there exists $V_c \subset X$ with  $\mu (X\setminus V_c)<c$, 
the functions $s(i,j;x)$ and
$(\bar \Phi ^1 s(i,j; *))(x,e(k),t)$
are equiuniformly continuous (by $x \in V_c$ and by $i, j, k \in {\bf N}$)
on $V_c$. Choosing $c=c(n)=1/n$ 
and using $\delta $-nets in $V_c$ we get a sequence of polygonal
mappings $(W_n:$ $n)$ converging by its matrix elements in
$L^1(X,\mu ,{\bf R})$,
from Condition (i) follows that it may be chosen equicontinuous for
matrix elements $s(i,j;x)$, $ds(i,j;x)$
and $s(i,P_jx)$ by $i, j$ (the same is true for $U^{-1}$).
\par  Indeed, for $V_c$ there is $\delta >0$ such that $|s(i,j;x')-
s(i,j;x)|<c$ and $|(\bar \Phi ^1s(i,j;*))(x,e(k),t)- (\bar \Phi ^1s(i,j;*))
(x',e(k),t')|< c$ (see Theorem 2.11.I \cite{lu3}) for each $x, x' \in V_c$ with
$|x-x'|< \delta $, $|t-t'|
\le 1$ and $i, j, k \in \bf N$. Then we can choose in $V_c$ a finite
$\delta $-net $x_1,...,x_r$ and define $l_q(i,j;x)=s(i,j;x_q)+
(ds(i,j;x_q))(x-x_q)$ and apply the non-Archimedean variant of the Taylor
theorem (see \S 29.4\cite{sch1} and Theorem A.5).
\par  Then calculating integrals as above for $W_n$ with functions
$f$, using the Fatou theorem we get
the inequalities analogous to written in
(III) for $J$ and $J_{\rho }$ of the general form. From
$\nu (dx)/\mu (dx)>0$ $(mod$ $\mu )$ 
we get the statement of the theorem (see also Theorem 1 \S 26\cite{sko}
in the classical case), since 
$$\int_Xf(A(j,x))\mu (dx)=\int_X f(x) \rho _{\mu }(P_j(x-A(x)),x)
|det(A(j,x)')|_{\bf K}^{-1}\mu (dx),$$
 where $A=U$, $A(j,x)=U(j,x)$ or
$A=U^{-1}$, $A(j,x)=V(j,x)$.
\par   {\bf 3.26. Note.} For a linear invertible operator $U$
Condition 3.25(i) is satisfied automatically, (ii) follows from (iii),
hence from (iii) also follows (vii) (see the proof of Theorem 3.24).
\par   {\bf 3.27. Examples.} Let $X$ be a Banach space 
over the field $\bf K$ with the valuation group 
$\Gamma _{\bf K}=\Gamma _{\bf Q_p}$.
We consider a diagonal compact operator
$T=diag(t_j:$ $j\in {\bf N})$ in a fixed orthonormal basis
$(e_j)$ in $X$ such that $ker T:=T^{-1}0=\{ 0 \}$. Let $\nu _j(dx_j)=
C(\xi _j) exp(-|(x_j-x^0_j)/\xi _j|_p^q) v(dx_j)$ for the Haar measure
$v: Bf({\bf K}) \to [0, \infty )$. 
We choose constant functions $C(\xi _j)$ and $C'(\xi _j)$
such that $\nu _j$ be a probability measure, where
$x^0=(x^0_j:$ $j)
\in X$, $x=(x_j:$ $j) \in X$, $x_j \in \bf K$. Particularly, for
$q=2:\quad $ $\int_{\bf K}[(\chi _e(x_j-x^0_j)-\chi _e(-x_j+x^0_j))/(2i)]
\nu _j(dx_j)=0$, since $\nu _j$ is symmetric and $[\chi _e(x)-\chi _e(-x)]$ 
is the odd function, where
$\chi _e $ is the same character as in \S 2.6. Then
$\int_{\bf K}|x_j-x^0_j|^2 \nu _j(dx_j)=
|\xi _j|^2$ due to the formula of changing variables in \S 3.21
and differentiation by the real parameter $b$ of the following
integral  $I(b)=\int_{\bf K}exp(-|x|
_p^2b^2)v(dx)$, then we take $1/|\xi |=b \in \Gamma _{\bf K}$. Therefore,
$x_j^0$ and $|\xi _j|^2$ have in some respect a meaning looking like
the mean value and the variance.
\par With the help of products
$\bigotimes_{j=1}^{\infty }\nu _j(dx_j)$ 
as in \S 3.20 we can construct a probability
quasi-invariant measure $\mu ^T$ on $X$ with values in $[0,1]$ or $\bf C_s$, 
since $cl(TX)$ is compact in $X$ and $sp_{\bf K}(e_j:$ $j)=:H
\subset J_{\mu }$. From $\bigcap
_{\lambda \in B({\bf K},0,1)\setminus 0}cl(\lambda TX)=\{ 0\}$ we may infer
that for each $c>0$ there exists a compact $V_c(\lambda )
\subset X$ such that $\mu (X\setminus V_c(\lambda ))<c$ 
and $\bigcap_{\lambda \ne 0}
V_c(\lambda )=\{ 0\}$, consequently, $\lim_{|\lambda | \to 0}
\int_Xf(x) \mu ^{\lambda T}(dx)=f(0)=\delta _0(f)$, hence
$\mu ^{\lambda T}$
is weakly converging to $\delta _0$ whilst $|\lambda | \to 0$ for the
space of bounded continuous functions $f: X \to \bf R$.
\par From Theorem
3.4 we conclude that from $\sum_{j=1}^{\infty }|y_j/\xi _j|^q_p < \infty $
it follows $y \in J_{\mu ^T}$. Then for a linear transformation
$U: X\to X$ from $\sum_j|\tilde e_j(x-U(x))/\xi _j|_p^2 < \infty $
it follows that $x-U(x) \in J_{\mu }$ and a pair $(x-U(x),x) \in dom(\rho
(a,z))$. Moreover, for $\rho $ corresponding to $\mu ^T$ conditions
(v) and (vi) in \S 3.25 are satisfied. Therefore, for such $y$ and
$S \in Af(X,\mu )$ 
a quantity $|\mu (ty+S)- \mu (S)|$ is of order of smallness $|t|^q$
whilst $t \to 0$, hence they are pseudo-differentiable of order $b$
for $0<Re(b)<2$ (see also \S 4 below).
\par Analogs of real-valued Wiener measures on spaces of functions 
from a non-Archimedean Banach space into $\bf C$ are given in
\cite{byvo}, that are quasi-invariant under some suitable choice.
\section{ Pseudo-differentiable measures.}
\par {\bf 4.1. Definition and notes.} 
A function $f: {\bf K}\to \bf R$ is called
pseudo-differentiable of order $b$, if there exists
the following integral:
$PD(b,f(x)):=\int_{\bf K}[(f(x)-f(y)) \times g(x,y,b)]
dv(y)$. We introduce the following notation $PD_c(b,f(x))$ for
such integral by $B({\bf K},0,1)$ instead of the entire $\bf K$.
Where $g(x,y,b):=\mid x-y\mid ^{-1-b}$ 
with the corresponding Haar measure $v$
with values in $\bf R$, where $b \in {\bf C}$ and 
$|x|_{\bf K}=p^{-ord_p(x)}$.
\par  Obviously, the definitions of differentiability of measures
can not be transferred from \cite{smbo,dal} onto the case considered here.
This is the reason why the notion
of pseudo-differentiability is introduced here. 
A quasi-invariant measure $\mu $ on $X$
is called pseudo-differentiable for $b \in \bf C$,
if there exists $PD(b,g(x))$ for $g(x):=\mu (-xz+S)$ for each
$S \in Af(X,\mu )$ with $|m|(S)<\infty $
and each $z \in J^b_{\mu }$, where $J^b_{\mu }$ is a $\bf K$-linear
subspace dense in $X$. For a fixed $z \in X$ such measure is called
pseudo-differentiable along $z$.
\par  For a one-parameter subfamily of operators $B({\bf K},0,1)\ni t \to
U_t:X \to X$ quasi-invariant measure $\mu $ is called 
pseudo-differentiable for $b \in \bf C$,
if for each $S$ the same as above there exists $PD_c(b,g(t))$ for a function
$g(t):=\mu (U_t^{-1}(S)$, where $X$ may be also a topological group
$G$ with a measure quasi-invariant relative to a dense subgroup
$G'$ (see \cite{lu3,luum983,lubp2}).
\par {\bf 4.2.} Let $\mu $, $X$, and $\rho $ be the same as in Theorem 3.24
and $F$ be a non-Archimedean Fourier transform defined in \cite{vla,roo}.
\par {\bf Theorems.} {\it $(1) $ $g(t):=\rho (z+tw,x)j(t) \in L^1(v,
K\to {\bf C}):=V$,  moreover, there exists
$F(g) \in C^0({\bf K},{\bf C})$ and $\lim_{
\mid t\mid \to \infty } F(g)(t)=0$  for $\mu $ and $v$ with values in
$\bf R$, where
$z$ and $w \in J_{\mu }$, $t \in \bf K$, 
$j(t)$ is the characteristic function
of a compact subset $W \subset \bf K$. In general, may be
$k(t):=\rho (z+tw,x) \notin V$.
\par  $(2)$ Let $g(t)=\rho (z+tw,x)j(t)$
with clopen subsets $W$ in $\bf K$. Then there are $\mu $,
for which there exists
$PD(b,g(t))$ for each $b \in {\bf C}$. If
$g(t)=\rho (z+tw,x)$, then there are probability measures $\mu $,
for which there exists $PD(b,g(t))$ for each $b$
with $0<Re (b)$ or $b=0$.
\par $(3)$ Let $S \in Af(X,y)$, $|\mu |(S)<\infty $,
then for each $b \in U:= \{ b': Re \mbox{ }b'>0 \mbox{ or }b'=0 \}$
there is a pseudo-differentiable quasi-invariant measure $\mu $ .}
\par {\bf Proof.} We consider the following additive compact subgroup
$G_T:=\{ x \in X| \| x(j) \| \le p^{k(j)} \mbox{ for each }
j\in {\bf N} \}$ in $X$, where $T=diag \{ d(j)\in K:  \mbox{ }
|d(j)|=p^{-k(j)} \mbox{ for each } j \in {\bf N} \}$ is a compact
diagonal operator.
Then $\mu $ from Theorem 3.20 is quasi-invariant 
relative to the following additive subgroup
$S_T:=G_T+H$, where $H:=sp_{\bf K}\{ e(j): j \in {\bf N} \}$.
Moreover, for each $z \mbox{ and }w\in G_T$ and $R>0$ there is $c>0$
such that $\rho (z+uw,x)=\rho (z+sw,x)$ for $|u-s|<c$ and
$x \in B(X,0,R)$, if all functions $f_j$ in the proof of Theorem
3.20 are locally constant (that is, $f$ are defined on 
${\bf K}e_j\subset X$).
In general, for each $b \in \bf C$ we can choose a sequence
$h_j(x)$ with $\sum _{j=0}^{\infty } \sup_{x \in \bf K} (|h_j(x)/
(h_j(x)+f_j(x))|r_j(x))<c'$ for suitable fixed $c'>0$ and functions
$r_j: {\bf K} \to [0, \infty )$ with
$\lim_{|x| \to 0}r_j(x)=\infty $. Carrying out calculations and using
the fact that $\mu $ is the probability quasi-invariant measure we get 
the pseudo-differentiability
of $\mu $. Using the Riemann-Lebesgue theorem from \cite{vla}
we get the statement of this theorem for $F(g)$.
\par  {\bf 4.3. } Let $X$ be a Banach space over
$\bf K$, $b_0\in {\bf R}\cup \{ + \infty \} $
and suppose that the following conditions are satisfied:
\par $(1)$ $T: X\to X$ is a compact operator with $ker(T)=\{ 0 \}$;
\par $(2)$ a mapping $\tilde F$ from $B({\bf K},0,1)$ to $C_T(X):=\{ U: U 
\in C^1(X,X)$
and $(U'(x)-I)$ is a compact operator for each $x \in X,$
there is $U^{-1}$ satisfying the same conditions as $U \} $ is given;
\par $(3)$ $\tilde F(t)=U_t(x)$ and  $\bar \Phi ^1 U_t(x+h,x)$
are continuous by
$t$, that is, $\tilde F \in C^1(B({\bf K},0,1),C_T(X))$;
\par $(4)$ there is $c>0$ such that
$\| U_t(x)-U_s(x) \| \le \|Tx \|$ for each $x \in X$ and $|t-s| <c$;
\par $(5)$ for each $R>0$ there is a finite-dimensional over $\bf K$ subspace
$H \subset X$ and $c'>0$ such that $\|U_t(x)-U_s(x)\| \le \|Tx\|/R$
for each
$x \in X \ominus H$ and $|t-s| <c'$ with $(3-5)$ satisfying also for
$U^{-1}_t$. 
\par {\bf Theorem.} {\it 
On $X$ there are probability quasi-invariant measures $\mu $ which are
pseudo-differentiable for each $b\in \bf C$ with ${\bf R}\ni
Re (b)\le b_0$ relative to a family $U_t$.}
\par {\bf Proof.} From Conditions (2,3) it follows that there is
$c>0$ such that $|det(U_t'(x))|=|det(U_s'(x))|$ for $\mu $-a.e.
$x \in X$ and all $|t-s|<c$, where quasi-invariant and pseudo-differentiable
measures $\mu $ on $X$ relative to $S_T$
may be constructed as in the proof of Theorems 3.20 and 4.2.
From $(1-5)$ it follows that conditions of Theorem
3.25 are satisfied for each
$U_t$. From $(3,5)$ it follows that for each $R>0$ and $b>0$
there exists $c'>0$ such that $|\rho _{\mu }(x-U^{-1}_t(x),x
)-\rho _{\mu }(x-U^{-1}_t(x),x)| \le
b$ for each $|t-s|<c'$ and $x \in B(X,0,R)$. If during construction of
$\mu $ to use only locally constant $f_j(x_j)$ with $h_j=0$, then
we can take $\infty \ge b_0\ge 0$.
Indeed, let $T=SCDE$ be a decomposition from
the appendix and $\mu ^A(\tilde S)=\mu (E\tilde S)$ with
$A=E^{-1}$ (or $C$), $\tilde S\in Bf(X)$,
then $\rho _{\mu ^A}(x-U_t^{-1}(x),x)$
$=\rho _{\mu }(Ex-U_t^{-1}(Ex),Ex)$
and $c(x):=\|U_t^{-1}(Ex)-U_s^{-1}(Ex)\| \le \|TEx\|$ for $x \in X$
and $c(x) \le \|TEx\|/R$ for $x \in X \ominus H$ (for each $R>0$
there is $H \subset X$ with $dim_{\bf K}H< \infty $), due to
compactness of $E-I$ and $C-I$ we can choose $H$ such that $A^{-1}H=H$,
$AH=H$.  We use also images of quasi-invariant 
measure constructed with the help of compact diagonal
operator $(D'-I)$, where $|D_j-1|< |D'_j-1|p^{k(j)}$,
$\lim_j k(j)=- \infty $.
Then $\sup\{ \|U_t(x)-U_t(0)\|:
x \in B(X,0,R) \}$ $\le \| \bar \Phi ^1 U_t \| _{B(X,0,R)} \times
\| x \| $. From the fact that $H$ is finite-dimensional and the existence
of the orthogonal projector (that is, corresponding to the decomposition into
the direct sum) $\pi _H: X \to H$ it follows that $\| U_t^{a}(e_j)-
U_s^{a}(e_j)\| \le b_j$ for each $j \in \bf N$ and each
$|t-s| <c'$, where $\lim _{j \to \infty } b_j =0$, $a=1 $ or
$a=-1$. This guarantees pseudo-differentiability of $\mu $.
\par {\bf 4.4.} Let $X$ be a Banach space over $\bf K$, $\mu $ be
a probability quasi-invariant measure  $\mu : Bf(X) \to \bf R$,
that is pseudo-differentiable for a given $b$ with $Re(b)>0$,
$C_b(X)$ be a space of continuous bounded functions
$f: X \to \bf R$ with $\| f\| :=\sup_{x \in X} |f(x)|$. 
\par {\bf Theorem.} {\it For each $a \in J_{\mu }$  and
$f \in C_b(X)$ is defined the following integral:
$$(i) \mbox{ }l(f)=\int _K[\int_X f(x)[\mu (-\lambda a +dx)-
\mu (dx)] g(\lambda ,0,b) v(d \lambda )$$ and there exists
a measure $\nu : Bf(X) \to \bf C$ 
with a bounded variation (for $b \in \bf R$ this $\nu $ is a mapping
from $Bf(X)$ into $\bf R$) such that
$$(ii) \mbox{ } l(f)=\int_X f(x) \nu (dx),$$  where $v$ 
is the Haar measure on $\bf K$ with values in 
$[0, \infty )$,
moreover, $\nu $ is independent from $f$ and may be dependent on
$a \in J_{\mu }$. We denote $\nu =:\tilde D^b_a \mu $.}
\par {\bf Proof.} From Definition 4.1 and the Fubini and Lebesgue
theorems it follows that there exists
$$\lim_{j \to \infty } \int_{{\bf K}\setminus
B({\bf K},0,p^{-j})}[\int_X (f(x+ \lambda a)-f(x)) g(\lambda ,0,b)
\mu (dx)] v(d \lambda )=l(f),$$ that is $(i)$ exists. 
Let 
$$(iii)\quad l_j(V,f):=\int_{{\bf K}\setminus B({\bf K},0,p^{-j})} 
[\int_V f(x)( \mu (-\lambda a+dx)- \mu (dx)) g(\lambda ,0,b)] v(d \lambda ),$$
where $V \in Bf (X)$.
The proofs of \S \S 2 and 3 in Theorem 1 \S 21\cite{sko}
may be slightly transformed for the considered here case of the 
Banach space $X$ over $\bf K$ instead of a Hilbert space
over $\bf R$. Hence for each $c>0$
there exists a compact $V_c \subset X$ with $| \nu _{\lambda }|
(X\setminus V_c) <c$ for each $|\lambda |>0$, where $\nu _{\lambda }(A):=$
$\int_{{\bf K}\setminus B({\bf K},0,
|\lambda |)} [\mu (-\lambda 'a +A) - \mu (A)] g(\lambda ',0,b)
v(d \lambda ')$ for $A \in Bf(X)$. Also there are
$\delta >0$ and $V_c$ such that $- \lambda 'a+V_c \subset V_c$  and
$\| [X\setminus V_c)\bigtriangleup (- \lambda 'a+
(X\setminus V_c))] \| _{\mu }=0$
for each $|\lambda '|< \delta $ (see also Theorem 7.22 \cite{roo}),
where $A\bigtriangleup B:=(A\setminus B) \cup (B\setminus A)$.
\par  Each continuous linear functional on
$C_b(V)$ for compact $V$ has the form $(ii)$ (see \cite{boui},
A.5 \cite{sch1}),
that is, for $l_j(V,f)$ there exists a measure $\nu _{\lambda }(c,*)$.
In view of the Alaoglu-Bourbaki theorem each bounded subset
$W$ in $[C_b(V_c)]'$
is precompact (i.e., $cl(W)$ is compact) in the weak topology,
where $C_b(V_c)$ is separable. For each $c>0$ we have 
$$l(f)=\int_{V_c} f(x) \nu (c, dx)+O(c\times \sup_{x \notin V_c} |f(x)|)$$
such that $\nu (c,*)$ are measures on $Bf(V_c)$. Choosing
$V_c \subset V_{c'}$
for $c>c'$ and defining $\nu (A):=\lim_{c \to +0} \nu (c, A\cap V_c)$
for each $A \in Bf(X)$ in view of the Radon-Nikodym theorem about
convergence of measures (see \cite{cons}) we get $(i)$ with the help
of the theorem about extension of measures, since the minimal
$\sigma $-field $\sigma Bco(X)$ generated by $Bco(X)$ coincides with
$Bf(X)$.
\par  {\bf 4.5. Theorem.} {\it If $\mu $ is probability real-valued
quasi-invariant measure on a Banach space
$X$ over ${\bf K}=\bf Q_p$, $a \in J_{\mu }$ and $\mu $
is pseudo-differentiable of order $1$, then for each
$f \in L^{\infty }(X,\mu ,{\bf R})$ with $supp( f) \subset
B(X,0,R)$, $\infty >R=R(f)>0$
the following equality is accomplished:
$$(i) \mbox{ }\int_X[D_{\lambda }^{-1} f(x+ \lambda a)]|^{\lambda
=\beta }_{\lambda =0} (\tilde D_a^1 \mu )(dx)\mbox{ }=\int_X[f(x+\beta
a) - f(x)] \mu (dx),$$ moreover, pseudo-differential
$\tilde D_a^1 \mu :=\nu $
in the direction $a$ from Theorem 4.4 given by formula
$(ii)$ characterizes $\mu $ uniquely.}
\par {\bf Proof.} From the restrictions on $f$ it follows that
$f(x+\lambda
a )=:\psi (\lambda )\in \sf E'$, where $\sf E'$ is the topologically
adjoint space to the space $\sf E$ of basic locally constant functions
\cite{vla}.
By formula $(1.5)$ \S 9\cite{vla} for them there is the following equality
$D^aD^b\psi =D^bD^a \psi =D^{a+b}
\psi $ for $a+b \ne -1$, $a, b \in \bf R$ ( $D^a_{\lambda }$ denotes
a pseudo-differentiation of functions by $\lambda $). In view of Theorem in
\S 6.3\cite{vla} $\sf E'$ $=[\phi \in {\sf D'}:$ $ supp (\phi )$ is compact
$]$, where $\sf D$ $=[\phi \in {\sf E}:$ $supp(\phi )$ is compact $]$.
Since $[\xi ]|_0^{\beta }:=\xi (\beta ) - \xi (0)$ for $\xi (\lambda )=
D_{\lambda }^{-1}\psi (\lambda )$ and $(D^{-1}D \psi )(\lambda )=
\psi (\lambda )$, hence due to the definition of $D_{\lambda }^{-1}$
we have:
$$p^2(p+1)^{-1} \int_X [D_{\lambda }^{-1} f(x+\lambda a)]|_0^{\beta }
(\tilde D_a^1 \mu )(dx)=\lim_{\alpha \to -1} \int_{\bf Q_p}
\int_{\bf Q_p} (1-p^{\alpha })$$ 
$$(1-p^{-1-\alpha })^{-1}p^2(1+p)^{-1}
[\int_X(\mu (-a(\lambda ' +\lambda ")+dx)- \mu (-a\lambda ' +dx))$$
$$|\lambda "|_p^{-2} |\lambda - \lambda '|_p^{-1-\alpha } \mu (dx)]
v(d \lambda ') v(d \lambda ").$$
In view of the Fubini and Lebesgue theorems
we establish Formula (i). Taking all different $f \in L^{\infty }
(X,\mu ,{\bf R})$ such that $\psi (\lambda )=f(x+a \lambda )$ are locally
polynomial by $\lambda $ and using the non-Archimedean variant of the
Stone-Weiestrass theorem by Kaplanski \cite{sch1}
we get that $\mu $ is characterized by
the left side of $(i)$ uniquely.
\par {\bf 4.6.} (Properties of $\tilde D_a^1 \mu $ for real-valued $\mu $,
see \S \S 4.4. and 4.5, where ${\bf K}=\bf Q_p$.)
\par {\bf Proposition.1.} 
{\it If $\mu $ is pseudo-differentiable in a direction
$a\in J^1_{\mu }$, $a \ne 0$, then for each $d \in J^1_{\mu }$
a measure $\mu ^d(A):=
\mu (A-d)$, $A \in Bf(X)$, is also pseudo-differentianble in a direction
$a$, moreover, $\tilde D_a^1 (\mu ^d)=(\tilde D_a^1 \mu )_d$.
If $a \in J^b_{\mu }$, then $\tilde D^b_{\lambda a}\mu $
exists for each $\lambda \in \bf K$.}
\par {\bf Proof.} From \S \S 4.4 and 4.5 it follows that 
$$\int_Xf(x)(\tilde
D_a^{-1}\mu )^d (dx)= \int_{\bf K}\int_X [f(x+\lambda a)-f(x)]
|\lambda |_{\bf K}^{-2} \mu ^d(dx) v(d \lambda )$$ $$=\int_X f(x)(\tilde D_a^1
(\mu ^d))(dx).$$ The last statement follows from Definition 4.1
with the help of \S 3.21.
\par  {\bf Proposition.2.} { \it If $a \in J^1_{\mu }$, $ \nu =\tilde D_a^1 \mu $ and
$D_{\lambda }^{-1}\rho _{\nu }(\lambda a,x)$ $=:\phi (x) \in
L^1(X, \nu , {\bf R})$ by $x$, then $\mu ^{\lambda a} \ll \nu $.}
\par {\bf Proof.} In view of Formula $4.5.(i)$ the following
equality is accomplished: 
$$\int_Xf(x)\mu ^{\lambda a} (dx)=\int_X
f(x)\mu (dx)+\int_X f(x)[D_{\xi }^{-1} \rho _{\nu }(\lambda a,x)]
|_{\xi =0}^{\xi =\lambda } \nu (dx).$$
 Taking $f(x)=\chi _A(x)$,
where $\chi _A$ is a characteristic function of $A \in Bf(X)$
we get $\mu ^{\lambda a}(A)=\mu (A)$ for every $|\lambda |>0$.
Therefore, considering $|\lambda | \to \infty $ we get
$\mu (A)=0$, that is, $\mu \ll \nu $ and $\mu ^{\lambda a} \ll \nu $.
\par  {\bf Proposition.3.} If a non-negative measure $\mu $
is pseudo-differentiable in a direction
$a \in J^1_
{\mu }$ and $(\tilde D_a^1 \mu )\ll \mu $, then $l_{\mu }(a,x):=
(\tilde D_a^1\mu )(dx)/ \mu (dx)$ is called the logarithmic
pseudo-derivative of quasi-invariant measure $\mu $ in a direction $a$.
Then due to Formula $4.5.(i)$
there is the equality: 
$$\int_X[f(x+\lambda a)-f(x)] \mu (dx)=
\int_X[\tilde D_{\xi }^{-1}f(x+\xi a)]|_{\xi =0}^{\xi =\lambda }
l_{\mu }(a,x) \mu (dx).$$ We may slightly transform for the
non-Archimedean case also Theorem 1 \S 22\cite{sko} with the help
of the results of \S \S 4.2-4.4 above.
\par {\bf 4.7. Theorem.}
{\it Let $X$ be a Banach space over $\bf K$, $|*|=mod_{\bf K}(*)$
with a probability quasi-invariant measure $\mu : Bf(X) \to \bf R$ 
and it is satisfied Condition $3.21(i)$, 
suppose $\mu $ is pseudo-differentiable and
\par $(viii)$ $J_b{\mu } \subset T"J_{\mu }$, $(U_t:$ $ t \in B({\bf K},0,1))$
is a one-parameter family of operators such that
Conditions $3.25(i-vii)$ are satisfied
with the substitution of $J_{\mu }$ onto $J^b_{\mu }$
uniformly by $t \in B({\bf K},0,1)$, $J_{\mu } \supset T'X$, where
$T', T": X \to X$ are compact operators, $ker (T')=ker (T")=0$.
Moreover, suppose that there are sequences 
\par $(ix)$ $[k(i,j)]$ and $[k'(i,j)]$
with $i, j \in \bf N$, $\lim_{i+j \to \infty }k(i,j)=$ $\lim_{i+j\to
\infty }k'(i,j)=- \infty $ and $n \in \bf N$ such that
$|T"_{i,j}-\delta _{i,j}|<|T'_{i,j}-\delta _{i,j}|p^{k(i,j)}$,
$|U_{i,j}-\delta _{i,j}|< |T"_{i,j}-\delta _{i,j}|p^{k'(i,j)}$
and $|(U^{-1})_{i,j}-\delta _{i,j}|<|T"_{i,j}-\delta _{i,j}|p^{k'(i,j)}$
for each $i+j>n$, where $U_{i,j}=\tilde e_iU(e_j)$, $(e_j:$ $j)$ is
orthonormal basis in $X$. Then for each $f \in C_b(X)$ there is defined
$$(i)\mbox { }l(f)=\int_{B(K,0,1)}[\int_X f(x)[\mu (U^{-1}_t(dx))-
\mu (dx)]g(t,0,b)v(dt)$$ and there exists a measure $\nu : Bf(X) \to
\bf C$ with a bounded total variation [particularly, for
$b \in \bf R$ it is such that $\nu : Bf(X) \to \bf R$] and
$$(ii)\mbox{ }l(f)=\int_Xf(x)\nu (dx),$$ where $\nu $ is independent
from $f$ and may be dependent on $(U_t:$ $t)$, $\nu =:\tilde D^b_{U_*}\mu $.}
\par {\bf Proof.} From the proof of Theorem 3.25
it follows that there exists a sequence
$U^{(q)}_t$ of polygonal operators converging uniformly by
$t \in B({\bf K},0,1)$ to $U_t$ and equicontinuously
by indices of matrix elements
in $L^1$. Then there exists $\lim_{q\to \infty }
\lim_{j \to \infty }\int_{B({\bf K},0,1)
\setminus B({\bf K},0,p^{-j}}[\int_Xf(U^{-1}_t(x))-f(x)]g(t,0,b)
\mu (dx)]v(dt)$ for each $f \in C_b(X)$. From conditions $(viii, ix),$
the Fubini and Lebesgue theorems it follows that for $\nu _{\lambda }:=
\int_{B({\bf K},0,1)\setminus B({\bf K},0,|\lambda |)}[\mu (U^{-1}_t(A))-
\mu (A)]g(t,0,b)v(dt)$ for $A \in Bf(X)$ 
for each $c>0$ there exists a compact $V_c \subset X$ and $\delta >0$
such that $| \nu _{\lambda }|(X\setminus V_c)<c$.
Indeed, $V_c$ and $\delta >0$ may be chosen due to
pseudo-differentiability of $\mu $, \S \S 2.35, 3.9, 3.23,
Formula $(i)$, $3.21(i)$ and due to continuity and boundednessy
(on $B({\bf K},0,1) \ni t$) of $|det$ $U'_t(U^{-1}_t)(x))|_{\bf K}$
satisfying the following conditions $U^{-1}_t(V_c) \subset V_c$
and $\| (X\setminus V_c)
\bigtriangleup (U^{-1}_t(X\setminus V_c))\| _{\mu }=0$ for each
$|t|< \delta $, since $V_c=Y(j)\cap V_c$ are compact for every
$j$. At the same time closed subsets $Y(j) \subset X$ for $U^{(q)}_t$
may be chosen independent from $t$ and in \S 2.35 operator $S_c$ may
be chosen symmetric, $TX
\supset J_{\mu }\supset T'X$, where $T$ is a compact operator.
Evidently, conditions of type $(ix)$ are carried out for $V(j,x)$
and $U(j,x)$
uniformly by $j$. Repeating proofs 3.25 and 4.4 with the use of
Lemma 2.5 for the family $(U_t:$ $t)$ we get formulas $(i,ii)$.
\section{ Convergence of quasi-invariant and pseudo-differentiable
measures.}
\par  Different types of convergence of measures were considered
in \cite{smbo,cons,smfo,top1,top2}. The definitions
and theorems given below are modifications
of the corresponding chapters of
\cite{cons}, but additionally taking into account the properties
of quasi-invariance and pseudo-differentiability.
\par {\bf 5.1. Definitions, notes and notations.} Let $S$ be a normal
topological group with the small inductive dimension
$ind (S)=0$, $S'$ be a dense subgroup, suppose their topologies are
$\tau $ and $\tau '$ correspondingly, $\tau ' \supset \tau |_{S'}$.
Let $G$ be an additive Hausdorff left-$R$-module, where $R$ is a
topological ring, ${\sf R } \supset Bf(S)$ be a  $\sigma $-ring
for real-valued measures,
${\sf M(R,} G)$ be a family of measures with values in $G$,
${\sf L(R,} G,R)$ be a family of quasi-invariant measure 
$\mu : {\sf R} \to G$ with $\rho _{\mu }
(g,x)\times \mu (dx):=\mu ^{g^{-1}}(dx)=:\mu (gdx)$,
$R \times G \to G$ be a continuous left action of $R$ on $G$
such that $\rho _{\mu }(gh,x)=\rho _{\mu}(g,hx)
\rho _{\mu }(h,x)$
for each $g,h \in S'$ and $x \in S$. Particularly,
$1=\rho _{\mu }(g,g^{-1}x)
\rho _{\mu }(g^{-1},x)$, that is, $\rho _{\mu }(g,x) \in R_o$,
where $R_o$ is a multiplicative subgroup of $R$.
Moreover, $zy \in \sf L$ for $z \in R_0$ with $\rho _{z\mu }(g,x)=
z\rho _{\mu }(g,x)z^{-1}$ and $z \ne 0$. We suppose that topological
characters and weights $S$ and $S'$ are countable and each open
$W$ in $S'$ is precompact in $S$.  Let $\bf P"$ be a family of
pseudo-metrics in $G$ generating the initial uniformity such that for each
$c>0$ and $d \in \bf P"$ and $\{U_n\in {\sf R}: n \in {\bf N} \}$ with
$\cap \{ U_n: n \in {\bf N} \}= \{x\}$
there is $m \in \bf N$ such that $d(\mu ^g(U_n), \rho _{\mu }(
g,x)\mu (U_n))<cd(\mu (U_n),0)$ for each $n>m$, in addition,
a limit $\rho $ is independent $\mu $-a.e. on the choice of
$\{U_n: n\}$ for each $x \in S$ and $g \in S'$.
Consider a subring $R' \subset \sf R$,  $R' \supset Bf(S)$
such that $\cup \{A_n: n=1,...,N\} \in R'$
for $A_n \in R'$ with $N \in \bf N$ and $S'R' =R'$. Then
${\sf L(R},G,R;R')
:= \{ (\mu ,\rho _{\mu }(*,*)) \in {\sf L(R},G,R): \mu - \mbox{ }R'-
\mbox{ is regular and for each } s \in S \mbox{ there are } A_n \in R',
n \in {\bf N} \mbox{ with } s=\cap (A_n:n), \{s\} \in R'\}$.
\par  For pseudo-differentiable measures $\mu $ let $S" \subset S'$,
$S"$ be a dense subgroup in $S$, $\tau '|S"$ is not stronger than $\tau "$
on $S"$ and there exists a neighbourhood $\tau " \ni W" \ni e$ in which
are dense elements lying on one-parameter subgroups $(U_t:$
$t\in B({\bf K},0,1))$. 
We suppose that $\mu $ is induced from the Banach space
$X$ over  $\bf K$ due to a local homeomorphism of neighbourhoods
of $e$ in $S$ and $0$ in $X$ as for the case of groups of
diffeomorphisms \cite{lu1} such that is accomplished Theorem
4.7 for each $U_* \subset S"$ inducing the correspopnding transformations
on $X$. In the following case $S=X$ we consider $S'=J_{\mu }$ and
$S"=J^b_{\mu }$ with $Re(b)>0$ such that $M_{\mu }
\supset J_{\mu } \subset (T_{\mu }X)^{\sim }$, $J^b_{\mu }
\subset (T^{(b)}_{\mu }X)^{\sim }$ with compact operators
$T_{\mu }$ and $T_{\mu }^{(b)}$, $ker (T_{\mu })=ker (T^{(b)}_{\mu })=0$
and norms induced by the Minkowski functional $P_E$ for $E=T_{\mu }
B(X,0,1)$ and $E=T^{(b)}_{\mu }B(X,0,1)$ respectively. 
We suppose furter that for pseudo-differentiable measures
$G$ is equal to $\bf C$ or $\bf R$. We denote
$P({\sf R},G,R,U_*;R'):=[(\mu ,\rho _{\mu }, \eta _{\mu }):$ $(\mu ,
\rho _{\mu }) \in {\sf L(R},G,R;R'), \mu $ is pseudo-differentiable
and $\eta _{\mu }(t,U_*,A) \in L^1({\bf K},v,{\bf C})$],
where 
$$\eta _{\mu }(t,U_*,A)=j(t)g(t,0,b)[\mu ^h(U^{-1}_t(A)-\mu ^h(A)],$$
$j(t)=1$ for each $t \in \bf K$ for $S=X$; $j(t)=1$ for $t \in B({\bf K},0,1)$,
$j(t)=0$ for $|t|_{\bf K}>1$ for a topological group $S$ 
that is not a Banach space $X$ over $\bf K$,
$v$ is the Haar measure on $\bf K$ with values in $[0, \infty )$,
$(U_t:$ $t \in B({\bf K},0,1))$ is an arbitrary one-parameter subgroup.
On these spaces $\sf L$ (or $P$) the additional conditions are imposed:
\par (a) for each neighbourhood (implying that it is open) $U \ni 0 \in G$
there exists a neighbourhood $S\supset V\ni e$ and a compact subset $V_U$,
$e \in V_U \subset V$, with $\mu (B) \in U$ (or in addition
$\tilde D^b_{U_*}
\mu (B) \in U$) for each $B$, ${\sf R} \ni B
\in Bf(S \setminus V_U)$;
\par  (b) for a given $U$ and a neighbourhood
$R\supset D\ni 0$ there exists a neighbourhood $W$, $S'\supset
W \ni e$, (pseudo)metric $d \in P"$ and $c>0$ such that $\rho _{\mu }
(g,x)-\rho _{\mu }(h,x')
\in D$ (or $\tilde D^b_{U_*}(\mu ^g -\mu ^h)(A) \in U$ for $A \in
Bf(V_U)$ in addition for $P$)
whilst $g, h \in W$, $x,x' \in V_U$, $d(x,x') < c$, where (a,b) is
satisfied for all $(\mu ,\rho _{\mu }) \in \sf L$ (or $(\mu ,\rho _{\mu },
\eta _{\mu }) \in P$) equicontinuously in (a) on $V\ni U_t, U^{-1}_t$ and
in (b) on $W$ and on each
$V_U$ for $\rho _{\mu }(g,x)-\rho _{\mu }(h,x')$
and $\tilde D^b_{U_*}(\mu ^g-\mu ^h)(A)$.
\par These conditions are justified, since for a Gaussian
measure $\nu $ on a Hilbert space $Z$ over $\bf C$ or $\bf R$
(or measures given above for a non-Archimedean Banach space 
$Z$ over $\bf K$)
there exists a subspace $Z'$ dense in $Z$ such that
$\nu $ is quasi-invariant relative to
$Z'$. Moreover, in view of Theorem 26.2\cite{sko} (or Theorems 3.20,
3.25, 4.3 and 4.7 in the non-Archimedean case respectively)
there exists a subspace $Z"$ dense in $Z'$ such that
for each $\epsilon >0$ and each $\infty >R>0$ there are $r>0$ and $\delta >0$
with $|\rho _{\nu }(g,x)- \rho _{\nu }(h,y) | <\epsilon $ for each
$\| g- h \|_{Z"}+ \| x-y\|_Z<\delta $, $g, h \in B(Z",0,r)$, $x, y\in
B(Z,0,R)$, where $Z"$ is the Banach space over $\bf C$ or $\bf R$ or
$\bf K$ respectively.
For a group of diffeomorphisms of a Hilbert manifold over $\bf R$
or over a non-Archimedean Banach manifold we have an analogous continuity
of $\rho _{\mu }$ for a subgroup $ G" $ of the entire group $G$
(see \cite{lu1,lu3,lubp2,lu5}).
\par By $\sf M_o$ we denote a subspace in $\sf M$, satisfying (a).
Henceforth, we imply that $R'$ contains all closed subsets
from $S$ belonging to $\sf R$, where $G$ and $\sf R$ are complete.
\par   For $\mu : Bf(S) \to G$ by $L(S,\mu ,G)$  we denote the completion
of a space of continuous $f: S\to G$ such that $\| f\|_d:=\sup_{
h \in C_b(S,G)}d(\int_Sf(x)h(x)$ $\mu (dx), 0)< \infty $ for each
$d \in P"$, where $C_b(S,G)$ is a space of continuous bounded functions
$h: S\to G$. We suppose that for each sequence $(f_n:$ $n) \subset
L(S,\mu ,G)$ for which  $g \in L(S,\mu ,G)$ exists with
$d(f_n(x),0) \le d(g(x),0)$ for every $d \in P"$, $x$ and $n$,
that $f_n$ converges uniformly on each compact subset $V \subset S$
with $|\mu |(V)>0$. In the cases $G=\bf C$
it coincides with $L^1(S,\mu ,{\bf C})$
correspondingly, hence this supposition is the Lebesgue theorem.
By $Y(v)$ we denote $L^1({\bf K},v,{\bf C})$.
\par   Now we may define topologies and uniformities with the help
of corresponding bases (see below) on ${\sf L}\subset G^{\sf R}
\times R_o^{S'\times S}=:Y$ (or
$P \subset G^{\sf R}\times R_o^{S'\times S}\times
G^{S'\times K\times \sf R}=:Y$), $R_o\subset R\setminus \{0\}$. There are
the natural projections $\pi : \sf L$ $(or $P$)
\to M_o$, $\pi (\mu ,\rho _{\mu }(*,*)$ $(\vee $ ,$\eta _{\mu }))=\mu $,
$\xi: {\sf L}$ $(\vee $ $P) \to R^{S'\times S}$,
$\xi (\mu ,\rho _{\mu },(\vee$ $\eta _{\mu }))=\rho _{\mu }$,
$\zeta : P\to G^{S'\times K\times \sf R}$, $\zeta (\mu ,\rho _{\mu },
\eta _{\mu })=\eta _{\mu }$. Let $\sf H$ be a filter on $\sf L$ or $P$,
$U=U'\times U"$ or $U=U'\times U"\times U"'$, $U'$ and $U"$ be elements
of uniformities on $G$, $R$ and $Y(v)$ correspondingly,
$\tau ' \ni W\ni e$, $\tau \ni V \supset V_{U'} \ni e$, $V_{U'}$ is compact.
By $[\mu ]$ we denote $(\mu ,\rho _{\mu })$ for $\sf L$
or $(\mu ,\rho _{\mu },\eta _{\mu })$ for $P$, $\Omega :=\sf L$ $\vee $
$P$, $[\mu ](A,W,V):=[\mu ^g(A), \rho _{\mu }(g,x),$ $\vee $
$\eta _{\mu ^g}(t,U_*,A)|$ $g \in W, x \in V,$ $\vee $ $t \in K]$.
We consider $\sf A \subset R$, then
$$(1)\mbox{ }{\sf W(A,}W,V_{U'};U):=\{ ([\mu], [\nu ]) \in \Omega ^2|
([\mu ],[\nu ])(A,W,V_{U'}) \subset U \};$$
$$ (2)\mbox{ }{\sf W(S};U):=\{ ([\mu ],[\nu ])\in \Omega ^2|
\{ (B,g,x) : ([\mu ],[\nu ])(B,g,x))\in  U \} \in
{\sf S} \},$$ where  ${\sf S}$ is a filter on ${\sf R} \times
S' \times S^c$, $S^c$ is a family of compact subsets $V' \ni e$.
$$ (3) \mbox{ }{\sf W(F},W,V;U):=\{ ([\mu ],[\nu ]) \in \Omega ^2|
\{ B: ([\mu ],[\nu ])(B,g,x) \in U, g \in W, x \in V \}
\in {\sf F} \},$$ where ${\sf F}$ is a filter on $\sf R$
(compare with $\S $ 2.1 and 4.1\cite{cons});
$$(4) \mbox{ }{\sf W(A,G};U) := \{ ([\mu ],[\nu ]) \in \Omega ^2|
\mbox{ }\{ (g,x): ([\mu ],[\nu])(B,g,x) \in U, B \in {\sf A} \}
\in {\sf G} \},$$ where $\sf G$ is a filter on $S' \times S^c$;
suppose ${\sf U \subset R} \times \tau '_e \times S^c$, $\Phi $
is a family of filters on ${\sf R} \times S' \times S^c$ or
${\sf R}\times S'\times S^c\times Y(v)$ (generated by products of
filters $\Phi_{\sf R} \times \Phi _{S'} \times \Phi _{S^c}$
on the corresponding spaces), ${\sf U'}$ be a uniformity on
$(G,R)$ or $(G,R,Y(v))$,
${\sf F} \subset Y$. A family of finite intersections of sets
${\sf W}(A,U) \cap ({\sf F \times F})$  (see (1)), where
$(A,U) \in {\sf U \times U'}$  (or
$ {\sf W}(F,U) \cap ({\sf F \times F})$ (see (2)), where
$(F,U) \in (\Phi \times {\sf U'}$) generate by the definition
a base of uniformity of $\sf U$-convergence ($\Phi $-convergence
respectively)
on $\sf F$ and generate the corresponding topologies. For these
uniformities are used notations 
$$(i)\mbox{ }{\sf F_U} \mbox{ and }
{\sf F}_ {\Phi };\quad {\sf F}_{{\sf R}
\times W \times V}\mbox{ is for }{\sf F}\mbox{ with the uniformity of uniform
convergence}$$
 on ${\sf R} \times W \times V$, where $W \in \tau '_e$,
$V \in S^c$, analogously for the entire space $Y$;
$$(ii)\mbox{ }{\sf F}_A\mbox{ denotes the uniformity (or topology)
of pointwise convergence for}$$
$A \subset {\sf R} \times \tau '_e \times
S^c=:Z$, for $A=Z$ we omit the index (see formula (1)).
Henceforward, we use $\sf H'$ instead of $\mbox{ }\sl H$ in 4.1.24\cite{cons},
that is, ${\sf H'}(A,{\sf \tilde R})$-filter on $\sf R$ generated by the
base $[(L \in {\sf R}:$ $L \subset A \setminus K'):$
$K' \in \sf \tilde R,$
$K' \subset A]$, where ${\sf \tilde R} \subset \sf R$ and $\sf \tilde R$
is closed relative to the finite unions.
\par For example, let $S$ be a locally $\bf K$-convex space, $S'$ be a dense
subspace, $G$ be a locally $\bf L$-convex space, where
$\bf K, L$ are fields, $R=B(G)$ be a space of bounded linear operators on
$G$, $R_o=GL(G)$ be a multiplicative group of invertible linear
operators. Then others possibilities are: $S=X$ be a Banach space 
over $\bf K$, $S'=J_{\mu }$, $S"=J_{\mu }^b$
as above; $S=G(t)$, $S'\supset S"$ are dense subgroups,
$G=R$ be the field $\bf R$, $M$ be an
analytic Banach manifold over $\bf K \supset
Q_p$ (see \cite{lu1}). The rest of the necessary standard definitions
are recalled further when they are used.
\par {\bf 5.2. Lemma.} {\it Let $\sf R$ be a quasi-$\delta $-ring
with the weakest uniformity in which each $\mu  \in \sf M$ is uniformly
continuous and $\Phi \subset \hat \Phi _C({\sf R}, S' \times S^c)$. Then
${\sf L(R}, G,R,R')_{\Phi }$ (or $P({\sf R},G,R,U_*;R')_{\Phi }$) )
is a topological space on which $R_o$ acts continuously from the right.}
\par {\bf Proof.} We recall that $\hat \Phi :=\hat \Phi (X)$
denotes a family of filters $F$ on $X$ such that
for each mapping $f: \Phi \to {\sf B}(X)$ with $f(\Sigma )
\subset \Sigma$ for each $\Sigma \in \Phi $ there exists a finite
subset $\Psi \subset \Phi $ such that $\bigcup_{\Sigma \in \Psi }
f(\Sigma ) \in F$, where ${\sf B}(X)$ denotes a family of all subsets
in $X$, $\Phi _C(X)=:\Phi _C$ be a family of Cauchy filters on $X$.
In view of proposition 4.2.2\cite{cons} the space of measures
${\sf M(R},G;R')_{\Phi }$ is a topological left-$R$-module.
On the other hand, $\rho
_{\lambda \mu }(a,x)=\lambda \rho _{\mu }(a,x)\lambda ^{-1}$ for
$\lambda \in R_o$, from the continuity of the inversion
$\lambda \mapsto \lambda ^{-1}$
in $R_o$ and the multiplication in $R$ for each entourage of the diagonal
$U"$ in $R$
there exists an entourage of a diagonal $U"_1$ in $R$ and an open
neighbourhood $U_4 \ni \lambda $ in $R_o$ such that
$\beta U"_1 \beta ^{-1}\subset
U"$ for each $\beta \in U_4$. Choosing $\Psi \subset \Phi $ for a given
$F$ we find $U"_1(\Sigma )$ and $U_4(\Sigma )$ for each
$\Sigma \in \Psi $, then $\bigcap _{\Sigma \in \Psi }U"_1
(\Sigma )=:\tilde U"$ and $\tilde U_4=\bigcap_{\Sigma \in \Psi }
U_4(\Sigma )$ are entourages of the diagonal and generate a neighbourhood
of $\lambda $.
This shows the continuity by $\lambda $ for $L$ and $P$, since for
$P$ is satisfied: $PD(b, \lambda f(x))=
\lambda PD(b,f(x))$ due to definition 4.1 for pseudo-differentiable $f$.
\par {\bf 5.3. Proposition.} {\it (1). Let $\sf T$ be a
$\hat \Phi _4$-filter on ${\sf M_o(R,}G;R')$, $\{ A_n\}$ be a disjoint
$\Theta ({\sf R})$-sequence, $\Sigma $ be the elementary filter on $\sf R$
generated by $\{A_n:n \in {\bf N} \}$ and $\phi : {\sf M_o \times R}
\to G$ with $\phi (\mu ,A)=\mu (A)$. Then $\phi ({\sf T} \times \Sigma )$
converges to $0$. 
\par  (2). Moreover, let $\sf U$ be a base of neighbourhoods
of $e \in S'$, $\phi: {\sf L} \to G\times R$, $\phi (\mu ,A,g):=(\mu ^g(A),
\rho _{\mu }(g,x))$, where $x \in A$. Then $(0,1) \in
\lim \phi ({\sf T} \times \Sigma \times {\sf U})$.
\par  (3). If $T$ is a
$\Phi _4$-filter on $P({\sf R},G,R,U_*;R')$, \\
$\psi (\mu ,B,g,t,U_*)=
[\mu (B); \rho _{\mu }(g,x); \eta _{\mu ^g}(t,U_*,B)]$, \\
then
$(0,1,0) \in \lim \psi ({\sf T}\times \Sigma \times {\sf U})$ for each
given $U_* \in S"$, where $\Sigma $ and $\sf U$ as in (1,2).}
\par {\bf Proof.} We recall that for a set $X$, $\Theta $-net
in $X$ is a net $(I,f)$ in $X$ such that for each increasing sequence
$(i_n:$ $n \in {\bf N})\subset I$
is satisfied $(f(i_n): n) \in \Theta $. A sequence $(x_n: n
\in {\bf N})\subset X$ is called a $\Theta $-sequence, if
$f: {\bf N}\to X$, $f(n)=x_n$, $({\sf N}, f)$ is a $\Theta $-net.
Then $\Phi (\Theta ,X)=:\Phi (\Theta )$ denotes a family of filters
on $X$ of the form $f(F)$, where $(I,f)$ is a $\Theta $-net in $X$,
$F$ is a filter of sections $[(\lambda \in I:$ $\lambda \ge i):$
$i \in I]$ for a directed set $I$, $\hat \Phi (\Theta ):=
{\Phi (\Theta )}^{\hat { } }$. For a topological space $X$,
$\Theta _j(X)$ denotes a family of sequences having converging
subsequences for $j=1$ or having limit points whilst $j=2$. For a
uniform space $X$, $\Theta _3(X)$ denotes a family of sequences
having Cauchy subsequences $\Theta _4:=\Theta _2\cup \Theta _3$,
$\Phi _j(X):=\Phi (\Theta _j(X))$.
\par  Then (1) follows from 4.2.6\cite{cons}. (2) From the restrictions
(a,b) in 5.1 it follows that $\mu ^g(A)=\int _A \rho _{\mu }(g,x)
\mu (dx)$. If in addition $\lim \mu ({\sf F}(A;
R'))=\mu (A)$, then $\lim \mu ({\sf F}(X\setminus A,R'))=\mu (X\setminus A)$,
since by the definition 4.1.24\cite{cons} for each $\tau _G \ni D \ni 0$
for $A \in \sf R$ there exists $\tau _S \ni V \supset A$ with
$\mu (V)-\mu (A)\in D$, where ${\sf F}(A,R')$ is a filter generated by
$A$ and $R'$. Due to the Radonian property of $\mu $ and $\mu ^g$
we get $\lim_{{\sf T}\times \Sigma }
\mu ^g(A)=0$. 
\par  (3). Additionally to (1,2) it remains to verify that
$\eta _{\mu }$ converges. In view of Theorems 4.4 and 4.7 (or \S 5.1)
a pseudo-differential of order $b:$ 
$\tilde D^b_{U_*}\mu =\nu $
is a measure for pseudo-differentiable $\mu $. Due to (1) and
the $\hat \Phi _4$-condition we get $(0,1,0) \in \lim
\psi({\sf T}\times \Sigma\times {\sf U})$ for a given $U_*$,
since in \S \S 4.1 and 5.1 the integral is by $t \in \bf K$ for
the Banach space $X$ over $\bf K$ and by $B({\bf K},0,1)$ correspondingly for
$S$ that is not the Banach space.
\par {\bf 5.4. Proposition.} {\it Let $\sf H$ be a
$\hat \Phi _4 $-filter
on $\sf L$ (or $P$) with the topology $\sf F$ (see \S 5.1.(ii)),
$A \in \sf R$,
$\tau _G \ni U \ni 0$, ${\sf H'}(A,R') \in \Psi _f({\sf R})$. Then
there are $L \in \sf H$, $\tilde K \in R'$ and an element of the
uniformity $\sf U $ for ${\sf L}_{R'}$ or $P_{R'}$ such that
$\tilde {\bf K} \subset A$,
$L=[(\mu ,\rho _{\mu }(g,x)):$
$M:=\pi _{\sf M_o}(L) \ni \mu , \pi _{\tau '_e}(L)=:W \ni g $ (or
$(\mu ,\rho _{\mu },\eta _{\mu }(*,*,U_*))$ and additionally
$\tilde D^b_{U_*}\mu
=PD(b,\eta _{\mu }))]$, $e \in W \in \tau '$,
$\mu ^g(B)-\nu ^{h}(C) \in U$ (or in addition $(\tilde D^b_{U_*}
\mu ^g(B)) -(\tilde D^b_{U_*}\nu ^h(C)) \in U$)
for $\tilde K \subset B \subset A$, $\tilde K \subset
C \subset A$ for each $([\mu ],[\nu ]) \in \bar {\bar L}^2
\cap {\sf U}$, where $\bar {\bar L}:=cl(L,{\sf L}_{R'})$ (or
$cl(L,P_{R'})$), $\pi _{\sf M_o}$ is a projector from
$L$ into $\sf M_o$.}
\par {\bf Proof.} We recall that $\Psi ({\sf R}):=\hat
\Phi (\Theta ({\sf R}))$, where $\Theta ({\sf R})$ is a family of
sequences $(A_n:$ $n \in {\bf N}) \subset \sf R$
for which there exists $\Omega \in \Sigma ({\sf R})$ with
$(n \in \bf N$ $: A_n \in \Omega )$ is infinite,
$\Sigma ({\sf R}):=[(\bigcup
_{n \in J}A_n:$ $J \in {\sf B}({\bf N})):$ $(A_n:$ $n \in {\bf N})\in
\Gamma ({\sf R})]$, $\Gamma ({\sf R})$ is a family of disjoint sequences
$(A_n:$ $n \in {\bf N})\subset \sf R$ for which
$[\bigcup_{n \in J}A_n:$ $J \subset {\bf N}] \subset \sf R]$. A ring of
sets $Z$ is called a quasi-$\delta $-ring, if each disjoint sequence
$(A_n:$ $n \in {\bf N})$ from $\sf R$ the union of which
$\bigcup_nA_n=A$ is contained in a set $B \in \sf R$,
$A \subset B$, has a subsequence $(A_{n_j}:$ $j \in
{\bf N}) \in \Gamma ({\sf R})$.
\par  From Proposition 4.2.7\cite{cons} and the Radonian property
of measures we have $\mu (B) -\nu (C) \in
U'$ for each $(\mu ,\nu ) \in (\bar {\bar {{\sf M \times M}}})
\cap \pi _{\sf M_0}
({\sf U})$. Then for each element $D'$ of the uniformity on $R$
there are $d \in P"$, $c>0$ and $L$ such that $\rho _{\mu }(g,x)-
\rho _{\nu }(h,x') \in D'$ for $g,g' \in W$, $\mu ,\nu  \in M$,
$d(x,x') <c$,
since ${\sf H} \in \hat \Phi_4$, where $\bar {\bar {\sf M}}:=cl({\sf M},
E({\sf R},G;R')_{R'})$ is the closure of ${\sf M} \subset E({\sf R},G;R')$
in $E({\sf R},G;R')_{R'}$, $E({\sf R},G;R'):=[\mu \in E(
{\sf R},G):$ $\mu $ is $R'$-regular $]$, that is, $\mu (F(A,R'))$
converges to $\mu (A)$ for each $A \in R$ (simply regular, if
$R'$ consists of closed subsets and conditions in definition
11.34\cite{hew} are satisfied).
Then $F(A,R')$ is a filter on $R$ generated by a base $[(B \in
{\sf R}:$ $\tilde K \subset B \subset A):$
$\tilde K\subset A, \tilde K \in R']$, $E({\sf R},G)$ is a set of
exhaustive additive maps $\mu $, that is, $\mu (A_n)$
converges to $0$ for each $(A_n:$ $n \in {\bf N})\in \Gamma ({\sf R})$.
\par  Then $\mu ^g(B)-\mu (B) \in U'$, $\nu ^h(C)-
\nu (C) \in U'$ and for $3U' \subset U$ we get 5.4 for $\sf L$.
From Theorems 4.4 and 4.7 (or \S 5.1), the Egorov conditions and the Lebesgue
theorem we get 5.4 for $P$, since $\mu $ are probability measures
and ${\sf L}_{R'}$ (or $P_{R'}$) correspond to uniformity
from \S 5.1.(ii) with $A=R'\times
\tau _e'\times S^c$. Indeed, $\mu ^g(A)-\nu ^h(A)=
(\mu ^g(A)-\mu ^g(V_{U'}))+(\mu ^g(V_{U'})-\nu ^h(V_{U'}))+(\nu ^h(V_{U'})-
\nu ^h(A))$, $\mu ^g(A)=\int_A\rho _{\mu }(g,x)\mu (dx)$ for
each $A \in Bf(S)$, for each $\tau _G\ni U'\ni 0$ there exists a
compact subset $V_U' \subset A$
with $\mu ^g(B)
\in U'$ for each $B \in Bf(A\setminus V_{U'})$ and the same
for $\nu ^h$ (due to the Condition in 5.1 that $R'$ contains $Bf(S)$).
From the separability of $S$, $S'$ and the equality of their
topological weights to $\aleph _0$, restrictions 5.1(a,b) it follows
that there exists a sequence of partitions $Z_n=[(x_m,A_m):$
$m, x_m \in A_m]$
for each $A \in Bf(S)$, $A_i\cap A_j=\emptyset $
for each $i \ne j$, $\bigcup_mA_m=A$, $A_m \in Bf(S)$,
such that $\lim_{n\to \infty }(\mu ^g(A)-\sum_j\rho _{\mu }(g,x_j)\mu
(A_j))=0$ and the same for $\nu $, moreover, for $V_{U'}$ each
$Z_n$ may be chosen finite. Then there exists $W \in \tau '_e$ with
$W\times
(S\setminus V^2)\subset (S\setminus V)$, $\tau _e\ni V \subset V^2$,
$\nu ^g(B)$ and
$\mu ^g(B) \in U'$ for each $B \in Bf(S\setminus V^2)$
(for $G=\bf R$) and $g \in W$ (see 5.1.(a)). Then from $A=[A\cap
(S\setminus V^2)]\cup [A\cap V^2]$ and the existence of compact
$V'_{U'}
\subset V$ with $\mu (E) \in U'$ for each $E \in Bf(V \setminus V'_{U'})$ 
and the same for $\nu $, moreover, $(V'_{U'})^2$ is also
compact,
it follows that $\mu ^g(B)-\nu ^h(C)\in U$ for $9U' \subset U$, since
$R' \supset Bf(S)$, where $W$ satisfies the following condition
$\mu ^g(V'_{U'})-\nu ^h(V'_{U'}) \in U'$ for $V'_{U'} \subset V^2$
due to 5.1.(b), $\mu (B)-\nu (C) \in U'$, $WV'_{U'} \subset (V'_{U'})^2$
due to precompactness of $W$ in $S$. Since pseudo-differentiable
measures are also quasi-invariant, hence for them Proposition
5.4 is true.
\par  Now let $[\mu ] \in \lim H$, $A \in Bf(S)$, then $\eta _{\mu }
\in \lim \zeta (H)$ in $Y(v)$ and there exists a sequence
$\eta _{\mu _n}$ such that $\int_{\bf K}\eta _{\mu _n}(\lambda ,U_*,A)
v(d \lambda )=\tilde D^b_{U_*}\mu _n(A)$ due to \S \S 4.4, 4.7 or 5.1 and
$\lim _{n \to \infty }\tilde D^b_{U_*}\mu _n(A)=$ $\int_{\bf K}\eta _{\mu }
(\lambda ,U_*,A)v(d \lambda )=:\kappa (A)$ due to the Lebesgue theorem.
From $\eta _{\mu }(\lambda ,U_*,A\cup B)=\eta (\lambda ,U_*,A)+
\eta (\lambda ,U_*,B)$ for $A\cap B=\emptyset $, $B \in Bf(S)$
and the Nikodym theorem \cite{cons} it follows that $\nu (A)$ is the
measure on $Bf(S)$, moreover,
$\kappa (A)=\tilde D^b_{U_*}\mu (A)$. Since $\mu ^g(A)=\int_A\rho _{\mu }
(g,x)\mu (dx)$ for $A \in Bf(S)$ for $g\in S'$, then 
$$\eta _{\mu ^g}
(\lambda ,U_*,A)=j(\lambda )g(\lambda ,0,b)[\mu ^g(A)-\mu ^g(U^{-1}
_{\lambda }A)]$$ 
$$=j(\lambda )g(\lambda ,0,b)\int_A\rho _{\mu }
(g,x)[\mu (dx)-\mu ^{U_{\lambda }}(dx)]$$ and in view of the Fubini theorem
there exists 
$$\tilde D^b_{U_*}\mu ^g(A)=\int_A[\int_K\rho _{\mu }
(g,x)j(\lambda )g(\lambda ,0,b)[\mu (dx)-\mu ^{U_{\lambda }}(dx)]v(d
\lambda ),$$ 
where $j(t)=1$ for $S=X$ and $j(t)$ is the characteristic
function of $B({\bf K},0,1)$ for $S$ that is not the Banach space $X$. 
Then $\mu $-a.e.
$\tilde D^b_{U_*}\mu ^g(dx)/ \tilde D^b_{U_*}\mu (dx)$ coincides with
$\rho _{\mu }(g,x)$ due to 5.1.(a,b), hence, $(\tilde D^b_{U_*}
\mu ^g, \rho _{\mu ^g})$ generate the $\Phi _4$-filter in $L$ arising
from the $\hat \Phi _4$-filter in $P$. Then we estimate
$\tilde D^b_{U_*}(\mu ^g-\nu ^h)(A)$ as above $\mu ^g(A)-\nu ^h(A)$.
Therefore, we find for the $\Phi _4$-filter corresponding $L$,
since there exists $\delta >0$
such that $U_{\lambda }\in W$ for each $|\lambda |<\delta $. For
$\Phi _4$-filter we use the corresponding finite intersections
$W_1\cap ...\cap W_n=W$, where $W_j$ correspond to the
$\Phi _4$-filters $H_j$.
\par {\bf 5.5. Corollary.} {\it If $\{ {\sf H'}(A,R') : A \in
{\sf R} \} \subset \Psi _f({\sf R})$, $\sf T$ is a
$\hat \Phi _4$-filter
in $\sf L$, $\sf U$ is an element of uniformity in $\sf L$ (or $P$),
then there are $L \in \sf T$ and an element $\sf V$ of the uniformity
in ${\sf L}_{R'}$ (or $P_{R'}$) such that $\bar {\bar L}^2
\cap {\sf V} \subset \bar {\bar L}^2 \cap \sf U$.}
\par   {\bf Proof.} In view of Corollary 4.2.8\cite{cons} and
Theorem 5.4
for ${\sf U_M}=\pi _{\sf M}\times \pi _{\sf M}\sf U$ there exists
an entourage of the diagonal, $V_{\sf M}$ in $E({\sf R},G;R')_{R'}$,
such that $\bar {\bar M}^2\cap V_{\sf M}\subset \bar {\bar M}^2\cap {\sf U_M}$,
where $\pi _{\sf M}:{\sf L}\to \sf M$ (or $P\to \sf M$) is a projector,
${\sf M(R},G)$ is a set of measures on $\sf R$ with values in $G$,
${\sf M(R},G;R'):=[\mu \in {\sf M(R},G):$ $\mu $ is $R'$-regular $]$.
Since ${\sf M(R},G;R')_{R'} \subset E({\sf R},G;R')_{R'}$, then
there exists an entourage of the diagonal, $V$ in ${\sf L}_{R'}$ (or $P_{R'}$),
such that $\pi _{\sf M}\times \pi _{\sf M}V \subset V_{\sf M}\cap \bar
{\bar M}^2$.
\par {\bf 5.6. Lemma.} {\it Let $\sf R$ be a quasi-$\sigma $-ring
directed by the inclusion, $\sf U$ be an upper directed subset in $\sf R$,
$\sf H$ be a  $\hat \Phi _4$-filter in $\sf L$ (or $P$),
$0\in U \ni \tau _G$. Then there are $M \in \sf H$, $A \in \sf U$,
an element of the uniformity $U'$ in ${\sf L}_
{R'}$ (or $P_{R'}$) such that $\mu ^g(B)-\nu ^h(C) \in U$ (or
in addition $\tilde D^b_{U_*}\mu ^g(B)-\tilde D^b_{U_*}\nu ^h(C)
\in U$) for each $([\mu ],[\nu ])
\in \bar {\bar L}^2 \cap U'$, $g,h \in W$, $W$ is defined by a
projection of $U'$
onto $S'$, $B,C \in \sf U$ with $A \le B$, $A \le C$.}
\par {\bf Proof.} For each $C$ and $E\in \sf U$ by the definition
there are $F \in \sf U$ with $C \le F$ and $E\le F$. There are
open subsets $P$
and $V$ in $S$ and $W \in \tau '_e$ with $WP \subset P^2 \subset V$
such that $(\mu -\nu )(B) \in U$ for each $B \subset S\setminus P$
and $(\mu ,\nu ) \in
\bar {\bar M^2} \cap \pi _{\sf M_o}(U')$. Indeed, $\pi _{\sf M_o}\sf H$
is a base of a filter $\sf T$ in $\sf M_o$. Therefore, $g(B \cap
S\setminus P^2) \subset S\setminus P$ for $g \in W$, consequently,
$[\mu ^g-\nu ^g](B \cap S\setminus P^2) \in U$. In view of \S 5.1.(a,b)
for $D, V$ there are $W$, $c>0$, $d \in P"$ such that $\rho _{\mu }(g,x)
-\rho _{\nu }(h,x') \in D$
for $g, h \in W$ and $d(x,x') <c$, $x, x' \in V_c$, consequently,
$\mu ^g(B)-\nu ^h(C)=
[\mu ^g(B \cap P^2)-\nu ^h(C \cap P^2)]+[\mu ^g(B \setminus P^2)-\nu ^h(
C \setminus P^2)] \in 3(DU+U)$. Modifying the proof of \S 4.2.9 \cite{cons}
we get the statement of this lemma for $\sf L$. For $P$
an estimate of $\tilde D^b_{U_*}
\mu ^g(B)-\tilde D^b_{U_*}\nu ^h(C)$ may be done analogously to
the proof of Proposition 5.4.
\par {\bf 5.7. Theorem.} {\it Let $\sf H$ be a $\hat \Phi _4$-filter
in $\sf L$ (or $P$), $\{A_n:n \in {\bf N} \} \in \Gamma ({\sf R})$,
$\tau _G \ni U \ni 0$. Then there are $L \in {\sf H}$, $M \in
{\sf B}({\bf N})$ and an element $\sf U $ of the uniformity in $\sf L$
such that $\mu ^g(\cup (A_n: n \in M')-\nu ^h(\cup (A_n: n \in M") \in U$
(or in addition $(\tilde D^b_{U_*}\mu ^g(\bigcup_{n \in M'})-
\tilde D^b_{U_*}\nu ^h(\bigcup_{n \in M"}A_n)) \in U$
for each $([\mu ],[\nu ]) \in \bar L^2 \cap {\sf U}$
and $M', M" \in {\sf B}({\bf N})$ with $M \cap M'=M\cap M"$. If
$\{ {\sf H'}(A,R') \} \subset \Psi _f({\sf R})$, then $\sf U$ may be
chosen as an element of the uniformity in ${\sf L}_{R'}$ (or $P_{R'}$),
where $\bar L:=cl(L,{\sf L})$ (or $\bar L:=cl(L,P)$).}
\par {\bf Proof.} We recall that ${\sf B}_f({\bf N})$ denotes the
family of finite subsets in $\bf N$, $\Psi _f({\sf R}):=\hat \Phi
(\Theta _f({\sf R}))$, $\Theta _f({\sf R})$ is a family of sequences
$(A_n:$ $n \in {\bf N}) \subset \sf R$ for which there exists
$\Omega \in \Sigma _f({\sf R})$ with $card[n \in {\bf N}:$
$A_n \in \Omega ]=\aleph _0$,
$\Sigma _f({\sf R}):=[(\bigcup_{n \in M}A_n:$ $M \in {\sf B}_f({\bf N})):$
$ (A_n:$ $n \in {\bf N}) \in \Gamma ({\sf R})]$.
Let $\phi (M):=\cup \{ A_n: n \in M\}$,
$M \subset {\sf B}_f({\bf N})$, $\Sigma $ be a filter of neighbourhoods
of $\emptyset $
in ${\sf B}_f({\bf N})$, then $\phi ( \Sigma )\in \Psi _f({\sf R})$. In
view of \S 4.1.14 \cite{cons} we have $\mu (\phi (\Sigma )) \to 0$ for each
$\mu \in {\sf M(R},G,R')$.
We take a symmetric neighbourhood $V' \ni 0$ in $G$ with $3V' \subset U$.
Then there are  $L \in \sf H$ and $M \in {\sf B}_f({\bf N})$ such that
$\mu ^g(\phi (M'))-\nu ^h(\phi (M")) \in V'$ (or in addition
$\tilde D^b_{U_*}\mu ^g(\phi (M'))-\tilde D^b_{U_*}\nu ^h(\phi (M"))
\in V'$) for each $([\mu ],[\nu ])
\in  \bar L^2 \cap \sf U$ and $M', M" \in {\sf B}({\bf N}
\setminus M)$ (see Lemma 5.6). We choose $\sf U$ such that
$\mu ^g(\phi (M_o))-\nu ^h(\phi (M_o))
\in V'$ (or in addition $\tilde D^b_{U_*}\mu ^g(\phi (M_o))-
\tilde D^b_{U_*}\nu ^h(\phi (M_o)) \in V'$)
for each $([\mu ],[\nu ]) \in \sf U$ and $M_o
\subset M$, consequently, $\mu ^g(\phi (M'))-\nu ^h(\phi (M")) \in
3V' \subset U$ (or additionally $\tilde D^b_{U_*}\mu ^g(\phi (M'))-
\tilde D^b_{U_*}\nu ^h(\phi (M")) \in 3V'$)
for $M \cap M'=M \cap M"$, the last statement follows from Corollary
5.5.
\par {\bf 5.8. Corollary.} {\it Let $(A_n: n \in \bf N) \in \Gamma
({\sf R})$, $\sf H$ be a $\hat \Phi _4$-filter in $\sf L$ (or $P$),
$\tau _G
\ni U \ni 0$, then there are $(L,M) \in {\sf H \times B}_f({\bf N})$
such that $\mu ^g(\cup (A_n: n \in M')) \in U$ (or additionally
$\tilde D^b_{U_*}\mu ^g(\bigcup_{n \in M'}A_n) \in U$)
for $[\mu ]\in
\bar L$ and $M' \in {\sf B}({\bf N}\setminus M)$.}
\par   {\bf Proof.} Let us take a fixed $\nu $ and $W$
from \S 5.7 with $\nu (\bigcup_{n \in M'}A_n)\in \tilde U$ (or additionally
$\tilde D^b_{U_*}\nu (\bigcup_{n \in M'}A_n) \in \tilde U$) and
$(\mu ^g-\nu )(\bigcup_{n \in M'}A_n) \in \tilde U$), where $2\tilde U
\subset U$, $0 \in \tilde U \in \tau _G$, $g \in W$.
\par {\bf 5.9. Theorem.} {\it Let ${\sf U} \subset {\sf R}\times \tau _e'
\times S^c=:Z$ be such that $id:{\sf L_U} \to {\sf L}$ (or $id:P_{\sf U}
\to P$) be uniformly
$\Phi _4$-continuous, $\{ {\sf H'} (A,R'): A \in {\sf R} \}
\subset \Psi _f({\sf R})$. Then 
\par (a) $\sf L$ (or $P$) is
$\Phi _4$-closed in
$G^{\sf R} \times R_o^{S' \times S}$ (or $G^{\sf R}\times
R_o^{S'\times S}\times G^{S'\times K\times \sf R}$),
$\sf M_o$ is  $\Phi _4$-closed in $G^{\sf R}$;
\par (b) if $G$ and $R$ are $\Phi _i$-compact, then $\sf L_U$
(or $P_{\sf U}$) is $\Phi
_i$-compact, where $i \in (1,2,3,4)$;
\par (c) if $G\times R_o$ is sequentially complete,
then $\sf L_U$ (or $P_{\sf U}$) is also; 
\par (d) if $(0,0)$ is the
$G_{\delta }$-subset in $G\times R$, then $\sf L_U$ (
or $P_{\sf U}$) is $\Phi _2$-compact; if additionally
$G\times R$ is sequentially complete, then $\sf L_U$ (or
$P_{\sf U}$) is $\Phi _4$-compact;
\par (e) $\sf L_U$ (or $P_{\sf U}$) is Hausdorff.}
\par {\bf Proof.} We recall that a subset $A \subset E$
of a topological or of a uniform space $E$ is called
$\Phi _i$-closed (or compact), if for each $\Phi _i(A)$-filter
$F$ is satisfied $\lim F \subset A$ (or $\lim F \ne \emptyset $
respectively, that is, this definition of compactness differs
from the usual). In view of Theorem 4.2.14\cite{cons} and Proposition
5.4 is accomplished 5.9(a). 
\par From \S \S 2.1.14,
1.8.11\cite{cons} and Theorem 5.7, $\Phi _i$-compactness of $\sf M$
and $\sf M_o$,
also from the completeness of $R_o$, $L^1(K,v,{\bf C})$ it
follows (b). 
\par Then (c) follows from (a) and \S 1.8.7\cite{cons};
\par (d) follows from (c) and \S 1.6.4\cite{cons}, since $G\times R$ is
$\Phi _2$-compact.
From ${\sf M(R},G;R')_{\sf U}$ (see \S 4.2.14\cite{cons}),
$R_o
^{S'\times S}$ and $L^1(K,v,{\bf C})$ being
Hausdorff it follows (d).
\par {\bf 5.10. Theorem.} {\it Let $\Phi $ be a set of $\hat
\Phi _4$-filters on $\sf L$ (or $P$), also let $\phi : {\sf L \times R}
\times S' \times S (\vee \times {\bf K})\to
G \times R_o (\vee \mbox{ }\times Y(v))$  be such that ($\phi
(\mu ,\rho _{\mu }(*,*)(\vee \mbox{  }\eta _{\mu }),A,g,x):=
(\mu (A),\rho _{\mu }(g,x)(\vee \mbox{ }\eta _{\mu ^g}(t,U_*,A))$). 
If $\{ {\sf H'}(A,R')| A \in {\sf R} \} \subset \Psi _f({\sf R})$,
then a mapping $\phi :{\sf R} \times S' \times S (\vee \times {\bf K})\to
G^{{\sf M_o(R)},G;R'}_{\Phi '} \times
R_o^{S' \times S}(\vee \mbox{ }\times Y(v))$ gives $R'$-regular quasi-invariant measure
(or in addition a pseudo-differentiable measure)
$(A,g,x) \to \phi (*,A,g,x)$ (or $(A,g,x,t) \to \phi (*,A,g,x,t)$),
satisfying Conditions 5.1.(a,b),
where $\Phi ':= \pi _{\sf M_o} (\Phi )$, in $R_o^{S' \times S}$ the
topology corresponds to the topology in $\sf L$.}
\par   {\bf Proof} follows from \S \S 5.4 and 5.8.
\par  {\bf 5.11. Note.} Let it be a sequence
$\{[\mu _n] : n \in {\bf N} \}$ of quasi-invariant measures (or 
pseudo-differentiable measures)
converging uniformly in uniformity 5.1.(ii) and fulfilling Conditions
5.1.(a,b)
then in accordance with Corollary 5.8 $\{ (\mu _n)^g: n \in {\bf N} \} $
(or also $[\tilde D^b_{U_*}\mu ^g_n: n]$) is uniformly $\sigma $-additive
for each fixed $g \in S'$. Moreover, it is uniformly $\sigma $-additive by
$g \in W$ for each given $B \in
\sf R$ such that $gB \subset V$ for suitable open $W$ in $S'$
and $V$ in $S$.
For $\sf L$ this means that for each $0\in D \in \tau _R$ and 
$e\in U \in \tau _G$ there are $W \ni e$,
$d \in P"$, $c>0$, $n$ and $V \ni e$, a compact subset $V_U$,
$e \in V_U\subset V$
with $\mu _m(C) \in U$ (or $\tilde D^b_{U_*}\mu _m(C) \in U$ in
addition for $P$)
for $C\in {\sf R}$ and  $C \subset S\setminus V_U$, with
$\rho _{\mu _m}(g,x)-\rho _{\mu _j}(h,x')
\in D$ (or in addition $\tilde D^b_{U_*}(\mu ^g_m-\mu ^h_j)(A) \in U$
for $A \in Bf(V_U)$ for $P$) whilst $g,h \in W$
for each $x,x'\in V_U$ with $d(x,x') <c$  and $m,j > n$.  In view
of Theorem 5.9
there exists $\lim _{n \to \infty }
(\mu _n,\rho _{\mu _n}(g,x))= (y,d(y;g,x)\in \sf L$,
that is, a quasi-invariant measure (or 
pseudo-differentiable measure $\lim_n [\mu ]_n =[\mu ]\in P$).
Therefore, they are analogs of the Radon theorem for quasi-invariant 
measure and pseudo-differentaible measures.
\section {Appendix.}
Suppose $X=c_0(\omega _0,{\bf K})$ is a Banach space over 
a locally compact non-Archimedean infinite field $\bf K$
with a non-trivial valuation and  $I$ is a unit operator on $X$.
If $A$ is an operator on $X$, then in some basis of $X$ we have
an infinite matrix $(A_{i,j})_{i,j \in \bf N}$, so we can consider
its transposed matrix $A^t$. If in some basis the following equality 
is satisfied $A^t=A$, then $A$ is called symmetric.
\par {\bf A.1. Lemma.} {\it Let $A: X\to X$
be a linear invertible operator with a compact operator $(A-I)$.
Then there exist an orthonormal basis
$(e_j:$ $j\in {\bf N})$ in $X$, invertible linear
operators $C, E, D: X\to X$ with compact opeartors $(C-I)$, $(E-I)$, $(D-I)$
such that $A=SCDE$, $D$ is diagonal, $C$ is lower triangular and $E$ is
upper triangular, $S$ is an operator transposing a 
finite number of vectors from an orthonormal basis 
in $X$. Moreover, there exists $n \in \bf N$ and invertible
linear operators $A', A": X\to X$ with compact operators $(A'-I)$,
$(A"-I)$ and $(A'_{i,j}-\delta _{i,j}=0)$
for $i$ or $j>n$, $A"$ is an isometry and there exist their determinants
$det(A')det(A")=det(A)$, $|det(A")|_{\bf K}=1$, $det(D)=det(A)$.
If in addition $A$ is symmetric, then $C^t=E$ and $S=I$.}
\par  {\bf Proof.} In view of Lemma 2.2\cite{sch2} for each $c>0$
there exists the following decomposition $X=Y\oplus Z$ into
$\bf K$-linear spaces such that $\| (A-I)|_Z \| <c$, where
$dim_{\bf K}Y=m< \aleph _0$. In the orthonormal basis 
$(e_j:$
$j)$ for which $sp_{\bf K}(e_1,...,e_m)=Y$ for $c \le 1/p$ we get
$A=A'A"$ with $(A-I)|_Z=0$, $|A"_{i,j}-\delta _{i,j}| \le c$ for each
$i, j$ such that $(A'_{i,j}-\delta _{i,j})=0$ for $i$ or $j>n$, where
$n \ge m$
is chosen such that $|A_{i,j}-\delta _{i,j}|\le c^2$ for $i>n$ and
$j=1,...,m$,
$A_{i,j}:=e_i^*(Ae_j)$, $e_i^*$ are vectors $e_i$
considered as linear continuous functionals $e_i^*
\in X^*$. Indeed, $(A_{i,j}:$  $i \in {\bf N})=Ae_j\in X$
and $\lim_{i \to \infty }A_{i,j}=0$ for each $j$. From the form of
$A"$ it follows that $\| A"e_j-e_j\| \le 1/p$ for each $j$,
consequently, $\| A"x\|=\| x\|$ for each $x \in X$. Since
$A"=(A')^{-1}A$, $(A-I)$ and $(A'-I)$ being compact, hence $(A"-I)$
is compact together with $(A^{-1}-I)$, $((A')^{-1}-I)$ and
$((A")^{-1}-I)$. Moreover, there exists $\lim_{k \to \infty }
det(A)_k =det(A)$ $=\lim_k
det((A')_k(A")_k)$ $=\lim_kdet(A')_kdet(A")_k$ $=det(A')det(A")$,
where $(A)_k:=(A_{i,j}:$ $i,j \le k)$. This follows from the
decompositions
$X=Y_k\oplus Z_k$ for $c=c(k) \to 0$ whilst $k \to \infty $. This
means that for each
$c(k)=p^{-k}$ there exists $n(k)$ such that $|A_{i,j}-
\delta _{i,j} |<c(k)$, $|A'_{i,j}-\delta _{i,j}|<c(k)$ and $|A"_{i,j}-
\delta _{i,j}|<c(k)$ for each $i$ or $j >n(k)$, consequently,
$|A{{1...n(k)i_1...i_q}\choose {1...n(k)j_1...j_q}} -A{{1...n(k)}
\choose {1...n(k)}}\delta _{i_1,j_1}...\delta _{i_q,j_q}| <c(k)$,
where $A{{i_1...i_r}\choose {j_1...j_r}}$ is a minor corresponding to
rows $i_1,...,i_r$ and columns $j_1,...,j_r$ for $r, q \in \bf N$. From
the ultrametric inequality it follows that $|det(A")-1| \le 1/p$,
hence $|det(A")|_{\bf K}=1$, $det(A")_k \ne 0$ for each $k$, $det(A')_k
=det(A')_n$ for each $k \ge n$. Using the decomposition of
$det(A')_n$
by the last row (analogously by the column) we get
$A'_{n,j} \ne 0$ and a minor
$A'{{1...n-1}\choose {1...j-1,j+1,n}} \ne 0$. Permuting the columns
$j$ and $n$ (or rows) we get as a result a matrix $(\bar A')_n$ with $\bar
A'{{1...n-1}\choose {1...n-1}} \ne 0$. Therefore, by the
denumeration of the basic vectors we get $A'{{1...k}\choose {1...k}}
\ne 0$ for each $k=1,...,n$, since $|det(\bar A')_n|=|det(A')_n|$.
\par   Therefore, there exists the orthonormal basis 
$(e_j:$ $j)$ such that $A{{1...j}\choose
{1...j}} \ne 0$ for each $j$ and $\lim_j A{{1...j}\choose {1..j}}
=det(A)\ne 0$. Applying to $(A)_j$ the Gaussian decomposition and using
compactness of $A-I$ due to formula (44) in \S II.4\cite{gan}, which is
valid in the case of $\bf K$ also, we get $D=diag(D_j:$
$j\in {\bf N})$, $D_j=A{{1...j}\choose {1...j}}/A{{1...j-1}\choose
{1...j-1}}$; $C_{g,k}=A{{1,...,k-1,g}\choose {1,...,k-1,k}}/A{{1...k}
\choose {1...k}}$; $E_{k,g}=A{{1,...,k-1,k}\choose {1,...,k-1,g}}/
A{{1...k}\choose {1...k}}$ for $g=k+1, k+2,...$, $k \in \bf N$.
Therefore, $(C-I)$, $(D-I)$, $(E-I)$ are the compact operators,
$C_{i,j}$, $D_j$,
$E_{i,j} \in \bf K$ for each $i, j$. Particularly, for $A^t=A$ ($A^t$
denotes the transposed matrix for $A$) we get $E_{k,g}=C_{g,k}$.
\par {\bf A.2. Notes. 1.} An isometry operator $S$ does not influence
on the results given above by the measures, since $A|Z$ and $S^{-1}A$ on $X$
have decompositions of the form $CDE$, where $X\ominus Z$ is
a finite-dimensional subspace of $X$ over $\bf K$, hence $X\ominus Z$
is locally compact for the locally compact field $\bf K$. 
Therefore, on $X\ominus Z$
there exists the Haar measure 
with values in $\bf R$ or $\bf Q_s$ with $s\ne p$,
$\bf K\supset Q_p$.
\par {\bf A.2.2.} For compact $A-I$ we can construct the following
decomposition $A=BDB^tC$, where $BB^t=I$, $CC^t=I$, $D$ is diagonal,
$B$, $D$ and $E$
are operators on $c_0(\omega _0, {\bf C_p})$ with $\lim_{i_j \to \infty }
(B_{i,j}-\delta _{i,j})=\lim_j(D_j-1)$ $\lim_{i+j \to \infty }(C_{i,j}
-\delta _{i,j})=0$. But, in general, matrix elements of $B, C, D$ are
in $\bf C_p$ and may be in $\bf C_p\setminus \bf K$, since from
the seqular equation $det(A-\lambda I)$ even for symmetric
matrix $A$ over $\bf K$ in general may appear $p^{1/n}$ for $n \in
\bf N$, but $\bf R$ is not contained in $\bf K$ and $\bf C_p$
is not locally compact. This decomposition is not used
in the present article for the construction of quasi-invariant measures
on a Banach space $X$
over $\bf K$, since on $\bf C_p$ there is not any non-trivial invariant
measure (or even quasi-invariant relative to all shifts from $\bf C_p$)
and may be $B_{i,j}$, $D_j$, $C_{i,j} \in \bf C_p\setminus \bf K$.
Instead of it (and apart from \cite{sko}) we use the decomposition
given in the Lemma A.1.
\par {\bf A.3. Remarks.} Let $\bf K$ be a  
non-Archimedean infinite field with a non-trivial valuation.
Let $X$ and $Y$ be normed spaces
over $\bf K$, and $F: U\to Y$ be a function, where
$U\subset X$ is an open subset.
The function $F$ is called differentiable, if for each
$t \in \bf K$, $x \in U$ and $h \in X$ such that $x+th \in U$ there exists
$DF(x,h):=\{ dF(x+th)/dt \mid t=0 \} =\lim_{t\to 0, t\ne 0} 
\{ F(x+th)-F(x) \}/t$  and
$DF(x,h)$ linear by $h$, that is, 
$DF(x,h)=:F'(x)h$, where $F'(x)$ is a bounded linear operator
(derivative). Let  $\Phi ^1 F(x;h;t):=\{ F(x+th)-
F(x)\} /t$ for each $t\ne 0$, $x\in U$, $x+th\in U$, $h\in X$. 
If this function $\Phi ^1 F(x;h;t)$ has a continous extension
$\bar \Phi ^1 F$ on $U\times V\times S$,
where $U$ and $V$ are open neighbourhoods of $x$ and $0$ in $X$, 
$S=B({\bf K},0,1)$ and
$\| {\bar \Phi }^1 F(x;h;t) \| :=\sup \{ \| {\bar \Phi }^1 F(x;h;t) \| /
\| h\|: x \in U, 0\ne h  \in V, t
\in S \}< \infty $  and ${\bar \Phi }^1 F(x;h;0)=F'(x)h$, then 
$F$ is called continuously differentiable on $U$, 
the space of all such $F$ is denoted by
$C^1(U,Y)$, where $B(X,y,r):= \{ z \in X: \| z-y \|_X \le r \}$.
By induction we define $\Phi ^{n+1} F(x;h(1),...,
h(n+1);t(1),...,t(n+1) \}:= \{ \Phi ^n F(x+t(n+1)h(n+1);h(1),..,h(n);
t(1),...,t(n))- \Phi ^n F(x;h(1),...,h(n);t(1),...,t(n)) \}/t(n+1)$
and the space $C^n(U,Y)$ \cite{lu3,luseabm}.
\par The family of all bijective surjective mappings of
$U$ onto itself of class $C^n$ is called the diffeomorphism group
and it is denoted by $Diff^n(U):=C^n(U,U)\cap Hom(U).$
\par {\bf A.4. Theorem.} {\it Let $\bf K$ be a spherically complete 
non-Archimedean infinite field with non-trivial valuation,
$F \in C^n(U,Y)$ with $n\in {\bf N}:=\{ 1,2,.. \} $, 
$F: X \to Y$, $X$ and $Y$ be Banach spaces, $U$ be an open neighbourhood of
$y \in X$, $F(X)=Y$ and $[F'(y)]^{-1}$ be continuous, where $n\in \bf Z$,
$0\le n<p$ when $char ({\bf K})=:p>0$, $n\ge 0$ when $char ({\bf K})=0$. 
Then there exists $r>0$ and a locally inverse operator
$G=F^{-1}$:
$V:=B(Y,F(y),\| F'(y) \| r)\to B(X,y,r)$, moreover, $G \in C^n(V,X)$ and
$G'(F(y))=[F'(y)]^{-1}$.}
\par {\bf Proof.} Let us consider an open subset $U \subset X$,
$[y,x] :=[z:$ $z=y+t(x-y),$ $|t| \le 1,$ $t \in {\bf K}] \subset
U$ and $g: Y \to \bf K$ be a continous $\bf K$-linear functional, 
$f(t):=g(F(y+t(x-y))$. Such family of non-trivial $g$ separating points of 
$X$ exists due to the Hahn-Banch theorem 
\cite{nari,roo}, since $\bf K$ is spherically complete.
Then there exists $f'(t)=g(F'(y+t(x-y))(x-y)),$
$f \in C^n(B({\bf K},0,1),{\bf K})$, consequently, for $F'(z) \ne 0$,
$z \in [x,y]$ there exists $b>0$ such that for $g(F'(z)) \ne 0$
we have $g(F(x+t(x-y)))-g(F(x))|=|g(F'(z)(x-y))|\times |t|$
for each $|t|<b$. From $|g(u(z))-g(u(x))|=|g(z)-g(x)|$ and the
Hahn-Banach it follows, that $u$ is the local isometry, where
$u(z):=[F'(y)]^{-1}(F(z))$, $z$, $x \in B(X,y,r)$, $r$ is chosen
such that $\| (\bar \Phi ^1u(x,h,t))-h \| < \| h \| /2$ for each
$x \in B(X,y,r),$ $h \in B(X,0,r)$, $t \in B({\bf K},0,1)$. 
Therefore, $F(B(X,y,r)) \subset B(Y,F(y),\| F'(y)\| r)=:S$. 
Applying to the function $H(z):=z-[F'(y)]^{-1}(F(z)-q)$ for $q \in S$ 
the fixed point theorem we get $z'$ such that 
$F'(z')=q$, since $\| H(x)-H(z)\| <b
\times \| x-z\| .$ From the uniform continuos by $(g,h)$
differentiability $f_{g,h}(t):=g(F(y+th)$ by $t$ it follows
continous differentiability of $G$.
\par  {\bf A.5. Theorem.} {\it Let $f \in C^n(U,Y)$, where $n \in \bf
N$, $Y$ is a Banach space, $n\in \bf Z$, $0\le n<p$ when
$p:=char ({\bf K})>0$, $n\ge 0$ when $char ({\bf K})=0$. 
Then for each $x$ and $y\in U$ is accomplished the following formula:
$f(x)=f(y)+\sum_{j=1}^{n-1}f^{(j)}(y)(x-y)^j/j!+ R_n(x,y)(x-y)^{n-1}$, where
$f^{(j)}(y)(x-y)^j=(\bar \Phi ^jf)$ $(y; x-y,...,x-y;0,...,0)\times j!,$
$f^{(j)}(y)h^j:=$ $f^{(j)}(y)$ $(h,...,h)$, $f^{(j)}(y): U\to 
L_j(X^{\otimes j},Y)$,
$R_n(x,y): U^2\to L_{n-1}(X^{\otimes (n-1)},Y)$ 
with $R_n(x,y)=o(\| x-y\|)$, where $L_j(X^{\otimes j},Y)$ is
the Banach space of continuous polylinear operators 
from $X^{\otimes j}$ to $Y$, $U$ is open in $X$.}
\par  {\bf Proof.} For $n=1$ this formula follows from the definition
of $\bar \Phi ^1f$. For $n=2$ let us take $R_2(x,y)=\bar \Phi
^2f(x;y,y)$. Evidently, $C^{n-1}(U,Y) \supset
C^n(U,Y)$. Let the statement be true for $n-1$,
where $n \ge 3$. Then from $\bar
\Phi ^{n-1}f(y+t_n(x-y);x-y,...,x-y; t_1,...,t_{n-1})=$ $\bar \Phi
^{n-1}f(y;x-y,...,x-y;t_1,..., t_{n-1})+$ $(\bar \Phi ^1(\bar \Phi
^{n-1}f(y;x-y,...,x-y;t_1,..., t_{n-1}))(y;x-y;t_n))\times t_n$ 
and continuity of
$\bar \Phi ^n$ it follows that $R_n(x,y)(x-y)^{n-1}=\bar \Phi
^nf(y;x-y,...,x-y;0,...,0,1)-$ $f^{(n)}(y)(x-y)^n/n!=$ $\alpha
(x-y)(x-y)^{n-1}$, where $\lim_{x\to y, \mbox{ }x\ne y}\alpha
(x-y)/\|x-y\|=0$, that is, $\alpha (x-y)= o(\|x-y\|)$.

\par Address: Theoretical Department,
\par Institute of General Physics,
\par Russian Academy of Sciences,
\par Str. Vavilov 38, Moscow, 119991 GSP-1, Russia
\end{document}